\documentclass[12pt,reqno]{amsart}
\usepackage{amssymb}
\usepackage{amsmath}
\usepackage{bm}
\usepackage{verbatim}
\makeatletter
\@addtoreset{equation}{section}
\makeatother

\newtheorem{thm}{Theorem}[section]
\newtheorem{lem}[thm]{Lemma}
\newtheorem{prop}[thm]{Proposition}
\newtheorem{cor}[thm]{Corollary}

\newtheorem{rem}[thm]{Remark}

\newcommand{\mc}[1]{{\mathcal #1}}

\newcommand{\mb}[1]{{\mathbf #1}}
\newcommand{\bb}[1]{{\mathbb #1}}
\newcommand{\reff}[1]{(\ref{#1})}
\newcommand{\lf}[1]{\left #1}
\newcommand{\rg}[1]{\right #1}

\title[Logarithmic Sobolev for Inhomogeneous Zero Range]{Logarithmic Sobolev Inequality for the
Inhomogeneous Zero Range Process }
\author{Hanna K. Jankowski}

\address{\noindent Department of Statistics, University of
Washington
\newline
e-mail:  \rm \texttt{hanna@stat.washington.edu} }

\date{June 24, 2006}

\begin{document}

\maketitle

\begin{abstract}
We prove that the logarithmic Sobolev constant for the
inhomogeneous symmetric nearest neighbour zero range process on a
cube of size $N^d$ grows as $N^2$.  We apply this result to the
inhomogeneous process which arises in the study of the homogeneous
version of the zero range interacting particle system with
colours.
\end{abstract}

\section{Introduction}

The logarithmic Sobolev inequality is a spectral bound which
provides much information about decay to equilibrium of the
dynamics of a stochastic process.

Consider a process governed by reversible dynamics described by a
generator $\mathcal L$, with semi-group $P_t$ and an invariant
measure $\mu$. The Dirichlet form is defined as $D_\mu(f)= \mu[f
(- \mathcal L)f]$. A logarithmic Sobolev inequality is a statement
which says that the entropy, $H(f|\, \mu)=\mu[f\log f]$, is
bounded by a constant times the Dirichlet form
\begin{eqnarray}\label{ls_genform}
H(f|\,\mu)\leq  C_{\mathrm{LS}} \ D_\mu(\sqrt{f}),
\end{eqnarray}
for all densities $f$. Additionally, the logarithmic Sobolev
inequality implies exponential decay of both the $L^2(\mu)$ norm
and entropy of $P_tf$.

 The Poincar\'{e} inequality is defined as
the bound, uniform in $f$,
\begin{eqnarray}\label{sg_genform}
\mbox{Var}_\mu[f] \leq  C_{\mathrm{SG}} \ D_\mu(f),
\end{eqnarray}
where $\mbox{Var}_\mu[f]$ is the variance of the function $f$ with
respect to the measure $\mu$. The inequality states that the
dynamics of the process have a spectral gap of order
$C_\mathrm{SG}^{-1} $ and hence there is exponential decay to
equilibrium in the $L^2$ sense. That is, we have that
$$\mbox{Var}_\mu[P_tf]\leq e^{-2t/C_\mathrm{SG}}\mbox{Var}_\mu[f],$$
where $P_t$ is the semi-group of the process.

An intermediate spectral bound may be established via the
following entropy dissipation inequality,
\begin{eqnarray}\label{ed_genform}
H(f|\,\mu)\leq  C_{\mathrm{ED}} \ \mu[f (-\mathcal L) \log f],
\end{eqnarray}
uniformly in positive functions $f$.   Because
$\partial_tH(P_tf|\,\mu)= \mu[P_tf (-\mathcal L) \log P_tf]$, this
implies that
$$H_\mu[P_tf]\leq e^{-2t/C_\mathrm{ED}}H_\mu[f],$$
again for positive functions.  One can also show that,
\begin{eqnarray}\label{relat_genform}
C_\mathrm{SG}\leq 2 C_\mathrm{ED}\leq \frac{1}{2}C_\mathrm{LS},
\end{eqnarray}
establishing a relationship between the three inequalities
\cite{persi}.

In the study of hydrodynamic scaling limits of interacting
particle systems, an understanding of the decay to equilibrium of
the dynamics is an important ingredient.  See for example
\cite{kl} for a review of the available methods. If the spectral
gap is of the order $N^{-2}$, we can establish a hydrodynamic
scaling limit for the process. However, when we consider the
fluctuations of this result in \emph{nonequilibrium}, the central
limit theorem multiplicative re-scaling by $\sqrt{N}$ requires
stronger tools. Chang and Yau \cite{yau92} developed a method to
prove nonequilibrium density fluctuations for the Ginzburg-Landau
model, which makes use of the logarithmic Sobolev inequality.
Indeed, our main interest in the inhomogeneous inequality stemmed
from the study of the nonequilibrium fluctuations of the
hydrodynamic scaling limit for the colour version of the zero
range process.

Spectral bounds for the zero range process have been studied
extensively in the literature.  The spectral gap of order $N^{-2}$
was established by Landim, Sethuraman and Varadhan for the
homogeneous symmetric nearest neighbour zero range process
\cite{lsv}.  They make the usual assumptions on the jump rate
function of the zero range process, $c(\cdot)$. Namely, they
assume Lipschitz growth of the rate function
\begin{eqnarray}\label{LGhom}
\sup_k|c(k+1)-c(k)| <\infty
\end{eqnarray}
as well as a weak monotonicity condition:
\begin{eqnarray}\label{Mhom}
\inf_k \{c(k+k_0)-c(k)\} >0,
\end{eqnarray}
for some integer $k_0$.  Assumption \reff{LGhom} is necessary to
ensure that the zero range process is well defined on the infinite
lattice \cite{a}.  Condition \reff{Mhom} rules out the cases, such
as the queueing system corresponding to $c(k)=\mathbb I (k\geq
1)$, where $C_\mathrm{SG}$ depends on the density of particles.

These are also the assumptions under which Dai Pra and Posta
showed the logarithmic Sobolev inequality in \cite{logsob1,
logsob2}. There they show that $C_\mathrm{LS}= CN^{2}$, where $C$
is independent of the particle density. Their approach is based on
the martingale method of Lu and Yau \cite{ly}.

Recently, Caputo and Posta \cite{cp} studied the case of the
inhomogeneous zero range process on the complete graph. As before,
allow the system to evolve on a cube of size $N^d$. The complete
graph setting means that particles are allowed to jump to any
other location of the cube with equal probability.  In the nearest
neighbour case, also known as local dynamics, particles make jumps
to one of their nearest neighbours.  Inhomogeneity means the the
rate at which the first particle leaves site $x$ depends on $x$,
and hence we now consider a family of rate functions $c_x(\cdot)$.

For the complete graph dynamics it is shown in \cite{cp} that
under the condition
\begin{eqnarray*}
\inf_{x,\,k} \{c_x(k+1)-c_x(k)\} > 0
\end{eqnarray*}
on the rate functions, the system has a spectral gap of constant
order.  This is known to imply a spectral gap of order $N^{-2}$
for the nearest neighbour model \cite{q}.  Under the additional
assumption
\begin{eqnarray*}
\sup_{x,\,k} \{c_x(k+1)-c_x(k)\} < \infty,
\end{eqnarray*}
Caputo and Posta also prove the entropy dissipation inequality
\reff{ed_genform} for the complete model.

In this article we consider the symmetric nearest neighbour
inhomogeneous zero range process.  Under these dynamics, the
particles move around a cubic subset of $\bb Z^d$ of size $N^d$.
Particles wait exponential time to make a jump, and then jump to
one of their closest neighbours with equal probability.   We show
that in this case the logarithmic Sobolev constant $C_\mathrm{LS}$
behaves like $C N^2$. It seems natural to study the problem under
the uniform versions of conditions \reff{LGhom} and \reff{Mhom}
\begin{eqnarray*}
(LG) & \sup_{k,\,x}|c_x(k+1)-c_x(k)|\leq  a_1 <\infty \\
(M)  & \inf_{k,\,x} \{c_x(k+k_0)-c_x(k)\} \geq a_2 > 0,
\end{eqnarray*}
for some fixed constants $a_1$, $a_2$ and integer $k_0$.  However,
for technical reasons, at this time we also need to make the
additional condition we now describe.

In the homogeneous case the grand canonical measures are product
measures with marginals indexed by the \emph{constant} density
$\rho$. When we move to the inhomogeneous case this is no longer
true;  the marginals are not spatially homogeneous. However, we
may now index the measure by the \emph{overall} density, which we
define as the average of the local densities $\rho_x$. Let
$\mu_{\Lambda,\rho}$ denote the grand canonical measure for the
process on the box $\Lambda$ with overall density $\rho$, and the
zero range particle configurations by $\eta$. We assume that
\begin{eqnarray}\label{extracond}
0<\inf_r \sqrt{r}\, \mu_{\Lambda, \frac{r}{|\Lambda|} }
\left(\sum_{x\in \Lambda}\eta(x)=r\right) \leq\sup_r \sqrt{r}\,
\mu_{\Lambda, \frac{r}{|\Lambda|} } \left(\sum_{x\in
\Lambda}\eta(x)=r\right) < \infty
\end{eqnarray}
for any size $|\Lambda|\geq 2$.  This condition is required only
in the proof of Lemma \ref{boundfromhell}.  One simple case when
this condition is satisfied occurs when we assume that there
exists a positive constant $\theta$ and universal $K_0$ such that
\begin{eqnarray}\label{egLLR}
c_x(k)=\theta k, \,\,\,\,\,\,\, \forall \,\, x \,\, \mbox{and}
\,\, k\geq K_0.
\end{eqnarray}
See Remark \ref{condLLR} for more details.

The inhomogeneous zero range process arises naturally in the
following setting.  Consider first the homogeneous case and assign
one of $k$ colours to each of the particles.  The configurations
of particles of each colour, considered jointly, form a Markov
process with a family of invariant measures.  Next, we single out
one of the colours, and condition on the configuration of the
remaining particles.  The invariant measures are still product
measures which can be seen as a special case of the invariant
measures for the inhomogenous model.  One can easily show that
assumptions \reff{LGhom} and \reff{Mhom} on the rate function
$c(\cdot)$ in the original homogeneous process imply conditions
(LG) and (M) for the induced inhomogeneous rates.  The additional
condition can be attained if we assume, for example, that
$c(k)=\theta k$ for all $k$ sufficiently large, as this implies
\reff{egLLR}, which in turn implies \reff{extracond}.

The relationship of the inhomogeneous process with the colour
homogeneous version was our main interest in writing down this
result. Notice that the system evolution for one colour in the
colour version of the process is \emph{not} the same as the
evolution of one colour conditioning on, or fixing, the remaining
particles. However, because the logarithmic Sobolev inequality is
a \emph{static} property, the inhomogeneous setting provides
useful spectral bounds regardless. Using this idea, we were able
to establish nonequilibrium fluctuations for the colour zero range
process \cite{me}.

The proof of \reff{ls_genform} for the inhomogenous zero range
process, with $C_\mathrm{SG} = CN^2$,   is the main result of this
paper. Our approach is a direct extension of the work of Dai Pra
and Posta \cite{logsob1,logsob2} for the homogeneous version.


\section{Notation and Main Results}

Throughout this paper we shall use the following notation to
denote the mean and covariance on a probability triple $(\Omega,
\mathcal{F}, \mu)$:
$$\mu[f]:= \int f d\mu, \,\,\,\, \mu[f;g]:=\mu[(f-\mu[f])(g-\mu[g])].$$
For a sub-$\sigma$-field $\mathcal{G}$ of $\mc F$, the conditional
mean and covariance is defined similarly by
\begin{eqnarray*}
\mu[f|\mathcal{G}]&=& \int f(\cdot) \mu(\cdot|\mc G),\\
\mu[f;g|\mathcal{G}]&=&\mu[(f-\mu[f|\mathcal{G}])[g-\mu[g|\mathcal{G}]]|\mathcal{G}].
\end{eqnarray*}

The entropy $H(f|\mu)$ is defined as
$$H(f|\mu)= \mu[f \log f]-\mu[f]\log \mu[f]$$
for a nonnegative $f$, and we may sometimes use the notation
$H_{\mu}(f)$.  Notice that for a density $f$ the entropy simply
becomes $\mu[f\log f]$.  The conditional entropy is defined as
$$H[f|\mathcal{G}]= \mu[f \log f|\mathcal{G}]-\mu[f|\mathcal{G}]\log\mu[f|\mathcal{G}].$$

Given a function $h$ and a set $\Lambda$, we will write
$$AV_{z\in\Lambda}h(\eta(z))$$
to denote the sample average of $h(\eta(z))$; that is,
$\frac{1}{|\Lambda|}\sum_{z\in\Lambda}h(\eta(z)).$

\bigskip

\noindent \textbf{Nearest Neighbour Inhomogeneous Zero Range.} The
inhomogeneous zero range process is a continuous time Markov
process where particles perform random walks with varying rates.
The particles move around some subset $\Lambda$ of $\mathbb{Z}^d$,
and the rate at which particles make a jump depends on the total
number of particles at the same site. Thus the particles form a
system of continuous-time interacting random walks.  The name
``zero range" comes from the notion that each particle is
interacting only with the particles at the same site, and hence
the interaction has ``no range".

We are interested in the evolution of the number of particles at
each site. To this end, let $\eta(x)$, 
denote the number of particles at site $x$ in $\Lambda$. The
function $\eta$ is an element of the space $\mc X={\mathbb
N}^{\Lambda}$.  To indicate that we are referring to the function
$\eta$ restricted to some subset $\tilde{\Lambda}\subset\Lambda$
we will use the notation $\eta_{\tilde{\Lambda}}$. For each $x \in
\Lambda$ fix a rate function $c_x:\bb{N}\mapsto\bb{R}$ such that
$c_x(0)=0$ and it is strictly positive otherwise.  A particle at
site $x$ waits independently for an exponential amount of time
with rate $c_x(\eta(x))/\eta(x)$ and then jumps from its current
position to one of its nearest neighbours $y$.  To maintain
symmetry the particle chooses either neighbour with equal
probability. Note that this implies that the first particle to
jump from site $x$ does so at rate $c_x(\eta(x))$.

If a particle moves from site $x$ to site $y$ the configuration
$\eta$ changes to $\eta^{x,y}$ where
$$(\eta^{x,y})(z) \; =\; \left\{
\begin{array}{ll}
\eta (z) - 1  & \hbox{if $z=x$}, \\
\eta (z) + 1  & \hbox{if $z=y$}, \\
\eta (z)      & \hbox{otherwise}.\;
\end{array} \right.
$$
Again, this change occurs at rate $c_x(\eta(x))$. The time
evolution of the configuration $\eta$ forms a Markov process and
we may write its generator $L$ as
\begin{align}\label{def:gen}
(L f)(\eta) &= \frac{1}{2}\sum_{x\sim y\in\Lambda}
c_x(\eta(x)) [f(\eta^{x,y})-f(\eta)]\\
&= \frac{1}{2}\sum_{x\sim y\in\Lambda} c_x(\eta(x))
\nabla_{x,y}f\notag
\end{align}
where $x\sim y$ denotes nearest neighbours of $\mathbb Z^d$ (or
$\Lambda$).

Notice that the dynamics we have described preserve the total
number of particles as the system evolves through time.  For a
configuration $\eta$, we shall denote the total number of
particles as $R=R(\eta)$, and a realization of this random
variable as $r$. Thus, for the case where $\Lambda$ is finite and
the total number of particles is $r$, the dynamics describe an
irreducible Markov process on a finite state space
$\mc{X}_r=\{\eta \in \mc X | R=r\}$. The stationary measure for
this process is denoted by $\nu_{\Lambda,r}$, and is proportional
to
$$\Pi_{x\in\Lambda}\frac{1}{c_x(\eta(x))!},$$ where we define the
factorial $c_x(k)!$ to be $c_x(k)\times c_x(k-1)\times\cdots
\times c_x(1)$, with the convention that $c_x(0)!=1$.

The canonical ensembles  $\nu_{\Lambda,r}$ satisfy the detailed
balance condition, and the system is hence reversible.   That is,
whenever $\eta(x)>0$, we have
\begin{eqnarray}\label{detbal}
c_x(\eta(x))\nu_{\Lambda,r}(\eta)=c_y(\eta(y)+1)\nu_{\Lambda,r}(\eta^{x,y}).
\end{eqnarray}
This allows us to write the Dirichlet form
$D_{\Lambda,r}(f)=\nu_{\Lambda,r}[f (-L)f]$ in the more convenient
form
$$D_{\Lambda,r}(f)=\frac{1}{2}\sum_{x\sim y \in \Lambda}\nu_{\Lambda,r}[c_x(\eta_x)\{\nabla_{x,y}f\}^2].$$

We next consider the grand canonical measures.  To this end, fix
$\varphi$ in $(0,\infty)$, and define $\mu_{\Lambda,\varphi}$ as
the product measure with marginals
$$\mu_{\{x\},\varphi}(\eta_x=k)=\frac{\varphi^k}{c_x(k)!\,Z_x(\varphi)},$$
where $Z_x(\varphi)=\sum_{k\geq0}\frac{\varphi^k}{c_x(k)!}$ is the
normalizing factor. $Z_x(\varphi)$ is also called the partition
function. The grand canonical measures continue to satisfy the
detailed balance condition.  The canonical ensembles are equal to
the grand canonical measure conditioned on the total number of
particles,
$\nu_{\Lambda,r}(\eta)=\mu_{\Lambda,\varphi}(\eta|R=r)$.

Let $\rho_x$ denote the density of particles at site $x$,
$\rho_x=\mu_{\Lambda,\varphi}(\eta_x)$.  We use the notation
$\rho_\Lambda$ to denote the average of the function $\rho_x$,
$$\rho_\Lambda = \frac{1}{|\Lambda|}\sum\rho_x.$$  For a fixed $\Lambda$,
we shall often simplify this notation to $\rho$. Note that the identity
$\varphi=\mu_{\Lambda,\varphi}[c_x(\eta_x)]$ continues to hold for
the inhomogeneous system.   Also, $\rho=\rho(\varphi)$, considered
as a function of $\varphi$, is strictly increasing. We will use
the notation $\sigma^2_x=\sigma^2_x(\varphi)$ to denote the
variance at $\eta_x$, $\mu_{\Lambda,\varphi}[\eta_x;\eta_x]$. We
also define
\begin{enumerate}
\item the average density
$$\sigma^2_\Lambda(\varphi)=\frac{1}{|\Lambda|}\sum\sigma^2_x(\varphi),$$
\item the $k^{th}$ moment at site $x$,
$m_k^x=E_{\mu_{\{x\},\varphi}}[(\eta_x-\rho_x)^k]$, and
\item the Fourier transform of the marginal
$$\hat{\mu}^x_\varphi(t)=E_{\mu_{\{x\},\varphi}}[e^{it(\eta_x-\rho_x)/\sigma_x}].$$
\end{enumerate}

For the homogeneous model , it is standard practice to index the
product measures $\mu_{\Lambda,\varphi}$ by $\rho$ instead of
$\varphi$.  This is natural as there exists a one-to-one map
between the two quantities, and $\rho$ may be interpreted as the
density, a quantity easily seen to be preserved by the system
dynamics. However, we will not do this here. Nonetheless, an
invertible relationship continues to hold in our setting. For any
fixed set $\Lambda$, there is a one-to-one relationship between
$\rho$ and $\varphi=\varphi(\rho)$.  We may then, using the
canonical measure, fully recover the function $\rho_x$.

As stated in the introduction, we shall assume that the rate
functions $c_x(\cdot)$ satisfy Lipschitz growth (LG) and weak
monoticity (M), \emph{uniformly} in $x$. The two conditions imply
that there exist universal constants $c_1$ and $c_2$ so that
\begin{eqnarray}\label{cfacts}
0<c_1 \leq \frac{c_x(k)}{k} \leq c_2<\infty
\end{eqnarray}
 for all $k$ and
$x$.  Additionally, we assume
\begin{eqnarray}\label{extracomplete}
(E) \hspace{.7cm} 0<\inf_r \sqrt{r}\, \mu_{\Lambda,
\varphi(\frac{r}{|\Lambda|}) } \left(\sum_{x\in
\Lambda}\eta(x)=r\right) \leq\sup_r \sqrt{r}\, \mu_{\Lambda,
\varphi(\frac{r}{|\Lambda|}) } \left(\sum_{x\in
\Lambda}\eta(x)=r\right) < \infty\notag
\end{eqnarray}
for any size $|\Lambda|\geq 2$. This is simply a restatement of
\reff{extracond} using the notation developed above.  Again, we
direct the reader to Remark \reff{condLLR} for further details.

We are finally in the position to state the main result.

\begin{thm}\label{bigresult}
Assume that conditions (LG),(M) and (E) are satisfied by the
inhomogeneous zero range process.  Then the system defined on
$\Lambda \subset \bb Z^d$, a cube of size $N^d$, satisfies a
logarithmic Sobolev inequality with logarithmic constant of the
order of $N^2$. That is, there exists a constant $C>0$ such that
for any choice of $r$, $|\Lambda|\geq 2,$ and non-negative
function $f,$
$$H(f|\nu_{\Lambda, r}) \leq C N^2 D_{\Lambda, r}(\sqrt{f}).$$
The constant $C$ may depend on the dimension $d$ of the cube, but
it is constant in $N$ as well as $r$, the total number of
particles.
\end{thm}

\begin{rem}
The constant $C$ in the above inequality depends also on the
parameters of the model given by the assumptions (LG), (M) and
(E).  However, we choose not to keep track of the exact form of
the dependence.
\end{rem}
As we mentioned in the introduction, the logarithmic Sobolev
inequality implies the spectral gap.  However, the proof of
Theorem \ref{bigresult} makes use of the spectral gap for the zero
range process, and so it was necessary to prove the following
beforehand.

\begin{thm}[Spectral Gap]\label{specgap}
Assume that conditions (LG) and (M) hold uniformly for the
inhomogeneous zero range process.  Then for the system defined on
$\Lambda \subset \bb Z^d$, where $|\Lambda|=N^d$, there exists a
finite constant $C^*>0$, such that
\begin{align*}
\nu_{\Lambda,r}[f;f]\leq C^* N^2 D_{\Lambda,r}(f)
\end{align*}
for all $f\geq0$.  The constant $C^*$ may depend on the dimension
$d$, as well as the constants $a_1$, $a_2$ and $k_0$ of the
assumptions.
\end{thm}
Notice that for this result we do not need the additional
assumption (E).

\bigskip

\noindent \textbf{Connection to Zero Range with Colours.}  To
simplify notation we define the $k$-colour model for the case when
$k=2$. The extension to general $k$ is obvious.

First, consider the colour-less (or \emph{colour-blind})
homogeneous zero range process.  To define this, simply take the
previously described inhomogeneous system and add the requirement
that the jump rates satisfy $c_x(\cdot)=c(\cdot)$ for all $x$. The
invariant measures now become the product measures $\mu_{\Lambda,
\varphi}$ with spatially homogeneous marginals
\begin{eqnarray}\label{zeromarg}
\mu_\varphi(\eta(x)=k)=\frac{\varphi^k}{c(k)!}Z^{-1}(\varphi).
\end{eqnarray}
$Z(\varphi)$ is again the normalizing factor.

Next, imagine that this process is made up of two different
colours of particles. The particles are mechanically identical to
the regular zero range particles, but we now also keep track of
their colour as the system evolves.  Let $\eta_i(x)$ denote the
number of particles of colour $i$ at site $x$ in $\Lambda$. Notice
that the time evolution of $\bm\eta = \{\eta_1, \eta_2\}$ is a
Markov process with state space $\mc X^2 = {\mathbb
N}^{\Lambda}\times{\bb N}^{\Lambda}$. Define the colour rate
functions
\begin{eqnarray*}
c^1(k_1,k_2) &=& k_1 \frac{c(k_1+k_2)}{k_1+k_2},\\
c^2(k_1,k_2) &=& c(k_1+k_2)-c^1(k_1,k_2).
\end{eqnarray*}
The generator for the two-colour process is then
\begin{eqnarray} (L^\mathrm{colour} f)(\bm \eta) & =
& \sum_{i=1}^2 \sum_{x \sim y} c^i(\bm{\eta}(x)) [f(\bm
\eta^{x,y}_i))-f(\bm \eta)],\label{def:zr2}
\end{eqnarray}
where $\bm \eta_i^{x,y}$ denotes the configuration obtained from
$\bm \eta$ by moving one particle of colour $i$ from site $x$ to
site $y$. Note that if the function $f$ is ``blind" to the
particle colour, i.e. $f(\bm{\eta}) = f(\eta_1 + \eta_2)$, then
$L^\mathrm{colour}f$ is equivalent to the generator for the
homogeneous zero range process.

Fix $0<\varphi_1, \varphi_2<\infty$.  The grand canonical measures
for the two colour process are the product measures with marginals
$$ \mu_{\varphi_1,\varphi_2}(\eta_1(x)=k, \eta_2(x)=m) = {k+m
\choose k} \frac{\varphi_1^k \ \varphi_2^m}{c(k+m)!}\  \
Z^{-1}(\varphi_1+\varphi_2),$$  where $Z(\varphi)$ is the same
partition function as in \reff{zeromarg}.  Further calculations
show that
\begin{eqnarray*}
\mu_{\varphi_1,\varphi_2}( \eta_1(x)=k | \eta_2(x)=0 ) &=&
\frac{\varphi_1^k}{c(k)!}Z^{-1}(\varphi_1)\\
\mu_{\varphi_1,\varphi_2}(\eta_1(x)+\eta_2(x)=n)&=&\frac{\varphi^n}{c(n)!
\
} Z^{-1}(\varphi_1+\varphi_2)\\
\mu_{\varphi_1,\varphi_2}(\eta_1(x)=k ,\eta_2(x)=n-k|\eta(x)=n )
&= &{n \choose
k}\left(\frac{\varphi_1}{\varphi_1+\varphi_2}\right)^k\left(\frac{\varphi_2}{\varphi_1
+ \varphi_2 }\right)^{n-k}.
\end{eqnarray*}
Notably, we obtain that
\begin{eqnarray*}
\mu_{\varphi_1,\varphi_2}( \eta_1(x)=k | \eta_2(x)=m ) &=&
\frac{\varphi_1^k}{c_m(k)!}Z_m^{-1}(\varphi_1)\label{def:condmeas},
\end{eqnarray*}
with $$c_m(k)=\frac{kc(k+m)}{k+m}$$ and we denote the associated
partition function as $Z_m(\cdot)$.  Equivalently, we can say that
the grand canonical measures for the first colour conditioned on
the configuration of the second colour are product measures with
marginals
\begin{eqnarray*}
\tilde \mu (\eta_1(x)=k)= \frac{\tilde\varphi^k}{\tilde
c_x(k)!}\tilde Z_x^{-1}(\tilde\varphi),
\end{eqnarray*}
where $\tilde\varphi=\varphi_1$,
$$\tilde c_x(k)=k\,\frac{c(k+\eta_2(x))}{k+\eta_2(x)},$$
and $\tilde Z_x(\tilde\varphi)$ is the normalizing constant.  We
hence obtain the invariant measures for an inhomogeneous process.
It is not difficult to show that assuming conditions \reff{LGhom}
and \reff{Mhom} imply uniform (LG) and (M) conditions for the
rates $\tilde c_x(\cdot)$.  Assuming that, for example,
$c(k)=\theta k$ for all sufficiently large $k$, gives also the
additional condition~(E).

Notice that this relationship holds regardless of the number of
colours originally considered.

\section{Outline of Proof}\label{sec:outline}

The proof of Theorem \ref{bigresult} is divided into two main
sections. First, we prove that the logarithmic Sobolev constant is
independent of the number of particles $r$.  Second, we obtain
sharper bounds which give the desired $N^2$ scaling for large
enough $N$.   The combination of these two results implies Theorem
\ref{bigresult}. In each case we use induction in the size of
$\Lambda$, where each induction step doubles the size of the cube.
For ease of presentation, we write out the proof for the case when
$d=1$.  Similar arguments as those presented in the proof of the
spectral gap (Section \ref{sec:specgap}) extend the argument to
the general case.

Suppose then that $|\Lambda|=2N$, so that we may write
$\Lambda=\Lambda_1\cup\Lambda_2$ where
$|\Lambda_1|=|\Lambda_2|=N$, and the two subsets are disjoint. We
shall denote the number of particles on $\eta_{\Lambda_i}$ as
$R_i$, with $R=R_1+R_2$.  We thus have
\begin{eqnarray}
H(f|\nu_{\Lambda, r})&=&\nu_{\Lambda, r}[H(f|\nu_{\Lambda,
r}(\cdot|R_1=r_1))]\notag\\
&&\hspace{2cm}+\,H(\nu_{\Lambda,
r}[f|R_1=r_1]|\nu_{\Lambda,r})\notag\\
&\leq &\nu_{\Lambda,r}[H(f|\nu_{\Lambda_1, r_1})]+\nu_{\Lambda,
r}[H(f|\nu_{\Lambda_2, r-r_1})]\notag\\
&&\hspace{2cm}+\,H(\nu_{\Lambda,
r}[f|R_1=r_1]|\nu_{\Lambda,r})\label{line1}
\end{eqnarray}
because
$\nu_{\Lambda,r}(\cdot|R_1=r_1)=\nu_{\Lambda_1,r_1}\otimes\nu_{\Lambda_2,r-r_1}$.
Let $\kappa(N,r)$ be the smallest constant such that
\begin{align}\label{def:kappa}
H(f|\nu_{\Lambda', r'})\leq \kappa(N,r) D_{\Lambda',r'}(\sqrt{f})
\end{align}
for all volumes $|\Lambda'|\leq N$ and $r'\leq r$ particles. The
induction hypothesis provides us with the bound on the first half
of \reff{line1} (recall that here $|\Lambda|=2N$)
\begin{eqnarray}
H(f|\nu_{\Lambda, r})&\leq &\kappa(N,r)
\nu_{\Lambda,r}\left[D_{\Lambda_1,r_1}(\sqrt{f})+D_{\Lambda_2,r-r_1}(\sqrt{f})\right]
\notag\\
&&\hspace{1cm}+H(\nu_{\Lambda, r}[f|
R_1=r_1]|\nu_{\Lambda,r})\notag\\
&\leq &\kappa(N,r)D_{\Lambda, r}(\sqrt{f})+H(\nu_{\Lambda,
R}[f|R_1=r_1]|\nu_{\Lambda,r}).\label{line2}
\end{eqnarray}
It thus remains to estimate the last term to obtain the
appropriate bounds.

In Section \ref{sec:ind}, we prove the following initial bound
\begin{align}\label{impline}
H(\nu_{\Lambda, R}[f|R_1=r_1]|\nu_{\Lambda,r})\leq
C(N)(1+\kappa(N,r))D_{\Lambda,r}(\sqrt{f})+C(N)\nu_{\Lambda,r}[f],
\end{align}
where $C(N)$ is a constant depending only on $N$. Together with
\reff{line2} this implies that
\begin{eqnarray}\label{line:new1}
H(f|\nu_{\Lambda, r})\leq C(N)\left(
\nu_{\Lambda,r}[f]+D_{\Lambda, r}(\sqrt{f})+\kappa(N,r)
D_{\Lambda, r}(\sqrt{f})\right).
\end{eqnarray}
Define the function $\tilde f =
(\sqrt{f}-\nu_{\Lambda,r}[\sqrt{f}])^2$. In Lemma
\ref{Rothaus:ineq} we give a proof of the inequality
\begin{eqnarray}\label{rothaus}
H(f|\nu_{\Lambda,r})\leq H(\tilde
f|\nu_{\Lambda,r})+2\nu_{\Lambda,r}[\sqrt f;\sqrt f],
\end{eqnarray} due to Rothaus \cite{Rot}.  This, together with
\reff{line:new1} evaluated at $\tilde f$,
 gives the new bound
$$H(f|\nu_{\Lambda, r})\leq C(N)\left( \nu_{\Lambda,r}[\sqrt{f};\sqrt{f}]+D_{\Lambda, r}(\sqrt{f})+\kappa(N,r) D_{\Lambda, r}(\sqrt{f})\right). $$
By applying the spectral gap result of Theorem \ref{specgap} we
may bound this again by
$$H(f|\nu_{\Lambda, r})\leq C(N)\left( D_{\Lambda, r}(\sqrt{f})+\kappa(N,r) D_{\Lambda, r}(\sqrt{f})\right). $$
for some different constant $C(N)$. From this it follows that
\begin{eqnarray*}
\kappa(2N,r) &\leq& C(N)(1+\kappa(N,r)).
\end{eqnarray*}
In Section \ref{sec:ind} we prove the initial induction step
\begin{eqnarray}\label{inistep}
\sup_{r}\kappa(2,r)<\infty.
\end{eqnarray}
These two bounds imply that the logarithmic Sobolev constant is
independent of $r$, the number of particles. 

The second part in the proof consists of tightening \reff{impline}
above to
\begin{align}\label{impline2}
H(\nu_{\Lambda, R}[f|R_1=r_1]|\nu_{\Lambda,r})\leq C N^2
D_{\Lambda,r}(\sqrt{f})+C\nu_{\Lambda,r}[f] +
\kappa(N,r)D_{\Lambda,r}(\sqrt{f})
\end{align}
for some $C>0$ and all large enough $N$.  After eliminating the
term $\nu_{\Lambda,r}[f]$ as before, this allows us to conclude
that for large enough $N$ we have the relationship
$$\kappa (2N,r)\leq \kappa (N,r)+ CN^2,$$
for some new constant $C$.  From this, by induction, we conclude
Theorem \ref{bigresult}.

The tighter bounds are discussed in section \ref{sec:tight}. The
better estimates are obtained by improved bounds on the
covariances appearing in Proposition \ref{rough1} of Section
\ref{sec:ind}.  These are possible because of the local limit
theorems established in Section \ref{sec:prelim}.  The local limit
theorems and resulting moment bounds are the main tools
established in this section.

Section \ref{birthdeath} looks at the spectral gap and logarithmic
Sobolev inequality for several birth and death processes.  These
are key in establishing the bounds on the second half of
\reff{line2}.

In section \ref{sec:specgap}, we prove Theorem \ref{specgap} for
any dimension $d$.  As mentioned previously, a similar approach
generalizes the arguments of Sections \ref{sec:ind} and
\ref{sec:tight} to higher dimensions.

\bigskip

\noindent\textbf{Remark on notation.}  In what follows, for
reasons of presentation, we use the notation $\eta_x$ in lieu of
$\eta(x)$.  As we no longer deal with the colour version of zero
range, there is no inconsistency to resolve.

\section{Preliminary Results\label{sec:prelim}}

In this section we establish certain tools used throughout this
paper.  We give a proof of \reff{rothaus} and discuss stochastic
monotonicity on the canonical ensembles which will be used in
Section \ref{birthdeath}.  Most notably, we establish uniform
local limit theorems for very small, small and large values of the
parameter $\varphi$.  These results are also used to give several
bounds on the zero range moments.

\begin{lem}\label{Rothaus:ineq}
Let $\mc X$ denote a Polish space. For any nonnegative function
$\phi:\mc X \mapsto \bb R$ and probability measure $\mu$ on $\mc
X$ the following inequality holds:
$$H(\phi|\mu)\leq H(\tilde \phi|\mu)+2\mu[\sqrt{\phi};\sqrt{\phi}\,],$$
where $\tilde \phi = (\sqrt{\phi}-\mu[\sqrt{\phi}\,])^2$.
\end{lem}

\begin{proof}
It is enough to consider continuous bounded functions $f$.  The
inequality follows if we can prove for any such $f$ and constant
$k$:
\begin{eqnarray}\label{line:rot1}
H(f^2|\mu)+2 \mu[f;f] \geq H((f+k)^2|\mu).
\end{eqnarray}
The proof of this for the case of $\mu$ Lebesgue measure appears
in \cite{Rot}.  To begin we define the following quantities
$$\tilde H_k(t)= H((tf+k)^2|\mu)$$
and
$$\tilde C(t)= \mu[tf; tf].$$  We calculate the first two
derivatives of the $\tilde H_k(t)$ to be
\begin{eqnarray*}
\partial_t \tilde H_k(t) = \int \{2(tf+k)f \cdot \log (tf+k)^2 \}d\mu
-\int \{2(tf+k)f \}d\mu \cdot \log \int (tf+k)^2 d\mu
\end{eqnarray*}
and
\begin{eqnarray*}
\partial_t^2 \tilde H_k(t) &=& \int \{2f^2\log (tf+k)^2\}d\mu - \int
2f^2d\mu \cdot \log \int (tf+k)^2d\mu\\
 &&+ \int (2f)^2d\mu - \frac{(\int \{2(tf+k)f\}d\mu )^2}{\int
 (tf+k)^2d\mu}.
\end{eqnarray*}
Notice that in the above there is a potential problem whenever
$(tk+f)=0$; in this case we may make the convention that $\log
(tk+f)=0$ and the above still integrates to the correct thing.

Notice also that we have
\begin{eqnarray}\label{line:rot2}
\left. \tilde H_k(t)\right|_{t=0} = 0 \hspace{1cm} \mbox{and}
\hspace{1cm} \left.\partial_t \tilde H_k(t)\right|_{t=0} = 0,
\end{eqnarray}
and that the same holds for the function $\tilde C(t)$.  Fix the
function $f$ and the constant $k$ and let $\tilde H(t)=\tilde
H_k(t)|_{k=0}$. Inequality \reff{line:rot1} may be re-written as
\begin{eqnarray*}
\tilde H(t)+2 \tilde C(t) \geq \tilde H_k(t),
\end{eqnarray*}
for all $t\geq 0$.  The inequality holds at $t=0$ and also after
taking one derivative in $t$ and evaluating at $t=0$ by
\reff{line:rot2}, as in both cases both sides are simply zero.
Hence, by integrating twice, we will obtain \reff{line:rot1} if we
can show that
\begin{eqnarray*}
\partial_t^2\tilde H(t)+2 \partial_t^2\tilde C(t) \geq \partial_t^2\tilde H_k(t),
\end{eqnarray*}
holds for all $t$.  By our previous calculations this is the same
as showing the following holds:
\begin{eqnarray*}
&&\hspace{-1in}\int \{2f^2\log (tf)^2\}d\mu - \int 2f^2d\mu \cdot
\log \int
(tf)^2d\mu \\
&\geq& \int \{2f^2\log (tf+k)^2\}d\mu - \int 2f^2d\mu \cdot \log
\int (tf+k)^2d\mu\\
&&\hspace{1cm} - \hspace{1cm}\frac{(\int \{2(tf+k)f\}d\mu
)^2}{\int
 (tf+k)^2d\mu}.
\end{eqnarray*}
This last inequality holds because by the variational formula of
the entropy we may deduce that for any $f$ and $g$ we have that
\begin{eqnarray}
\int f^2 \log f^2 d\mu - \int f^2 d\mu \cdot \log \int f^2 d\mu
\geq \int f^2 \log g^2 d\mu - \int f^2 d\mu \cdot \log \int g^2
d\mu.\notag\\\label{line:rot3}
\end{eqnarray}
Indeed, it is enough to consider $f$ such that $\int f^2 d\mu =1$.
Recall that the variational formula of entropy is
$$H(\phi|\mu)=\sup_{h \in C_b(\mc X)} \left\{ \int h \,\phi\, d\mu - \log \int e^{h} d\mu \right\}.$$
To obtain \reff{line:rot3} we need only choose $h=\log g^2$.
\end{proof}

The following result is proved for the homogeneous model in
\cite{lsv}.  It remains valid for the inhomogeneous model.

\begin{lem}\label{monotone1}
There exists a constant $B=B(a_1,a_2,k_0)$ such that
$$\nu_{\Lambda,r}\leq \nu_{\Lambda, r + M },$$
for all $M\geq B|\Lambda|$.
\end{lem}

This is an equivalent statement to the following theorem (see, for
example, \cite{lig} ).

\begin{lem}\label{monotone2}
There exists a constant $B=B(a_1,a_2,k_0)$ such that if $M\geq B
|\Lambda|$ there exists a measure $Q$ on $\mc X\times \mc X$ which
is concentrated on the configurations $\{\eta,\xi\}$ with
$\eta\leq\xi$, that is,
$$Q(\eta\leq\xi)=1,$$
and the marginals of $Q$ are $\nu_{\Lambda,r}\,\,\mbox{and}\,\,
\nu_{\Lambda, r + M }.$

\end{lem}

\begin{proof}
To clarify the proof we assume that the $k_0=2$ in assumption (M),
that is, $c_x(k)-c_x(j)\geq a_2$ for $k-j\geq 2$.  Consider the
following version of the complete zero range process with
generator
$$\mc L=\sum_{x,y\in \Lambda} c_x(\eta_x)[\nabla_{x,y}f].$$
The measures $\nu_{\Lambda,r}$ are ergodic and reversible for the
complete process as well.

Consider two configurations $\eta$ and $\xi$ where $R(\eta)=r$ and
$R(\xi)=r+M$ and $\eta\leq\xi$.  Our goal is to define a coupled
process $\{\eta(t),\xi(t)\}$ with initial configuration
$\{\eta,\xi\}$, which preserves the order $\eta(t) \leq \xi(t)$ at
all times $t$ and where the marginals of $\eta(t)$ and $\xi(t)$
evolve according to the dynamics defined by the generator $\mc L$.
To this end we define the subsets of $\Lambda$
$$b_0=\{x; \eta_x-\xi_x=0\}\,\,\mbox{and}\,\,b_1=\{x; \xi_x-\eta_x=1\}.$$

We may now begin to consider the coupling.  We need to describe
the jumps of the $\eta$ and $\xi$ particles so that for any two
configurations $\eta$ and $\xi$ with $\eta \leq \xi$ the order is
preserved after any possible jump.

First we make the particles on the sites $b_0$ make a jump
together.  Because \linebreak $c_x(k)-c_x(j)\geq a_2$ for $k-j\geq
2$ for all sites of $\Lambda \setminus b_0\cup b_1$ we may couple
all of the $\eta$ particles with a $\xi$ particle, and we are left
with $M-|b_1|$ ``free" $\xi$ particles. The only problem in the
definition of a coupling occurs if an $\eta$ particle jumps from a
$b_1$ site to a $b_0$ site. This event occurs at rate of at most
$a_1|b_0||b_1|$. We may compensate for this behaviour using our
``free" $\xi$ particles. The rate at which these particles make a
jump to $b_0$ is
$$b_0 \sum_{x \in \Lambda
\setminus b_0\cup b_1} c_x(\xi_x)-c_x(\eta_x)$$ which is bounded
below by
\begin{eqnarray*}
b_0 \sum \left\{a_2 \lfloor\frac{\xi_x-\eta_x}{k_0}\rfloor-a_1k_0
\right\}&\geq & b_0 \sum \left\{a_2
\left(\frac{\xi_x-\eta_x}{k_0}-1\right)-a_1k_0 \right\}\\
&\geq& \frac{a_2 b_0}{k_0}(M-|b_0|) - (a_2|b_0|+a_1k_0)|\Lambda
\setminus b_0\cup b_1|\\
&\geq& \frac{a_2 b_0}{k_0} M - b_0 |\Lambda|\left\{a_2
+\frac{a_2}{k_0}+a_1k_0\right\}.
\end{eqnarray*}
This will be greater than $a_1|b_0||b_1|$ as long as $M\geq
B|\Lambda|$ for $B=\frac{a_1}{a_2}k_0(k_0+1)+k_0+1.$ Finally,
because the rate at which the ``free" $\xi$ particles jump to
$b_0$ sites is greater than the rate at which the uncoupled $\eta$
particles do so, we may couple these jumps while preserving the
correct marginal dynamics.  All other particles are allowed to
jump freely.  This is exactly what we need so that our joint
process preserves the order $\eta(t) \leq \xi(t)$ at all times
$t$.
\end{proof}

\subsubsection{Local Limit Theorems.}

We begin with some moment bounds.

\begin{prop}\label{sigmabounds}
\begin{enumerate}
\item For all $x$ in $\Lambda$
$$0 < c_1 \leq \frac{\varphi}{\rho_x}\leq c_2 < \infty.$$
\item There exist constants
$0<\tilde{c}_1 \leq \tilde{c}_2 < \infty$ such that for all $x$ in
$\Lambda$
$$ \tilde{c}_1 \leq \frac{\sigma^2_x(\varphi)}{\varphi} \leq \tilde{c}_2.$$
The constants $\tilde c_1$ and $\tilde c_2$ depend only on the
values $a_1$, $a_2$ and $k_0$.
\end{enumerate}
\end{prop}

The first inequality is a simple consequence of \reff{cfacts}.
The proof of the second inequality appears in \cite{lsv}.  As the
bounds are bounds on the single site marginals, they also apply in
this setting.  Because the bounds depend only on the constants
$a_1$, $a_2$ and $k_0$ they apply uniformly to all $x$.  Notice
that we do not require the additional assumption (E).  The same is
true of the following:
\begin{prop}\label{hatbounds}
For each $\varphi_0>0$ and $\bar{k}$ in $\mathbb{N}$
\begin{enumerate}
\item There exists a finite constant $K_0$ such that
$$\sup_{\varphi\geq\varphi_0} \frac{m^x_{2k}(\varphi)}{\sigma^{2k}(\varphi)} \leq K_0,\,\,\, for \,\, 1\leq k \leq \bar{k}.$$
\item For every $\delta>0$, there exists $C(\delta)<1$ such that
$$\sup_{\varphi\geq\varphi_0} \sup_{\delta\leq|t|\leq\pi\sigma(\varphi)}|\hat\mu^x_\varphi(t)|\leq C(\delta).$$
\item There exists a $\kappa>0$ so that
$$\sup_{\varphi\geq\varphi_0}\int_{|t|\leq\pi\sigma(\varphi)}|\hat\mu^x_\varphi(t)|^\kappa\leq C<\infty.$$
\end{enumerate}
\end{prop}

From the above and a simple calculation we also obtain

\begin{cor}\label{phideriv}
There exists a finite, positive constant $\tilde c$ such that for
all $\Lambda$
$$\tilde c^{-1} \leq \varphi'(\rho)=\frac{\varphi}{\sigma^2(\rho)}\leq \tilde c .$$
\end{cor}

We now turn to the local limit theorems.  Recall the definition of
the Hermite polynomial of degree $m$, for $m\geq0$:
$$H_m(x)=(-1)^m \exp {(\ \frac{x^2}{2})} \frac{d^m}{dx^m} \exp {(-\frac{x^2}{2})}.$$
Let $g_0(x)$ denote the density of a standard normal random
variable, and, for $j\geq1$, define
\begin{align}
g_j(x)=g_0(x) \sum_{} H_{j+2a}(x) \prod_{m=1}^j
\frac{1}{k_m!}\left(\frac{\kappa_{m+2}}{(m+2)!\,\,\sigma^{m+2}}\right)^{k_m}
\end{align}
where the sum is taken over all nonnegative integer solutions
$\{k_l\}_{l=1}^j$ and $a$ of $k_1+2k_2+\hdots+jk_j=j$ and
$k_1+k_2+\hdots+k_j=a$, and $\kappa_m$ denotes the $k^{th}$
cumulant of the distribution.

In what follows, we assume that $|\Lambda|=N$.

\begin{thm}\label{normalLLT-2}
\begin{enumerate}
\item For all $\varphi_0>0$ and $J \in \mathbb{N}$, there exist
finite constants $E_0=E_0(\varphi_0,J)$ and $A=A(\varphi_0,J)$
such that
\begin{align*}
\left|\sqrt{N\sigma^2}\mu_{\Lambda,\varphi}\left[\sum_{x\in\Lambda}\eta(x)=N\rho+\sigma\sqrt{N}z\right]-\sum_{j=0}^{J-2}\frac{1}{N^{j/2}}g_j(z)\right|
\leq \frac{E_0}{(\sigma^2 N)^{(J-1)/2}}
\end{align*}
uniformly over $z$ and over all parameters $A/N \leq \varphi \leq
\varphi_o.$  In the above, $\sigma = \sigma_{\Lambda}(\varphi)$.
\item  For any $\varphi_0>0$ and any $k_1 \in \bb{N}$, there exists a constant
$E_1=E_1(k_1,\varphi_0)$ and $n_1=n_1(k_1,\varphi_0)$ such that
\begin{align*}
\left|\sqrt{N\sigma^2}\mu_{\Lambda,\varphi}\left[\sum_{x\in\Lambda}\eta(x)=N\rho+\sigma\sqrt{N}z\right]-\sum_{j=0}^{k_1-2}\frac{1}{N^{j/2}}g_j(z)\right|
\leq \frac{E_1}{N^{(k_1-1)/2}}
\end{align*}
for all $N> n_1$ uniformly over $z$ and $\varphi>\varphi_0$.
Again, in the above, $\sigma = \sigma_{\Lambda}(\varphi)$.
\end{enumerate}
\end{thm}

\begin{proof}[Sketch of Proof.]
The proof of this result is discussed in \cite{lsv} for the
homogeneous case. It relies on repeating the usual local limit
theorem argument (see for example \cite{gk} or \cite{pet}), while
checking that the bounds are valid uniformly in the parameter
$\varphi$ in the two settings. These bounds rely on moment
estimates, which, due to conditions (LG) and (M) are identical in
both the homogeneous and nonhomogeneous models.

Fix $J$ in $\bb N$. Following the usual local limit theorem
approach as in \cite{pet} we write
\begin{eqnarray*}
\sqrt{N\sigma^2}\mu_{\Lambda,\varphi}\left[\sum_{x\in\Lambda}\eta_x=N\rho+\sigma\sqrt{N}z\right]-\sum_{j=0}^{J-2}\frac{1}{N^{j/2}}g_j(z)
\end{eqnarray*}
as
\begin{eqnarray*}
\int_{-\pi \sqrt n \sigma}^{\pi \sqrt n \sigma} e^{-i t z}
f_N(t)dt - \int e^{-itz}u_{k,N}(t)dt,
\end{eqnarray*}
where $u_{k,N}(t)=\sum_{j=0}^{J-2}\frac{1}{N^{j/2}}\int
e^{itx}g_j(x)dx$ and $f_N(t)$ is the characteristic function of
$(\sqrt{N}\sigma)^{-1} \sum_{x\in\Lambda}(\eta_x-\rho_x)$.  We may
bound this quantity by the sum of 4 integrals
\begin{eqnarray*}
I_1 &=& \int_{|t|<N^{1/6}} |f_N(t)-u_{k,N}(t)|dt\\
I_2 &=& \int_{N^{1/6}\leq|t|<T_N} |f_N(t)|dt\\
I_3 &=& \int_{T_N\leq |t|<\pi \sqrt{N}\sigma} |f_N(t)|dt\\
I_4 &=& \int_{|t|>N^{1/6}} |u_{k,N}(t)|dt,
\end{eqnarray*}
where $T_N = \frac{1}{4}N^{3/2}\sigma^{3/2}\sum_{x=1}^N m^3_x$.
The details on how to bound $I_1$, $I_2$ and $I_4$ are given in
\cite{pet}: Theorem 12 of chapter 7 and Lemmas 11 and 12 of
chapter 6. We obtain, for $i=1,2,4$, the bound
$$I_i \leq \frac{C}{(\sigma^2(\varphi)N)^{(J-1)/2}},$$
where $C$ is a positive constant not depending on $\varphi$.
Notice that under the second regime we have that
$\sigma^2(\varphi)>C\varphi_0$, for all $\varphi>\varphi_0$.   The
main difficulty lies in obtaining the appropriate uniform bounds
on $I_3$. Because of the moment bounds at the beginning of this
section we have for some $b$
\begin{eqnarray*}
|I_3|\leq \sqrt{N}\sigma \int_{b<|t|<\pi}\Pi_{x\in\Lambda
}|\mu^x_{\varphi}(t)|dt.
\end{eqnarray*}
The result follows if we can show that $|\mu^x_{\varphi}(t)|$
stays strictly below 1 uniformly in $\varphi$ and $x$. In part (1)
of the theorem this follows from the following bound
$$|\mu^x_{\varphi}(t)|^2-1\leq C \varphi_0 \varphi(\cos t -1).$$  In
part (2) of the theorem the necessary bound follows from part (2)
of Proposition \ref{hatbounds}.  This gives us exponential decay
on the term $I_3$.
\end{proof}

The above are both Edgeworth expansions for a lattice
distribution, which are valid uniformly for the family of measures
$\mu_{\varphi}$. Note that we require that the average density is
not ``too small". This should not be surprising, as in the case
where $\rho<A/N$ there is at most a finite number of particles,
and so the probability of observing a ``success" decreases as
$N\rightarrow\infty$ and a Poisson limit theorem holds instead of
a Gaussian limit theorem.

\begin{lem}\label{poissonLLT-2}
For  every $A \in \mathbb{N}\backslash\{0\}$ there exists a
constant $A_0$ such that
$$\sup_{r \leq A}\left|\mu_{\Lambda,\varphi_\Lambda(\frac{r}{N})}\left(\sum_{x\in\Lambda}\eta_x=k\right)-\frac{r^k}{k!}e^{-r}\right| \leq \frac{A_0}{N}$$
for any $k\in\mathbb{N}$, where $N=|\Lambda|$.
\end{lem}

\begin{proof}
This lemma is proved in \cite{logsob1} for the homogeneous system.
We extend it here to the inhomogeneous case by using the general
version of the Poisson limit theorem as proved in \cite{B}, for
example.
\begin{eqnarray*}
\mu_{\Lambda, \varphi}(R=k)&=&\mu_{\Lambda,\varphi}(R=k|\max
\eta_x \leq 1
)\mu_{\Lambda,\varphi}(\max \eta_x \leq 1 )\\
&&+\mu_{\Lambda,\varphi}(R=k|\max \eta_x
>1)\mu_{\Lambda,\varphi}(\max \eta_x >1 )
\end{eqnarray*}
Since we choose $\varphi=\varphi_\Lambda(r/N)$ and because $r$ is
bounded above we have that \linebreak $\mu_{\Lambda,\varphi}(\max
\eta_x
>1 )=O(\frac{1}{N})$.  We thus need to show that
$$\mu_{\Lambda,\varphi}(R=k|\max \eta_x \leq 1 )=\frac{r^k}{k!}e^{- r}+O(N^{-1}).$$
Define
\begin{eqnarray*} \tilde r_x &=&
\frac{\mu_{\Lambda,\varphi}[\eta_x\bb I_{ \eta_x \leq 1
}]}{\mu_{\Lambda,\varphi}(\eta_x \leq 1 )},
\end{eqnarray*}
and notice that this is equal to
\begin{eqnarray}
\tilde r_x &=&\frac{\rho_x-\mu_{\Lambda,\varphi}[\eta_x\bb I_{
\eta_x \leq 1
}]}{\mu_{\Lambda,\varphi}(\eta_x\leq 1)}\notag\\
&=&\frac{\rho_x+ O(N^{-2})}{1+O(N^{-2})}.\label{pp2}
\end{eqnarray}
Next consider the interval $[0,1]$.  Define $I_0(p)=[0,1-p]$ and
$I_1(p)=[1-p,1]$.  Also define $J_0=[0,e^{-p}]$ and
$J_m=[e^{-p}\frac{p^m}{m!},e^{-p}\frac{p^{m+1}}{(m+1)!}]$ for
$m\geq 1$.  Notice that because $1-p \leq e^{-p}$ we have
\begin{eqnarray}\label{pp1}
J_1(p) \subset I_1(p).
\end{eqnarray}
Consider a sequence of independent, Uniform $[0,1]$ random
variables $U_x$. Define
$$V_x = \left\{\begin{array}{cc}
        0 & \mbox{if \,\,$U_x \in I_0(\tilde r_x)$}\\
        1 & \mbox{otherwise.}
        \end{array}\right.$$
Also define
$$W_x = i \,\,\,\mbox {if} \,\,\,\,\, U_x \in J_i(\tilde r_x).$$
First notice that in the system we have defined we have
$$P\left(\sum_{x \in \Lambda} V_x =k\right)=\mu_{\Lambda,\varphi}(R=k|\max \eta_x \leq 1 ),$$
whereas
$$P\left(\sum_{x \in \Lambda} W_x =k\right)= e^{-\tilde r} \frac{(\tilde r)^k}{k!},$$
where we let $\tilde r = \sum_{x\in\Lambda} \tilde r_x$.  Because
of \reff{pp2} we have
$$P\left(\sum_{x \in \Lambda} W_x =k\right)= e^{- r} \frac{r^k}{k!} + O(N^{-1}).$$
Also, because of \reff{pp1}, we have
\begin{eqnarray*}
P(W_x \neq V_x) &=& P(U_x \in I_1(\tilde r_x) \setminus J_1(\tilde
r_x) )\\
&\leq& \tilde r_x - e^{-\tilde r_x}\tilde r_x \leq \tilde r_x^2.
\end{eqnarray*}
Thus
\begin{eqnarray*}
P\left(\sum_{x \in \Lambda} W_x \neq \sum_{x \in \Lambda}
V_x\right) \leq \sum_{x \in \Lambda} \tilde r_x^2 = O(N^{-1}).
\end{eqnarray*}
The result follows once we note that
\begin{eqnarray*}
P\left(\sum_{x \in \Lambda} V_x= k\right)&=&P\left(\sum_{x \in
\Lambda} W_x =k\right) + P\left(\left\{\sum_{x \in \Lambda} V_x
=k\right\} \bigcap \left\{\sum_{x \in \Lambda}
W_x \neq \sum_{x \in \Lambda} V_x \right\}\right)\\
&&-\,\,\,P\left(\left\{\sum_{x \in \Lambda} W_x =k\right\} \bigcap
\left\{\sum_{x \in \Lambda} W_x \neq \sum_{x \in \Lambda}
V_x\right\} \right).
\end{eqnarray*}
\end{proof}

The above lemmas imply the following result.  It follows by direct
computation of the conditional probability and applying the bounds
on the grand canonical measures obtained above.

\begin{cor}\label{subprob}
Fix $\delta \in (0,1)$; there exists an $n_0$ and $A>0$ such that
for any $\Lambda_0 \subset \Lambda$ with $|\Lambda|>n_0$ and
$\frac{|\Lambda_0|}{|\Lambda|}\leq \delta$
$$\nu_{\Lambda,r}(\eta_{\Lambda_0})\leq A \mu_{\varphi_\Lambda(\frac{r}{|\Lambda|})}(\eta_{\Lambda_0})$$
for all $r>0$.
\end{cor}

\begin{proof}
Let $R_0$ denote the total number of particles in $\Lambda_0$ and
$R_1$ are the total number of particles in $\Lambda \setminus
\Lambda_0$, respectively. We set $R_0(\eta_{\Lambda_0})=r_0$.  Let
$|\Lambda|=N$ and \linebreak $\varphi=\varphi_\Lambda(r/N)$.  We
begin by computing
$$\nu_{\Lambda,r}(\eta_{\Lambda_0})= \mu_{\Lambda,\varphi}(\eta_{\Lambda_0}) \frac{\mu_{\Lambda,\varphi}(R_1=r-r_0)}{\mu_{\Lambda,\varphi}(R=r)}.$$
We thus need to bound
$$\frac{\mu_{\Lambda,\varphi}(R_1=r-r_0)}{\mu_{\Lambda,\varphi}(R=r)}$$
uniformly.  This follows from the local limit theorems just
described.

\bigskip

\noindent \textbf{Case 1: $r\leq A$.}

\medskip

By Proposition \ref{poissonLLT-2} we have

\begin{eqnarray*}
\frac{\mu_{\Lambda,\varphi}(R_1=r-r_0)}{\mu_{\Lambda,\varphi}(R=r)}&\leq&
\frac{1}{\mu_{\Lambda,\varphi}(R=r)} \\
&\leq& \frac{1}{\frac{r^r e^{-r}}{r!} - \frac{A_0}{N}}\\
&\leq& \frac{2}{ \inf_{0<r\leq A}\frac{r^r e^{-r}}{r!}},
\end{eqnarray*}
for some large $N_1$ and all $N\geq N_1$.

\bigskip

\noindent \textbf{Case 2: $A \leq r \leq \rho_0 N$.}

\medskip

In this case we write
\begin{eqnarray*}
\frac{\mu_{\Lambda,\varphi}(R_1=r-r_0)}{\mu_{\Lambda,\varphi}(R=r)}&=&
\sqrt{\frac{N}{N-|\Lambda_0|}}\frac{\sqrt{\sigma^2 (N-|\Lambda_0|)}\mu_{\Lambda,\varphi}(R_1=r-r_0)}{\sqrt{\sigma^2 N}\mu_{\Lambda,\varphi}(R=r)}. \\
\end{eqnarray*}

By the first part of Theorem \ref{normalLLT-2}, setting $J=2$, we
have
\begin{eqnarray*}
\sqrt{\sigma^2
(N-|\Lambda_0|)}\mu_{\Lambda,\varphi}(R_1=r-r_0)\leq
\frac{1}{\sqrt{2\pi}}+\frac{E_0}{\sqrt{\sigma^2(N-|\Lambda_0|)}}\leq\frac{1}{\sqrt{2\pi}}+\frac{E_0}{\sqrt{A
(1-\delta)N }}
\end{eqnarray*}
and
\begin{eqnarray*}
\sqrt{\sigma^2 N}\mu_{\Lambda,\varphi}(R=r)\geq
\frac{1}{\sqrt{2\pi}}-\frac{E_0}{\sqrt{\sigma^2N}}\geq\frac{1}{\sqrt{2\pi}}-\frac{E_0}{\sqrt{A
N }}.
\end{eqnarray*}
This implies that
\begin{eqnarray*}
\frac{\mu_{\Lambda,\varphi}(R_1=r-r_0)}{\mu_{\Lambda,\varphi}(R=r)}&\leq&
\frac{1}{\sqrt{1-\delta}}\frac{\frac{1}{\sqrt{2\pi}}+\frac{E_0}{\sqrt{A
(1-\delta)N }}}{\frac{1}{\sqrt{2\pi}}-\frac{E_0}{\sqrt{A N }}}\\
&\leq&\frac{4 E_0}{\sqrt{A}(1-\delta)},
\end{eqnarray*}
for a sufficiently large $N_1$ and all $N\geq N_1$.

\bigskip \noindent \textbf{Case 3: $r > \rho_0 N$.}

\medskip

Again we write
\begin{eqnarray*}
\frac{\mu_{\Lambda,\varphi}(R_1=r-r_0)}{\mu_{\Lambda,\varphi}(R=r)}&=&
\sqrt{\frac{N}{N-|\Lambda_0|}}\frac{\sqrt{\sigma^2 (N-|\Lambda_0|)}\mu_{\Lambda,\varphi}(R_1=r-r_0)}{\sqrt{\sigma^2 N}\mu_{\Lambda,\varphi}(R=r)} \\
\end{eqnarray*}

This time by the second part of Proposition \ref{normalLLT-2} we
have
\begin{eqnarray*}
\sqrt{\sigma^2
(N-|\Lambda_0|)}\mu_{\Lambda,\varphi}(R_1=r-r_0)\leq
\frac{1}{\sqrt{2\pi}}+\frac{E_1}{\sqrt{N-|\Lambda_0|}}\leq\frac{1}{\sqrt{2\pi}}+\frac{E_1}{\sqrt{
(1-\delta)N }}
\end{eqnarray*}
and
\begin{eqnarray*}
\sqrt{\sigma^2 N}\mu_{\Lambda,\varphi}(R=r)\geq
\frac{1}{\sqrt{2\pi}}-\frac{E_1}{\sqrt{N}}\geq\frac{1}{\sqrt{2\pi}}-\frac{E_1}{\sqrt{
N }}
\end{eqnarray*}
This implies that
\begin{eqnarray*}
\frac{\mu_{\Lambda,\varphi}(R_1=r-r_0)}{\mu_{\Lambda,\varphi}(R=r)}&\leq&
\frac{1}{\sqrt{1-\delta}}\frac{\frac{1}{\sqrt{2\pi}}+\frac{E_1}{\sqrt{
(1-\delta)N }}}{\frac{1}{\sqrt{2\pi}}-\frac{E_1}{\sqrt{ N }}}\\
&\leq&\frac{4}{\sqrt{1-\delta}},
\end{eqnarray*}
for a sufficiently large $N_2$ and all $N\geq N_2$. To complete
the proof we choose \linebreak $n_0=\max\{N_0,N_1,N_2\}$.

\end{proof}

A direct consequence of this is the following result.
\begin{cor}\label{submom}
Under the conditions defined above, given a function $f$ whose
support is restricted to $\Lambda_0$, there exists a constant $C$
and an $n_0$ such that
\begin{enumerate}
\item $\nu_{\Lambda,r}[f]\leq C \mu_{{\varphi(\frac{r}{|\Lambda|})}}[|f|]$
\item $|\nu_{\Lambda,r}[f]-\mu_{{\varphi(\frac{r}{|\Lambda|})}}[f]|\leq C\sqrt{\mu_{\varphi(\frac{r}{|\Lambda|})}[f;f]}$
\end{enumerate}
uniformly for $|\Lambda|\geq n_0$ and all $r>0$.
\end{cor}
The first result follows directly from Corollary \ref{subprob}.
The second result follows from the first by plugging in
$f-\mu_\varphi(f)$ into the first and applying the Cauchy-Schwarz
inequality.

If instead we consider a function which is local, with support of
$\Lambda_0\subset\Lambda$,we may apply the local central limit
theorem in a similar manner as in the corollaries to obtain the
following result.

\begin{cor}\label{submom2}
Consider a function $f$ whose support is again $\Lambda_0$. Fix
$\rho_0>0$. In the formulae below let
$\varphi=\varphi_\Lambda\left(\frac{r}{|\Lambda|}\right)$.
\begin{enumerate}
\item There exists constants $C$ and $A$ and an $n_0$ such that
\begin{align*}
&|\nu_{\Lambda,r}[f]-\mu_{{\varphi}}[f]|\\&\leq
C\frac{|\Lambda_0|}{|\Lambda|}\left\{\frac{1}{\sigma^2_\Lambda(\varphi)}\mu_{\Lambda,\varphi}[|f-\mu_{\Lambda,\varphi}[f]|]+\frac{1}{\sigma_\Lambda(\varphi)}\sqrt{\mu_{\Lambda,\varphi}[f;f]}\right\}
\end{align*}
for any $|\Lambda|\geq n_0$ and all $r$ such that
$\frac{A}{|\Lambda|}\leq \frac{r}{|\Lambda|} \leq \rho_0$.
\item There exists a constant $C$ and an $n_0$ such that
$$|\nu_{\Lambda,r}[f]-\mu_{{\varphi}}[f]|\leq C\frac{|\Lambda_0|}{|\Lambda|}\sqrt{\mu_{\Lambda,\varphi}[f;f]}$$
uniformly for $|\Lambda|\geq n_0$ and all
$\frac{r}{|\Lambda|}>\rho_0$.
\end{enumerate}
We may choose $n_0$ to be the same in both of these cases.
\end{cor}

\begin{proof}
Consider the second case, and denote by $\xi$ the configuration
$\eta$ restricted to the subset $\Lambda_0$.  We continue using
notation from the previous result.  We may write
\begin{eqnarray}\label{submomproof}
\nu_{\Lambda,r}[f]-\mu_{\Lambda,\varphi}[f]\leq\sum_{\xi}\left|f(\xi)-\mu_{\Lambda,\varphi}[f]\right|
\mu_\varphi(\xi)\left\{\frac{\mu_\varphi[R_1=r-r_0]}{\mu_\varphi[R=r]}-1\right\}.
\end{eqnarray}
We therefore need to bound the difference inside the brackets to
complete the proof.  As before, we use the second part of the
local limit theorem with $J=3$.
\begin{align*}
&\hspace{2cm}\frac{\mu_\varphi[R_1=r-r_0]}{\mu_\varphi[R=r]}-1 \\
&\hspace{2cm}\leq \ C \{ \sqrt{\sigma^2
(N-|\Lambda_o|)}\mu_\varphi[R_1=r-r_0]-\sqrt{\sigma^2
N}\mu_\varphi[R=r]\}\\
&\hspace{2cm}\leq \ \ C \left\{\frac{E_0}{N} + g_0(0)-
g_0\left(\frac{r_0}{\sigma
\sqrt{N-|\Lambda_0|}}\right)-\frac{1}{\sqrt{N}}g_1\left(\frac{r_0}{\sigma
\sqrt{N-|\Lambda_0|}}\right)\right\}\\
&\hspace{2cm}\leq \ \ C \frac{C(E_0)|\Lambda_0|}{N}
\left\{1+\frac{r_0}{\sigma|\Lambda_0|}+\left(\frac{r_0}{\sigma|\Lambda_0|}\right)^2\right\}.
\end{align*}
We plug this estimate into \reff{submomproof} and apply
Cauchy-Schwarz inequality.  To finish we apply the bounds from the
first part of Proposition \ref{hatbounds} to get rid of the extra
moment terms.  A similar argument proves the first case, using the
first part of Theorem \ref{normalLLT-2}.
\end{proof}
Using the same approach, but setting $J=4$ in the local limit
theorem expansion, we obtain a further decomposition:
\begin{cor}\label{submom7}
Consider a function $f=f(\eta_x)$ for a fixed $x\in\Lambda$.  Fix
$\rho_0>0$. In the formulae below let
$\varphi=\varphi_\Lambda\left(\frac{r}{|\Lambda|}\right)$. There
exists constants $C$ and $A$ and an $n_0$ such that
\begin{align*}
&|\nu_{\Lambda,r}[f]-\mu_{{\varphi}}[f]-\frac{1}{2|\Lambda|}\{\mu_{\Lambda,\varphi}[f;\eta_x-\rho_x]+\mu_\varphi[f;(\eta_x-\rho_x)^2]\}|\\&\leq
C\frac{1}{|\Lambda|^{3/2}}\left\{\frac{1}{\sigma^2_\Lambda(\varphi)}\sqrt{\mu_{\Lambda,\varphi}[f;f]}\right\}
\end{align*}
for any $|\Lambda|\geq n_0$ and all $r$ such that
$\frac{A}{|\Lambda|}\leq \frac{r}{|\Lambda|} \leq \rho_0$.
\end{cor}

\bigskip

\section{Some Birth and Death Processes.\label{birthdeath}}

As we mentioned previously, in this section we establish spectral
gap and logarithmic Sobolev inequalities for several birth and
death chains which arise naturally in the study of the second half
of \reff{line2}.

\subsection{Some Spectral Gaps.} We begin by stating a result whose proof appears
as Lemma 4.3 in \cite{lsv}.

\begin{lem}\label{bdspecgap}
Let $Y_t$ be a birth and death process on $\{0,1,\ldots,r\}$ with
death rate $d(\cdot)$ and birth rate $b(\cdot)$.  Assume that
there exists a finite positive constant $J_2$ such that
\begin{eqnarray}\label{bdassume2}
\sup_k |b(k+1)-b(k)| \leq J_2,
\end{eqnarray}
and that there exist finite constants $J_0>0$ and $J_1>J_2$ such
that
\begin{eqnarray}\label{bdassume}
d(k)-d(j)\geq J_1 (k-j) -J_0,
\end{eqnarray}
for all $k\geq j$. Then the spectral gap for this process is
bounded below by a strictly positive constant $\lambda$ depending
on $J_0, J_1, J_2$ and $d^* = \min_{k\geq 1} d(k)$.
\end{lem}

From the above Lemma follow the next two results.

\begin{lem}\label{Kawasakispecgap}
Under assumptions (LG) and (M), there exists a constant
$B_0=B_0(a_1, a_2, k_0)$ such that
$$\mu_{\Lambda,\varphi}[\phi;\phi]\leq B_0 \ \mu_{\Lambda,\varphi}[c_x(\eta_x)\{\phi(\eta_x-1)-\phi(\eta_x)\}^2]$$
for all functions $\phi: \bb N \mapsto \bb R$ with
$\mu_{\Lambda,\varphi}[\phi^2] < \infty$.
\end{lem}

\begin{proof}
A straightforward calculation using the properties of the marginal
$\mu_\varphi$ shows that
\begin{align}\label{bddirformcalc}
&\mu_{\Lambda,\varphi}[c_x(\eta_x)\{\phi(\eta_x-1)-\phi(\eta_x)\}^2]\notag\\
&=
\frac{1}{2}\mu_{\Lambda,\varphi}[c_x(\eta_x)\{\phi(\eta_x-1)-\phi(\eta_x)\}^2]+\frac{\varphi}{2}\mu_{\Lambda,\varphi}[\{\phi(\eta_x+1)-\phi(\eta_x)\}^2]
\end{align}
 This is the Dirichlet form corresponding to the
birth and death chain with death rate $c_x(\cdot)$ and birth rate
$\varphi$. Since the birth rate is constant, and the death rate
satisfies assumption \reff{bdassume} with $J_1=a_2/k_0$ and
$J_0=a_2+a_1k_0$, the result follows.  Also, since the death rate
$c_x(k)\geq c_1 k$ for all $k$ we obtain a uniform lower bound
$d^*,$ which depends on $a_1$ and $a_2$.
\end{proof}

\begin{lem}\label{onesitesg}
Under the uniform assumptions (LG) and (M), there exists a
constant $B_0=B_0(a_1,a_2,k_0)$ such that
\begin{eqnarray}\label{line:bdsg}
\nu_{\Lambda,r}[f;f]\leq B_1
\nu_{\Lambda,r}[c_x(\eta_x)\{f(\eta_x-1)-f(\eta_x)\}^2]
\end{eqnarray}
for all functions $f: \bb N \mapsto \bb R$ with
$\nu_{\Lambda,r}[f^2] < \infty$, and all $\Lambda$ such that
$|\Lambda|\geq 1$ and $r\geq1$.
\end{lem}

\begin{proof}
Note that the marginals $\nu_{\Lambda,r}$ satisfy the relationship
$$\nu_{\Lambda,r}(\ \cdot\ |\eta_x=k)=\nu_{\Lambda_x,r-k}(\ \cdot\ ),$$
where $\Lambda_x = \Lambda \setminus \{x\}$. Using this fact, a
similar calculation to the one in \reff{bddirformcalc} reveals
that this corresponds to the birth and death chain with death rate
$c_x(\cdot)$ and birth rate $b_x(k)=AV_{y\sim
x}\nu_{\Lambda_x,r-k}[c_y(\eta_y)]$.  As above, we need only show
that the birth rate satisfies the necessary conditions.  That is,
we need to show that there exists a constant $C$ such that for any
$\Lambda$, $r$ and site $x$,
\begin{eqnarray*}
|\nu_{\Lambda,r+1}[c_x(\eta_x)]-\nu_{\Lambda,r}[c_x(\eta_x)]|\leq
C
\end{eqnarray*}
 We split this up
into three cases. Let $\rho=\frac{r}{|\Lambda|}$ and fix
$\rho_0>0$.  From Proposition \ref{submom2} we choose an $n_0$ and
$A$.

We first assume that $\frac{A}{|\Lambda|}\leq \rho \leq \rho_0$.
Choosing $\phi(\eta)=c_x(\eta_x)$ in the first part of Proposition
\ref{submom2} we have that
$|\nu_{\Lambda,r}[c_x(\eta_x)]-\varphi(\rho)|\leq
\frac{1}{|\Lambda|}C(\rho_0),$ for some constant C depending on
$\rho_0$. From this it follows that
\begin{eqnarray*}
|\nu_{\Lambda,r+1}[c_x(\eta_x)]-\nu_{\Lambda,r}[c_x(\eta_x)]|\leq
2C(\rho_0)/|\Lambda|.
\end{eqnarray*}

Next assume that $\rho>\rho_0$.  Here we may write
\begin{eqnarray*}
&\hspace{-8cm}|\nu_{\Lambda,r+1}[c_x(\eta_x)]-\nu_{\Lambda,r}[c_x(\eta_x)]|\\
&\leq \sum_{k\geq 0}
|c_x(k)-\varphi(\rho)|\mu_{\Lambda,\varphi}(k)\left\{\frac{\mu_{\Lambda,\varphi}(R=r)\mu_{\Lambda,\varphi}(R=r+1-k)-\mu_{\Lambda,\varphi}(R=r-k)\mu_{\Lambda,\varphi}(R=r+1)}{\mu_{\Lambda,\varphi}(R=r+1)\mu_{\Lambda,\varphi}(R=r)}\right\}.
\end{eqnarray*}
Using the first part of Proposition \ref{normalLLT-2} with the
expansion up to $J=3$, we bound the term inside the brackets by
$C\frac{k}{|\Lambda|}$, from which the desired bound follows.

To handle the last case, namely $r\geq A$ and $|\Lambda|\leq n_0$,
we use the monotonicity results of Lemma \ref{monotone2}.  We fix
$B$ sufficiently large as in the requirement of the lemma, and set
$M=B|\Lambda|$.  There exists a measure Q on $\mc X \times \mc C$
with marginals $\nu_{\Lambda,r}$ and $\nu_{\Lambda,r+M+1}$ such
that $\nu_{\Lambda,r}\leq \nu_{\Lambda,r+M+1}$. We thus have
\begin{eqnarray*}
&\hspace{-9cm}\left|\nu_{\Lambda,r+1}[c_x(\eta_x)]-\nu_{\Lambda,r}[c_x(\eta_x)]\right|\\
&\leq \ \
|\nu_{\Lambda,r+1}[c_x(\eta_x)]-\nu_{\Lambda,r+M+1}[c_x(\eta_x)]|
+\ \
|\nu_{\Lambda,r+M+1}[c_x(\eta_x)]-\nu_{\Lambda,r}[c_x(\eta_x)]|,
\end{eqnarray*}
which is smaller than $2a_1$ by assumption (LG) on the rates. This
last fact, together with a bound on the finitely many remaining
rates for $k\leq A$ proves that the birth rates satisfy
\reff{bdassume2}.   Notice again that the lower bound $d^*$ is
uniform in the sites.
\end{proof}

\subsection{A Logarithmic Sobolev Inequality.}

In the remainder of this section we prove a logarithmic Sobolev
inequality for yet another birth death process.  Up to now we have
considered processes formed by considering marginal dynamics on a
single site of zero range.  These single site marginals are used
mainly in the proof of Theorem \ref{specgap}.  As we explain in
Section \ref{sec:specgap}, the proof of the spectral gap is also
an induction argument.  The argument differs in that with each
induction step we \emph{add} one site to the set $\Lambda$.  In
the proof of the logarithmic Sobolev inequality we \emph{double}
the size of $\Lambda$ at each induction step.  Hence, we will need
a logarithmic Sobolev inequality which acts on the number of
particles moving between two subsets of $\Lambda$, and not on a
single site.  In the remainder of this section we \emph{do} need
to use the assumption (E).

Let $\gamma_1(r_1)=\nu_{\Lambda,R}(R_1=r_1)$. Recall that $R_1$ is
the random variable which counts the total number of particles in
subset $\Lambda_1$.  The function $\gamma_1(r_1)$ is a probability
measure on $\{0,1,\hdots,r\}$ that is reversible for the birth and
death process with generator
\begin{eqnarray}
L^{bd}\psi(r_1)&=&\left[\frac{\gamma_1(r_1+1)}{\gamma_1(r_1)}\wedge1\right]\{\psi(r_1+1)-\psi(r_1)\}\notag\\
&&+\left[\frac{\gamma_1(r_1-1)}{\gamma_1(r_1)}\wedge1\right]\{\psi(r_1-1)-\psi(r_1)\}\label{gen:bd}
\end{eqnarray}
and Dirichlet form
\begin{eqnarray*}
D^{bd}(\psi)=\sum_{r_1=1}^r
[\gamma_1(r_1)\wedge\gamma_1(r_1-1)]\{\psi(r_1)-\psi(r_1-1)\}^2
\end{eqnarray*}
The work of Miclo in \cite{m1} allows us to check that this birth
and death process satisfies a logarithmic Sobolev inequality of
its own:

\begin{prop}\label{bdlogsob}
The birth death process defined through the generator
\reff{gen:bd} satisfies a logarithmic Sobolev inequality
$$H(\psi|\gamma(r_1)) \leq C r D^{bd}(\sqrt{\psi})$$
for some constant $C>0$, independent of the sites in $\Lambda_1$,
for all $\psi \geq 0$.
\end{prop}

We divide the proof into several steps.  In \cite{m1}, necessary
and sufficient conditions on the rates for birth and death
processes are given so that a logarithmic Sobolev inequality
holds.  In \cite{cmr} it is proved that the Miclo conditions are
satisfied by probability measures satisfying certain exponential
bounds.  Together these results imply the following statement:

\begin{lem}\label{miccond}
For a birth and death process as described above, suppose that
there exists a constant $A_0$ such that for any integer $r>0$ we
can find $\bar{r}\in\{0,\ldots,r\}$ such that $A_0^{-1}\bar{r}\leq
r-\bar{r}\leq A_0 \bar{r}$ and
\begin{align}
\frac{\gamma_1(r_1+1)}{\gamma_1(r_1)}\leq
e^{-\left(\frac{r_1-\bar{r}}{A_0 \bar{r}}\right)} \ \ \ \ \mbox{for $r_1 \in \{\bar r+1, \ldots, r\},$}\label{bd:cond1}\\
\frac{\gamma_1(r_1-1)}{\gamma_1(r_1)}\leq
e^{-\left(\frac{r_1-\bar{r}}{A_0 \bar{r}}\right)}\ \ \ \ \mbox{for $r_1 \in \{0, \ldots, \bar r-1\},$}\label{bd:cond2}\\
\frac{1}{A_0 \sqrt{\bar{r}}}e^{-\left(\frac{r_1-\bar{r}}{A_0
\bar{r}}\right)^2}\leq \gamma_1(r_1)\leq \frac{A_0}{
\sqrt{\bar{r}}}e^{-\left(\frac{r_1-\bar{r}}{A_0
\bar{r}}\right)^2}\ \ \ \ \mbox{for $r_1 \in \{0, \ldots,
r\}$}.\label{bd:cond3}
\end{align}
Then, there exists a positive constant $A_1$ such that for any
positive function $\psi$ on $\{0,\ldots, r\}$,
$$ H(\psi|\gamma_1(\cdot))\leq A_1 r \sum_{k=1}^r
[\gamma_1(r)\wedge\gamma_1(r)][\psi(r)-\psi(r-1)]^2,
$$
for any integer $r>0.$
\end{lem}

The proof of the above appears in \cite{cmr} or \cite{logsob1}.

If we can show that our probabilities $\gamma_1$ satisfy
\reff{bd:cond1} through \reff{bd:cond3} we will have proved
Proposition \ref{bdlogsob}. It turns out that these conditions are
satisfied by a modified measure, $\gamma_1^{\varepsilon_0}$, which
is equivalent to $\gamma_1$.  By the standard comparison method
the logarithmic Sobolev inequality then follows for $\gamma_1$. We
proceed by defining a class of equivalent measures
$\gamma_1^{\varepsilon}$, and then finding a particular value of
$\varepsilon$, $\varepsilon_0$ so that \reff{bd:cond1} through
\reff{bd:cond3} are satisfied.

We begin with a technical result.  Assume that $\Lambda$ is of
size $N$.

\begin{lem}\label{boundfromhell}
$$0 < \inf_{\Lambda, r>0} \sigma_\Lambda\left(\frac{r}{N}\right)\sqrt{N}\mu_{\Lambda,
\varphi_\Lambda(\frac{r}{N})}(R=r)\leq \sup_{\Lambda, r>0}
\sigma_\Lambda\left(\frac{r}{N}\right)\sqrt{N}\mu_{\Lambda,
\varphi_\Lambda(\frac{r}{N})}(R=r) < \infty$$
\end{lem}
\begin{proof}
This result follows for most cases directly from the local limit
theorems \ref{normalLLT-2} and \ref{poissonLLT-2}.  It remains to
bound the only case not covered by these theorems:  $N\leq N_0$
and $r>\rho_0 N$, for some fixed $N_0$ and $\rho_0$.  This last
case is simply condition (E).
\end{proof}

\begin{rem}\label{condLLR}
As mentioned before,  there are several simpler conditions such
that condition (E) is satisfied.  One we have mentioned already,
namely, that there exists a $K_0$ such that
\begin{eqnarray*}
c_x(k)=\theta k, \,\,\,\,\,\,\, \forall \,\, x \,\, \mbox{and}
\,\, k\geq K_0.
\end{eqnarray*}
To show that (E) is satisfied here is a straightforward, albeit
lengthy, calculation using Stirling's formula.  A similar argument
also shows that (E) is satisfied if we assume that there exists a
large constant $K_0$, and two positive constants $r_1$ and $r_2$
such that for all $k\geq K_0$ the rate function $c$ satisfies for
all $x$
\[ c_x(k)=\left \{ \begin{array} {ll}
                            \theta_1 \, k & \mbox{if $k$ is odd,}\\
                            \theta_2 \, k & \mbox{if $k$ is even.}
                            \end{array}
                   \right. \]
In fact, any other similar pattern also works.
\end{rem}

We next define the modified measures $\gamma_1^\varepsilon$ with
$\varepsilon$ in $(0,1/4)$ . We will first show that these
measures are equivalent to $\gamma_1$.  We will then show that
there is a special choice of $\varepsilon$ such that the measures
$\gamma_1^\varepsilon$ satisfy all three conditions of Lemma
\ref{miccond}.   The modified measure $\gamma_1^{\varepsilon}$ is
defined as follows. Let $I_\varepsilon:=[\varepsilon r,
(1-\varepsilon) r]\cap \bb Z$, and set $\bar r = \lceil
r/2\rceil$.  For $r_1$ in $\{0,\ldots, r\}$ we define the function
\begin{align*}
H(r_1)=\log \left\{\frac{\mu_{\Lambda_1,
\varphi_{\Lambda_1}(\frac{r_1}{|\Lambda_1|})}(R_1=r_1)\mu_{\Lambda_2,\varphi_{\Lambda_2}(\frac{r-r_1}{|\Lambda_2|})}(R_2=r-r_1)}{\mu_{\Lambda_1,
\varphi_{\Lambda_1}(\frac{\bar{r}}{|\Lambda_1|})}(R_1=r_1)\mu_{\Lambda_1,\varphi_{\Lambda_2}(\frac{\bar{r}}{|\Lambda_2|})}(R_2=r-r_1)}\right\},
\end{align*}
and use it to define the normalizing constant
\begin{eqnarray*}
Z=\frac{\sum_{k \in I^\varepsilon}e^{-H(k)}}{\sum_{k \in
I^\varepsilon}  \gamma_1(k)}.
\end{eqnarray*}
We may now define the new measure
\begin{align}\label{newmeasure}
\gamma_1^{\varepsilon}(r_1)=\left \{ \begin{array} {ll}
                            e^{-H(r_1)}/Z & \mbox{if $r_1 \in  I_\varepsilon$}\\
                            \gamma_1(r_1) & \mbox{otherwise.}
                            \end{array}
                   \right.
\end{align}
We next use Lemma \ref{boundfromhell} to show that the two
measures $\gamma_1^{\varepsilon}$ and $\gamma_1$ are equivalent.

\begin{lem}\label{lemma:equivmeas}
For any fixed $0<\varepsilon<\frac{1}{4}$ there exists a positive
constant $C$ such that
$$\frac{1}{C}\leq \frac{\gamma_1^\varepsilon(r_1)}{\gamma_1(r_1)} \leq C,$$
for all $r>0$, $\Lambda$ and $r_1 \in \{0,\ldots, r\}$.
\end{lem}

\begin{proof}
We need only check bounds inside $I_\varepsilon$.  Define
$$\pi(r_1)=\frac{\mu_{\Lambda_1,
\varphi_{\Lambda_1}(\frac{r_1}{|\Lambda_1|})}(R_1=r_1)\mu_{\Lambda_1,
\varphi_{\Lambda_2}(\frac{r-r_1}{|\Lambda_2|})}(R_2=r-r_1)}{\mu_{\Lambda,
\varphi_{\Lambda}(\frac{r}{|\Lambda|})}(R=r)}.$$

As we may write
$$\frac{\gamma_1(r_1)}{\gamma_1^\varepsilon(r_1)}=\frac{\sum_{k \in I_\varepsilon} \gamma_1(k) \frac{\pi(r_1)}{\pi(k)}}{\gamma_1(k)}$$
to prove the result it is enough to bound the ratio
$\frac{\pi(r_1)}{\pi(k)}$ uniformly as $k$ and $r_1$ vary over
$I_\varepsilon$.  Define
$$B(n,m)^2=\frac{\sigma^2_{\Lambda_1}(\frac{m}{|\Lambda_1|})\sigma^2_{\Lambda_2}(\frac{r-m}{|\Lambda_2|})}{\sigma^2_{\Lambda_1}(\frac{n}{|\Lambda_1|})\sigma^2_{\Lambda_2}(\frac{r-n}{|\Lambda_2|})}$$
By Lemma \ref{boundfromhell} we have that there exists a positive
finite constant $C_0$ such that
$$\frac{1}{C_0}B(r_1,\bar r)\leq \frac{\pi(r_1)}{\pi(k)}\leq C_0 B(r_1,\bar r).$$
To finish notice the following
$$\frac{1}{C_1}4\epsilon (1-\epsilon)\leq \frac{r_1(r-r_1)}{C_1\bar r(r-\bar r)}\leq B^2(r_1,\bar r)\leq \frac{C_1 r_1(r-r_1)}{\bar r(r-\bar r)}\leq \frac{C_1}{4\epsilon (1-\epsilon)}.$$
The result follows.
\end{proof}

Our next goal is to show that the measure
$\gamma_1^\epsilon(\cdot)$ satisfies conditions \reff{bd:cond1}
through \reff{bd:cond3}.   We first need the following.
\begin{prop}\label{gammabounds}
$$C^{-1}\frac{r_1}{r-r_1+1}\leq\frac{\gamma_1(r_1-1)}{\gamma_1(r_1)}\leq C\frac{r_1}{r-r_1+1}.$$
\end{prop}

\begin{proof}
We begin with the following identity:
\begin{eqnarray*}
\mu_{\Lambda,\varphi}[R=r+1]&=&\frac{1}{r+1}\sum_{x\in
\Lambda}\mu_{\Lambda,\varphi}[\eta_x\cdot\bb I[R=r+1]]\\
&=&\frac{1}{r+1}\sum_{x\in
\Lambda}\mu_{\Lambda,\varphi}[\frac{\varphi\cdot
(\eta_x+1)}{c_x(\eta_x+1)}\bb I[R=r+1]].
\end{eqnarray*}
Using the uniform bounds on $\frac{c_x(k)}{k}$ we obtain for some
$B>0$
$$\frac{\varphi |\Lambda| \mu_{\Lambda,\varphi}[R=r]}{B (r+1)} \leq \mu_{\Lambda,\varphi}[R=r+1] \leq \frac{B \varphi |\Lambda|\mu_{\Lambda,\varphi}[R=r]}{r+1}.$$
Together with
$$\frac{\gamma_i(r_i-1)}{\gamma_i(r_i)}=\frac{\mu_{\Lambda_i,\varphi}[R_i=r_i-1]\mu_{\Lambda_j,\varphi}[R_j=r-r_i+1]}{\mu_{\Lambda_i,\varphi}[R_i=r_i]\mu_{\Lambda_j,\varphi}[R_j=r-r_j]},$$
for $i\neq j$.  This implies the result.
\end{proof}

Proposition \ref{gammabounds} implies that for all $\varepsilon$
smaller than $1/3$ we have that
\begin{eqnarray*}
\frac{\gamma_1^\varepsilon(r_1+1)}{\gamma_1^\varepsilon(r_1)}&\leq\frac{1}{2}&\hspace{1cm}\mbox{ for $r_1 \in [(1-\varepsilon)r,r-1]$}\\
\frac{\gamma_1^\varepsilon(r_1-1)}{\gamma_1^\varepsilon(r_1)}&\leq\frac{1}{2}&\hspace{1cm}\mbox{ for $r_1 \in [1,\varepsilon r]$} \\
\end{eqnarray*}
for any $r>0$, and any $\Lambda$ with $|\Lambda|\geq 2$.  We next
prove the following lemma, which together with the above statement
implies that conditions \reff{bd:cond1} and \reff{bd:cond2} are
satisfied.
\begin{lem}\label{lemma:H}
For any $\varepsilon \in (0,1/4)$ there exists a positive constant
$A_0$ such that for any $r>0$
\begin{eqnarray*}
\frac{r_1-\bar r}{A_0 \bar r}\leq H(r_1+1)-H(r_1)\leq
A_0\frac{r_1-\bar r}{\bar r} & \hspace{1cm}\mbox{for any $r_1$ in
$I_\varepsilon$ and $r_1>\bar r$}\\
 \frac{r_1-\bar r}{A_0 \bar r}\leq
H(r_1-1)-H(r_1)\leq A_0\frac{r_1-\bar r}{\bar r} &
\hspace{1cm}\mbox{for any $r_1$ in $I_\varepsilon$ and $r_1<\bar
r$}
\end{eqnarray*}
\end{lem}
\begin{proof}
A careful calculation reveals
$$\partial_z H(z)= \log\varphi_{\Lambda_1}(z/m)-\log\varphi_{\Lambda_2}(r-z/m),$$
where $m=N/2=|\Lambda_1|$.  We differentiate again to obtain
$$\partial_z^2 H(z)=\frac{1}{\sigma_{\Lambda_1}^2(z/m)m}-\frac{1}{\sigma^2_{\Lambda_2}((r-z)/m)m}.$$
Using the bounds from Proposition \ref{sigmabounds} we have that
there exists a positive constant $B$ so that
$$\frac{1}{B}\left\{\frac{1}{z}+\frac{1}{r-z}\right\}\leq \partial_z^2 H(z) \leq B \left\{\frac{1}{z}+\frac{1}{r-z}\right\}$$
for all $z$ in $[\bar r, (1-\varepsilon)r]$, and thus we can find
a constant $B_2=B_2(\varepsilon)$ so that
$$\frac{1}{B_2}\left\{\frac{1}{\bar r}\right\}\leq \partial_z^2 H(z) \leq B_2 \left\{\frac{1}{\bar r}\right\}.$$
Integrating once from $\bar r$ to $z$ and then again from $z$ to
$z+1$ we obtain the first part of the inequality.  We repeat the
argument to obtain the other direction.
\end{proof}

It remains to prove the third condition.

\begin{lem}
There exists $\varepsilon_o \in (0,1/4)$ and a positive constant
$A_0$ such that
$$\frac{1}{A_0
\sqrt{\bar{r}}}e^{-\left(\frac{r-\bar{r}}{A_0 \bar{r}}\right)}
\leq \gamma_1^{\varepsilon_o}(r)\leq \frac{A_0}{
\sqrt{\bar{r}}}e^{-\left(\frac{r-\bar{r}}{A_0 \bar{r}}\right)}$$
\end{lem}

\begin{proof}
We split the proof into several steps. We will also make use of
the fact that the measures satisfy conditions \reff{bd:cond1} and
\reff{bd:cond2}.

\smallskip \noindent \textbf{Step 1.}  We write the arguments below for $r_1>\bar r$;
the argument in the opposite direction is the same.
\begin{eqnarray*}
\log\frac{\gamma_1^{\varepsilon}(r_1)}{\gamma_1^{\varepsilon}(\bar
r)}&=&\sum_{k=\bar r}^{r_1-1}
\log\frac{\gamma_1^{\varepsilon}(k+1)}{\gamma_1^{\varepsilon}(k)}\\
&\leq & -\sum_{k=\bar r}^{r_1-1} \frac{k-\bar r}{A_0 \bar r}\\
&\leq& - \frac{(r_1-\bar r)^2}{2A_0\bar r} +\frac{1}{2A_0}
\end{eqnarray*}
for any $r_1$.  This implies that there exists an $A_1>0$ so that
$$\frac{\gamma_1^{\varepsilon}(r_1)}{\gamma_1^{\varepsilon}(\bar
r)}\leq A_1 e^{-\frac{(r_1-\bar r)^2}{A_1 \bar r}}.$$

\smallskip
\noindent \textbf{Step 2.} We repeat a similar argument using
Lemma \ref{lemma:H} to obtain
$$\frac{\gamma_1^{\varepsilon}(r_1)}{\gamma_1^{\varepsilon}(\bar
r)}\geq A_2 e^{-\frac{(r_1-\bar r)^2}{A_2 \bar r}},$$ for some
$A_2>0$ and for any $r_1$ in $I_\varepsilon$.  We are restricted
to $I_\varepsilon$ as that is where we obtain the necessary lower
bound.

\smallskip
\noindent \textbf{Step 3.}  We now extend the above to $r_1$
outside of $I_\varepsilon$.  We work in one direction first,
assuming that $r_1> (1-\epsilon)r$.  Let $\tilde r =
\lfloor(1-\epsilon)r\rfloor$,
\begin{eqnarray*}
\frac{\gamma_1^{\varepsilon}(r_1)}{\gamma_1^{\varepsilon}(\bar
r)}&=&\frac{\gamma_1(r_1)}{\gamma_1(\tilde
r+1)}\frac{\gamma_1(\tilde r+1 )}{\gamma_1^{\varepsilon}(\tilde
r)}\frac{\gamma_1^{\varepsilon}(\tilde r
)}{\gamma_1^{\varepsilon}(\bar r)} \\
&\geq& C^{-1}(\varepsilon)\frac{\gamma_1(r_1)}{\gamma_1(\tilde
r)}\frac{\gamma_1^\varepsilon(\tilde r
)}{\gamma_1^{\varepsilon}(\bar r)},
\end{eqnarray*}
by Lemma \ref{lemma:equivmeas}.  For the first fraction we have
Proposition \ref{gammabounds} which gives
\begin{eqnarray*}
\log \frac{\gamma_1(r_1)}{\gamma_1(\tilde
r)}&\geq& \sum_{k=\tilde r}^{r_1-1} \log \frac{r-k}{C(k+1)}\\
&\geq&- r\log C - \frac{r}{\epsilon}.
\end{eqnarray*}
Combining this with the results of step 2, we have that
\begin{eqnarray*}
\frac{\gamma_1^{\varepsilon}(r_1)}{\gamma_1^{\varepsilon}(\bar
r)}&\geq&\frac{A_2}{C}e^{-\frac{(\tilde r-\bar r)^2}{A_2 \bar
r}}e^{-\log C r - r/\varepsilon}.
\end{eqnarray*}
We select an $\varepsilon=\varepsilon_0$ sufficiently small so
that we obtain a positive constant $C_2$ so that
\begin{eqnarray*}
\frac{\gamma_1^{\varepsilon_0}(r_1)}{\gamma_1^{\varepsilon_0}(\bar
r)}&\geq&\frac{A_2}{C}e^{-C_2 r}\\
&\geq& \frac{A_2}{C}e^{-64C\frac{(r_1-\bar r)^2}{\bar r}}.
\end{eqnarray*}

Repeating the argument in the opposite direction, we obtain that
there exists a positive $A_3$ such that for any $r_1$
\begin{eqnarray*}
\frac{1}{A_3}e^{-A_3 \frac{(r_1-\bar r)^2}{\bar r}}\leq
\frac{\gamma_1^{\varepsilon_o}(r_1)}{\gamma_1^{\varepsilon_o}(\bar
r)}\leq A_3 e^{-A_3 \frac{(r_1-\bar r)^2}{\bar r}}.
\end{eqnarray*}
Note that if repeating the argument produces a smaller
$\varepsilon_0$, we simply take the smaller of the two.

\smallskip
 \noindent
\textbf{Step 4.} We next sum the above fractions in $r_1$ to
obtain the following bound
$$\frac{1}{A_4\sqrt{\bar r}}\leq \gamma^{\varepsilon_0}_1(\bar r)\leq \frac{A_4}{\sqrt{\bar r}}$$
for some constant $A_4$.  This together with the previous bound
implies condition \reff{bd:cond3}.
\end{proof}

\section{Independence of the number of particles \label{sec:ind}}

We now have the tools necessary to proceed with the first part of
the proof of the main result, Theorem \ref{bigresult}.  In this
section we finish the argument given in the outline given in
Section \ref{sec:outline}, which allows us to establish that the
logarithmic constant is independent of $r$, the total number of
particles.  There are two things we need to do in order to
establish this result.  We first need to establish the initial
induction step from line \reff{inistep}.  Second, we need to
obtain the bound in \reff{impline}:
\begin{align*}
H(\nu_{\Lambda, R}(f|R_1=r_1)|\nu_{\Lambda,r})\leq
C(N)(1+\kappa(N,r))D_{\Lambda,r}(\sqrt{f})+C(N)\nu_{\Lambda,r}[f],
\end{align*}
where $C(N)$ is a large positive constant depending on $N$. We
begin with the latter, and establish the initial induction result
in Proposition \ref{iniinduct}. The proof of \reff{impline} is
computationally intensive; hence, for ease of reading, we split it
up into several steps.  In step 1, we reduce the problem to
calculating bounds on two covariances. These estimates are
provided in steps 2 and 3. In step 4 we combine these bounds to
obtain the above result.  All of the work involved is essentially
identical to that for the homogeneous problem.  The only
differences lie in that many functions we now estimate depend on
the site, and hence we need to check that all the bounds hold
uniformly. However, as most of these bounds rely on single site
estimates, the uniform bounds (LG) and (M) on the rates are
sufficient for the results to hold.

\medskip

\noindent \textbf{Step 1: Initial Calculations. } In Section
\ref{birthdeath} we defined a specific birth and death process,
and showed that this process satisfies a logarithmic Sobolev
inequality. We begin by applying this result to the term
$H(\nu_{\Lambda, R}[f|R_1=r_1]|\nu_{\Lambda,r})$.
\begin{eqnarray*}
H(\nu_{\Lambda,R}[f|R_1]|\nu_{\Lambda,r})&=&H(\psi(R_1)|\nu_{\Lambda,r})\\
&=&H(\psi|\gamma_1(\cdot))\\
&\leq & C r \sum_{r_1=1}^r
[\gamma_1(r_1)\wedge\gamma_1(r_1-1)](\sqrt{\psi(r_1)}-\sqrt{\psi(r_1-1)})^2.
\end{eqnarray*}
Recall that $\gamma_1(k)=\nu_{\Lambda,r}(R_1=k).$
 We continue
\begin{eqnarray}
&&\hspace{-1cm}H(\nu_{\Lambda,R}[f|R_1]|\nu_{\Lambda,r})\notag\\
&\leq&C r
\sum_{r_1=1}^r \gamma_1(r_1)\wedge\gamma_1(r_1-1)\left(\sqrt{\nu_{\Lambda,r}(f|R_1=r_1)}-\sqrt{\nu_{\Lambda,r}(f|R_1=r_1-1)}\right)^2\notag\\
&\leq& C r\sum_{r_1=1}^r
\left\{\rule[-0.3cm]{0cm}{1cm}\left[\frac{\gamma_1(r_1)\wedge\gamma_1(r_1-1)}{\nu_{\Lambda,r}(f|R_1=r_1)\vee\nu_{\Lambda,r}(f|R_1=r_1-1)}\right]\right.\notag\\
&&\left.\hspace{1.5in}\times(\nu_{\Lambda,r}(f|R_1=r_1)-\nu_{\Lambda,r}(f|R_1=r_1-1))^2\rule[-0.3cm]{0cm}{1cm}\right\}\label{line3}
\end{eqnarray}
where the inequality
$(\sqrt{a}-\sqrt{b})^2\leq\frac{(a-b)^2}{a\vee b}$, for $a$ and
$b$ positive, was used in the above.
\begin{prop}\label{above1}
For every $f$ and $r_1=1, 2, \hdots, r$ we have
\begin{align*}
&\nu_{\Lambda,r}(f\mid R_1=r_1)-\nu_{\Lambda,r}(f\mid R_1=r_1-1)\\
&=\frac{\gamma_1(r_1-1)}{\gamma_1(r_1)}\frac{1}{r_1N}\left[\nu_{\Lambda,r}\left(\left.\sum_{x\in
\Lambda_1, y\in \Lambda_2 }h_x(\eta_x)c_y(\eta_y)\nabla_{y,x}f\right| R_1=r_1-1\right)\right.\\
&\hspace{1cm}\left.+\nu_{\Lambda,r}\left(\left.f;\sum_{x\in
\Lambda_1, y\in \Lambda_2 }h_x(\eta_x)c_y(\eta_y)\right|
R_1=r_1-1\right)\right],
\end{align*}
where $h_x(k)=\frac{k+1}{c_x(k+1)}$.  Moreover, by exchanging the
roles of $\Lambda_1$ and $\Lambda_2$ we obtain for $r_1=0,1, 2,
\hdots, r-1$ (equivalently, $r_2=1, 2, \hdots, r$ where
$r_2=r-r_1$), where $\gamma_2(r_2)=\nu_{\Lambda,r}(R_2=r_2)$
\begin{align*}
&\nu_{\Lambda,r}(f\mid R_1=r_1)-\nu_{\Lambda,r}(f\mid R_1=r_1-1)\\
&=\frac{\gamma_2(r-r_1)}{\gamma_2(r-r_1+1)}\frac{1}{(r-r_1+1)N}\left[\nu_{\Lambda,r}\left(\left.\sum_{x\in
\Lambda_1, y\in \Lambda_2 }h_y(\eta_y)c_x(\eta_x)\nabla_{x,y}f\right| R_1=r_1\right)\right.\\
&\hspace{1cm}\left.+\nu_{\Lambda,r}\left(\left.f;\sum_{x\in
\Lambda_1, y\in \Lambda_2 }h_y(\eta_y)c_x(\eta_x)\right|
R_1=r_1\right)\right].
\end{align*}
\end{prop}
We define
\begin{eqnarray*}
A_{1}(r_1)&=&\nu_{\Lambda,r}\left(\left.\sum_{x\in
\Lambda_1, y\in \Lambda_2 }h_x(\eta_x)c_y(\eta_y)\nabla_{y,x}f\right| R_1=r_1\right),\\
B_{1}(r_1)&=&\nu_{\Lambda,r}\left(\left.f;\sum_{x\in \Lambda_1,
y\in \Lambda_2 }h_x(\eta_x)c_y(\eta_y)\right| R_1=r_1\right),
\end{eqnarray*}
and
\begin{eqnarray}
A(r_1)&=&\left\{\begin{array}{ll}
\frac{\gamma_1(r_1-1)}{\gamma_1(r_1)}\frac{1}{r_1N}\,\,A_1(r_1-1)&\mbox{for
$r_1> \frac{r}{2}$}\\
\frac{\gamma_2(r-r_1)}{\gamma_2(r-r_1+1)}\frac{1}{(r-r_1+1)N}\,\,A_1(r_1)&\mbox{for
$r_1\leq \frac{r}{2}$}
\end{array}\right.\label{defA}\\
B(r_1)&=&\left\{\begin{array}{ll}
\frac{\gamma_1(r_1-1)}{\gamma_1(r_1)}\frac{1}{r_1N}\,\,B_1(r_1-1)&\mbox{for $r_1> \frac{r}{2}$}\\
\frac{\gamma_2(r-r_1)}{\gamma_2(r-r_1+1)}\frac{1}{(r-r_1+1)N}\,\,B_1(r_1)&\mbox{for
$r_1\leq \frac{r}{2}.$}
\end{array}\label{defB}\right.
\end{eqnarray}
We combine both representations above to obtain
$$\nu_{\Lambda,r}(f\mid R_1=r_1)-\nu_{\Lambda,r}(f\mid R_1=r_1-1)=A(r_1)+B(r_1),$$
Our next steps will be to obtain bounds on the terms $A$ and $B$.

\begin{proof}[Proof of Proposition \ref{above1}]
The proof is based on the following two calculations. \linebreak
For $y\in \Lambda_2,$
\begin{align*}
&\hspace{-1cm}\nu_{\Lambda,r}[f\bb I (\eta_x >0)|R_1=r_1]\\
&=\left \{
\begin{array}{ll}
                            \frac{\gamma_1(r_1-1)}{\gamma_1(r_1)}\nu_{\Lambda,r}\left[\frac{c_y(\eta_y)}{c_x(\eta_x+1)}f(\eta^{y,x})|R_1=r_1-1\right] & \mbox{if $x \in  \Lambda_1$}\\
                            \nu_{\Lambda,r}\left[\frac{c_y(\eta_y)}{c_x(\eta_x+1)}f(\eta^{y,x})|R_1=r_1\right] & \mbox{if $x \in  \Lambda_2$}
                            \end{array}
                   \right..
\end{align*}
Next, write
\begin{eqnarray*}
\nu_{\Lambda,r}[f|R_1=r_1]=-\frac{1}{r_1}\sum_{x\in\Lambda_1}\nu_{\Lambda,r}[\eta_x
\cdot\nabla^{x,y}f|R_1=r_1]+\frac{1}{r_1}\sum_{x\in\Lambda_1}\nu_{\Lambda,r}[\eta_x
f(\eta^{x,y})|R_1=r_1],
\end{eqnarray*}
and plug in the above calculation, along with
$h_x(\eta_x)=\frac{\eta_x+1}{c_x(\eta_x+1)}$ to obtain
\begin{eqnarray*}
&&\hspace{-1cm}\nu_{\Lambda,r}[f|R_1=r_1]\\
&=&\frac{\gamma_1(r_1-1)}{\gamma_1(r_1)r_1}\left(\nu_{\Lambda,r}\left[\left.c_y(\eta_y)\sum_{x\in\Lambda_1}h_x(\eta_x)\cdot\nabla^{y,x}f\right|R_1=r_1-1\right]\right.\\
&&\hspace{3cm}+\hspace{1cm}\left.\nu_{\Lambda,r}\left[\left.c_y(\eta_y)\sum_{x\in\Lambda_1}h_x(\eta_x)f\right|R_1=r_1-1\right]\right).
\end{eqnarray*}
Setting $f=1$ in the above formula we obtain that
$$\frac{\gamma_1(r_1-1)}{\gamma_1(r_1)r_1}\nu_{\Lambda,r}\left[\left.c_y(\eta_y)\sum_{x\in\Lambda_1}h_x(\eta_x)\right|R_1=r_1-1\right]=1,$$
and hence we have
\begin{eqnarray*}
&&\hspace{-1cm}\nu_{\Lambda,r}[f|R_1=r_1-1]\\
&=&\frac{\gamma_1(r_1-1)}{\gamma_1(r_1)}\left(\nu_{\Lambda,r}\left[\left.c_y(\eta_y)\sum_{x\in\Lambda_1}h_x(\eta_x)f\right|R_1=r_1-1\right]\right.\\
&&\hspace{3cm}-\hspace{1cm}\left.\nu_{\Lambda,r}\left[\left.c_y(\eta_y)\sum_{x\in\Lambda_1}h_x(\eta_x);f\right|R_1=r_1-1\right]\right).
\end{eqnarray*}

\end{proof}

\noindent\textbf{Step 2: bounds on $\mb A$.}  Suppose
$\Lambda=\Lambda_1\cup\Lambda_2$.  For any $i=1,2$ define
$\gamma_i(r_i)=\nu_{\Lambda,r}[R_i=r_i]$.  Notice that by symmetry
Proposition \ref{gammabounds} applies also to $\gamma_2$.

\begin{prop}\label{prop:A}
Recall the definition of $A$ from \reff{defA}. There exists a
constant $C>0$ such that
\begin{align*}
A^2(r_1)&\leq C \frac{N^2}{r} \ \nu_{\Lambda,r}(f\mid
R_1=r_1)\vee\nu_{\Lambda,r}(f\mid R_1=r_1-1)\\
&\hspace{2cm}\times\left[\frac{\gamma_1(r_1-1)}{\gamma_1(r_1)}D_{\nu_{\Lambda,r}(\cdot\mid
R_1=r_1)}(\sqrt{f})+D_{\nu_{\Lambda,r}(\cdot\mid
R_1=r_1)}(\sqrt{f})\right].
\end{align*}
\end{prop}

\begin{proof}
We work out the case where $r_1>r/2$, as the argument in the other
direction is identical.  Because $h_x \leq c_2$ and $\nabla_{y,x}f
=\nabla_{y,x}\sqrt{f}\,[\sqrt{f}(\eta^{y,x})+\sqrt{f}(\eta)], $ we
may use the Cauchy-Schwarz inequality to bound $A^2(r_1)$ by
\begin{align}
&c_2^2\left[\frac{\gamma_1(r_1-1)}{\gamma_1(r_1)}\frac{1}{r_1N}\right]^2\nu_{\Lambda,r}\left(\sum_{x\in
\Lambda_1, y\in \Lambda_2 }c_y(\eta_y)(\nabla_{y,x}\sqrt{f})^2\mid
R_1=r_1-1\right)\notag\\
&\hspace{1cm}\times\nu_{\Lambda,r}\left(\sum_{x\in \Lambda_1, y\in
\Lambda_2 }c_y(\eta_y)(f(\eta^{y,x})+f(\eta))\mid
R_1=r_1-1\right)\notag
\end{align}
We next change measure to move from $f(\eta^{y,x})$ back to $f$
with $x\in\Lambda_1$ and $y\in\Lambda_2$ using
$$\nu_{\Lambda,r}[c_y(\eta_y)f(\eta^{y,x})|R_1=r_1-1]=\frac{\gamma_1(r_1)}{\gamma_1(r_1-1)}\nu_{\Lambda,r}[c_x(\eta_x)f(\eta)|R_1=r_1]$$
which follows from the detailed balance condition \reff{detbal}.
This gives us a new bound of
\begin{align}
&c_2^2\left[\frac{\gamma_1(r_1-1)}{\gamma_1(r_1)}\frac{1}{r_1N}\right]^2\nu_{\Lambda,r}\left(\sum_{x\in
\Lambda_1, y\in \Lambda_2 }c_y(\eta_y)(\nabla_{y,x}\sqrt{f})^2\mid
R_1=r_1-1\right)\notag\\
&\times\sum_{x\in \Lambda_1, y\in \Lambda_2
}\left\{\frac{\gamma_1(r_1)}{\gamma_1(r_1-1)}\nu_{\Lambda,r}\left(c_x(\eta_y)f(\eta)\mid
R_1=r_1\right)+\nu_{\Lambda,r}\left(c_y(\eta_y)f(\eta)\mid
R_1=r_1-1\right)\right\}.\notag
\end{align}
 We bound $\sum_{x\in
\Lambda_1, y\in \Lambda_2 }c_y(\eta_y)$ with $c_2 (r-r_1+1)N$,
$\sum_{x\in\Lambda_1, y\in\Lambda_2}c_x(\eta_x)$ with $c_2 r_1 N$,
and use Proposition \ref{gammabounds} to obtain the new bound for
$A^2(r_1)$
\begin{align}
&c_2^2\left[\frac{\gamma_1(r_1-1)}{\gamma_1(r_1)}\frac{1}{r_1N}\right]\nu_{\Lambda,r}\left(\sum_{x\in
\Lambda_1, y\in \Lambda_2 }c_y(\eta_y)(\nabla_{y,x}\sqrt{f})^2\mid
R_1=r_1-1\right)\label{above2}\\
&\hspace{1cm}\times \left\{c_2\,C
\nu_{\Lambda,r}[f|R_1=r_1]\vee\nu_{\Lambda,r}[f|R_1=r_1-1]
\right\}.\label{above2.1}
\end{align}
We next bound $(\nabla_{y,x}\sqrt{f})^2$ with
$$C' N \sum_{e_z}(\sqrt{f(\eta^{x,e_{z+1}})}-\sqrt{f(\eta^{x,e_z})})^2,$$
for some positive constant $C'$, where the sum is over $e_z$:
sites which form a path from $y$ to $x$. We use this along with
repeated change of measure to bound \reff{above2} as before.
\begin{eqnarray*}
&&\hspace{-1cm}\left[\frac{\gamma_1(r_1-1)}{\gamma_1(r_1)}\frac{1}{r_1N}\right]\nu_{\Lambda,r}\left[\left.\sum_{x\in
\Lambda_1, y\in \Lambda_2
}c_y(\eta_y)(\nabla_{y,x}\sqrt{f})^2\right|
R_1=r_1-1\right]\\
&\leq&\left[C'\frac{\gamma_1(r_1-1)}{\gamma_1(r_1)}\frac{1}{r_1}\right]\sum_{x\in
\Lambda_1, y\in \Lambda_2 }\left\{\sum_{e_z \in
\Lambda_1}\nu_{\Lambda,r}\left[\left.c_{e_z}(\eta_{e_z})(\nabla^{e_{z},e_{z+1}}f)^2\right|
R_1=r_1-1\right]\right.\\
&&\hspace{3cm}+\frac{\gamma_1(r_1)}{\gamma_1(r_1-1)}\left.\sum_{e_z
\in\Lambda_2}\nu_{\Lambda,r}\left[\left.c_{e_z}(\eta_{e_z})(\nabla^{e_{z},e_{z+1}}f)^2\right|
R_1=r_1\right]\right\}\\
&\leq&\frac{2c_2N^2}{r_1}\sum_{x\sim
y}\left\{\rule[-0.3cm]{0cm}{1cm}\frac{\gamma_1(r_1-1)}{\gamma_1(r_1)}\nu_{\Lambda,r}[c_x(\eta_x)(\nabla^{x,y}f)^2\mid
R_1=r_1-1]\right.\\
&&\hspace{5cm}\left.+\nu_{\Lambda,r}[c_x(\eta_{x})(\nabla^{x,y}f)^2\mid
R_1=r_1]\rule[-0.3cm]{0cm}{1cm}\right\}\\
\end{eqnarray*}
This last bound together with \reff{above2.1} completes the proof.
\end{proof}

\vspace{0.2in} \noindent\textbf{Step 3: bounds on $\mb B$.}

\begin{prop}\label{Bprop}If $r_1\leq r/2$ then
\begin{eqnarray}
B^2(r_1)&\leq&
C(N)\frac{\gamma_2^2(r_2)}{\gamma_2^2(r_2+1)}\nu_{\Lambda,r}[f|R_1=r_1]\left[\frac{r_1}{(r_2+1)^2}\nu_{{\Lambda_2,r-r_1}}[H(f|\nu_{{\Lambda_1,r_1}})]\right.\notag\\
&&+\left.\frac{r_1^2}{(r_2+1)^3}\left\{\nu_{\Lambda,r}[f|R_1=r_1]+\nu_{{\Lambda_1,r_1}}[H(f|\nu_{{\Lambda_2,r-r_1}})]\rule[-0.1cm]{0cm}{0.6cm}\right\}\right]\label{line:B3}.
\end{eqnarray}
If $r_1> r/2$ then
\begin{eqnarray*}
B^2(r_1)&\leq&
C(N)\frac{\gamma_1^2(r_1-1)}{\gamma_1^2(r_1)}\nu_{\Lambda,r}[f|R_1=r_1-1]\left[\frac{r_2+1}{r_1^2}\nu_{{\Lambda_2,r-r_1+1}}[H(f|\nu_{{\Lambda_1,r_1-1}})]\right.\notag\\
&&+\left.\frac{r_1^2}{(r_2+1)^3}\left\{\nu_{\Lambda,r}[f|R_1=r_1-1]+\nu_{{\Lambda_1,r_1-1}}[H(f|\nu_{{\Lambda_2,r-r_1+1}})]\rule[-0.1cm]{0cm}{0.6cm}\right\}\right].
\end{eqnarray*}
\end{prop}

We will prove here the case where $r_1\leq r/2$. Recall that in
this instance $B(r_1)$ is equal to
$$\frac{\gamma_2(r-r_1)}{\gamma_2(r-r_1+1)}\frac{1}{(r-r_1+1)N}\nu_{\Lambda,r}\left[\left.f;\sum_{x\in \Lambda_1,
y\in \Lambda_2 }h_y(\eta_y)c_x(\eta_x)\right| R_1=r_1\right].$$ We
use the fact that
$\nu_{\Lambda,r}(\cdot|R_1=r_1)=\nu_{\Lambda_1,r_1}\otimes\nu_{\Lambda_2,r-r_1}$
and that $\sum_x c_x$ and $\sum_yh_y$ act on $\eta_{\Lambda_1}$
and $\eta_{\Lambda_1}$ respectively. Hence,
\begin{eqnarray*}
&&\hspace{-2cm}\nu_{\Lambda,r}\left[\left.f;\sum_{x\in \Lambda_1,
y\in \Lambda_2
}h_y(\eta_y)c_x(\eta_x)\right| R_1=r_1\right]\\
&=&\nu_{\Lambda_2,r-r_1}\left[\sum_{y\in \Lambda_2
}h_y(\eta_y)\cdot\nu_{\Lambda_1,r_1}\left[f;\sum_{x\in
\Lambda_1}c_x(\eta_x)\right]\right]\\
&&\ \ \ +\ \
\nu_{\Lambda_2,r-r_1}\left[\nu_{\Lambda_1,r_1}[f];\sum_{y\in
\Lambda_2}h_y(\eta_y)\right]\nu_{\Lambda_1,r_1}\left[\sum_{x\in\Lambda_1}c_x(\eta_x)\right].
\end{eqnarray*}
Using the consequences of assumptions (LG) and (M) we thus get the
following simple bound
\begin{align}
B^2(r_1)\leq2\frac{\gamma_2^2(r-r_1)}{\gamma_2^2(r-r_1+1)}\left\{\frac{1}{c_1^2(r_2+1)^2}\nu_{\Lambda_1,r_1}\left[\nu_{\Lambda_2,r-r_1}[f];\sum_{x\in
\Lambda_1}c_x(\eta_x)\right]^2\right.\notag\\
\left.+\frac{r_1^2c_2^2}{N^2(r_2+1)^2}\nu_{\Lambda_2,r-r_1}\left[\nu_{\Lambda_1,r_1}[f];\sum_{y\in
\Lambda_2}h_y(\eta_y)\right]^2\right\}.\label{line4}
\end{align}
To obtain bounds on the two remaining covariances we will make use
of the entropy inequality
\begin{align}
|\mu(f;g)|\leq\frac{\mu(f)}{s}\log(\mu(e^{s(g-\mu(g))})\vee\mu(e^{-s(g-\mu(g)}))+\frac{1}{s}H(f|\mu),\label{entineq}
\end{align}
valid for any $s>0$, and we will optimize over $s$ after obtaining
bounds on $\mu(e^{s(g-\mu(g))})$ for $\mu$ and $g$ of interest. In
the proof we will make use of the quantity
$r_x=\nu_{\Lambda,r}(\eta_x)$ and of the spectral gap from Theorem
\ref{specgap} .

\begin{prop}\label{roughexpbounds}
For a subset $\Lambda'$ of size $|\Lambda|=N'$ there exists a
constant $C(N')$ such that for every $N'>0$ and $t \in [-1,1]$,
the following hold
\begin{align*}
\log\,\,\nu_{\Lambda',r'}\left[e^{t(c_x(\eta_x)-\nu_{\Lambda,r}[c_x(\eta_x)])}\right]\leq
C(N') \,\,r' \,\,t^2\\
\nu_{\Lambda',r'}\left[e^{t\cdot r \cdot
\{h_x(\eta_x)-\nu_{\Lambda,r}[h_x(\eta_x)]\}}\right]\leq
C(N')e^{C(N')\{r't^2+\sqrt{r'}|t|\}}.
\end{align*}
\end{prop}
\begin{proof}
To simplify notation slightly we will denote $\Lambda'$, $N'$ and
$r'$ simply as $\Lambda$, $N$ and $r$.  We start with the first
inequality. Notice that if $r$ is bounded then a simple Taylor
series expansion proves the bound. We may hence assume that $r$ is
larger than any finite constant we need. We begin with with the
identity
\begin{align}\label{changeofmeasure}
\nu_{\Lambda,r}[c_x(\eta_x)f(\eta)]=\nu_{\Lambda,r}[c_x(\eta_x)]\nu_{N,r-1}[f(\eta+\delta_x)].
\end{align}
Define $h(t)=\nu_{\Lambda,r}[e^{tc_x(\eta_x)}]$.  It follows that
\begin{eqnarray*}
t h'(t)-h(t)\log h(t)&=&t
\nu_{\Lambda,r}[c_x(\eta_x)e^{tc_x(\eta_x)}]-\nu_{\Lambda,r}[e^{tc_x(\eta_x)}]\log\nu_{\Lambda,r}[e^{tc_x(\eta_x)}]\\
&\leq&t
\nu_{\Lambda,r}[c_x(\eta_x)e^{tc_x(\eta_x)}]-t\nu_{\Lambda,r}[c_x(\eta_x)]\nu_{\Lambda,r}[e^{tc_x(\eta_x)}]\\
&=&t\nu_{\Lambda,r}[c_x(\eta_x)]\left\{\nu_{N,r-1}[e^{tc_x(\eta_x+1)}]-\nu_{\Lambda,r}[e^{tc_x(\eta_x)}]\right\}\\
&\leq&c_2 t r_x
\left\{\nu_{N,r-1}[e^{tc_x(\eta_x+1)}]-\nu_{\Lambda,r}[e^{tc_x(\eta_x)}]\right\}.
\end{eqnarray*}
We next bound
$\nu_{N,r-1}[e^{tc_x(\eta_x+1)}]-\nu_{\Lambda,r}[e^{tc_x(\eta_x)}]$.
We split this bound into two parts.  By the inequality
$|e^x-e^y|\leq|x-y|e^{|x-y|}e^x$ we have
$$|\nu_{N,r-1}[e^{tc_x(\eta_x+1)}]-\nu_{N,r-1}[e^{tc_x(\eta_x)}]|\leq a_1te^{a_1t}\nu_{N,r-1}[e^{tc_x(\eta_x+1)}].$$
Next consider
\begin{eqnarray}
&&\hspace{-2cm}\left|\nu_{N,r-1}[e^{tc_x(\eta_x)}]-\nu_{\Lambda,r}[e^{tc_x(\eta_x)}]\right|\notag\\
&\hspace{-1cm}\leq&\left|\nu_{N,r-1}[e^{tc_x(\eta_x)}]-\nu_{N,r-M}[e^{tc_x(\eta_x)}]\right|+\left|\nu_{\Lambda,r}[e^{tc_x(\eta_x)}]-\nu_{N,r-M}[e^{tc_x(\eta_x)}]\right|\label{Qline}
\end{eqnarray}
The two pieces are now dealt with in the same way, by Lemma
\ref{monotone2} as long as $M$ is large enough, there exists a
coupling measure $Q$ on $\{\eta,\xi\}$ such that the marginal of
$\eta$ is $\nu_{\Lambda,r}$, the marginal of $\xi$ is
$\nu_{N,r-M}$, and $Q$ is concentrated on the configurations such
that $\xi\leq\eta$. Hence we have
\begin{eqnarray*}
&&\hspace{-2cm}|\nu_{\Lambda,r}[e^{tc_x(\eta_x)}]-\nu_{N,r-M}[e^{tc_x(\eta_x)}]| \\
&=&|Q[e^{tc_x(\eta_x)}-e^{tc_x(\xi_x)}]|\\
&\leq& t C' \nu_{\Lambda,r}[e^{tc_x(\eta_x)}],
\end{eqnarray*}
where $C'$ is some constant depending on $M$.  For the second term
in \reff{Qline} we obtain a bound of $t C'
\nu_{N,r-1}[e^{tc_x(\eta_x)}]$. We next replace
$\nu_{N,r-1}[e^{tc_x(\eta_x)}]$ with
$\nu_{\Lambda,r}[e^{tc_x(\eta_x)}]$.
\begin{eqnarray*}
\nu_{N,r-1}[e^{tc_x(\eta_x)}]&\leq&
e^{c_2t}\nu_{N,r-1}[e^{tc_x(\eta_x+1)}]\\
&=&\frac{e^{c_2t}}{\nu_{N,r-1}[c_x(\eta_x)]}\nu_{\Lambda,r}[c_x(\eta_x)e^{tc_x(\eta_x)}]\\
&\leq&\frac{c_2r}{c_1r_x}\nu_{\Lambda,r}[e^{tc_x(\eta_x)}]
\end{eqnarray*}
Putting all of this together, and recalling that $|t|\leq 1$, we
have
\begin{eqnarray*}
t h'(t)-h(t)\log h(t)&\leq&C''rt^2 h(t),
\end{eqnarray*}
for some constant $C''$, which is the same as
$$\frac{d}{dt} \frac{\log h(t)}{t}\leq C'' r .$$
Integrating in $t$ and noting that $\lim_{t\rightarrow
0}\frac{\log h(t)}{t}=\nu_{\Lambda,r}[c_x(\eta_x)]$ we obtain
$$\frac{\log h(t)}{t}\leq \nu_{\Lambda,r}[c_x(\eta_x)] + C''r t, $$
from which the first part of the Proposition follows.  The
constant $C''$ depends on $N$ via $M=BN$ of Lemma \ref{monotone2}.

\bigskip

We now turn to the proof of $$ \nu_{\Lambda,r}\left[e^{t\cdot r
\cdot \{h_x(\eta_x)-\nu_{\Lambda,r}[h_x(\eta_x)]\}}\right]\leq
C(N)e^{C(N)\{rt^2+\sqrt{r}|t|\}}.$$ By the change of measure
formula \reff{changeofmeasure} we may calculate the expectation of
\linebreak $h_x(\eta_x)=\frac{\eta_x+1}{c_x(\eta_x+1)}$, to be
$$\nu_{\Lambda,r}[h_x(\eta_x)]=\frac{\nu_{\Lambda,r+1}[\eta_x]}{\nu_{\Lambda,r+1}[c_x(\eta_x)]}=\frac{\tilde r_x}{\nu_{\Lambda,r+1}[c_x(\eta_x)]}.$$
Hence,
\begin{eqnarray*}
|h_x(\eta_x)-\nu_{\Lambda,r}[h_x(\eta_x)]| &\leq&
\frac{1}{c_x(\eta_x+1)\nu_{\Lambda,r+1}[c_x(\eta_x)]}\left\{\rule[-0.2cm]{0cm}{0.7cm}c_x(\eta_x+1)\left|\rule[-0.2cm]{0cm}{0.7cm}\eta_x+1-\tilde
r_x\right|\right.\\
&&\hspace{3cm}+\left.(\eta_x+1)\left|c_x(\eta_x+1)-\nu_{\Lambda,r+1}[c_x(\eta_x)]\rule[-0.2cm]{0cm}{0.7cm}\right|\rule[-0.2cm]{0cm}{0.7cm}\right\}\\
&\leq&\frac{C'}{\tilde
r_x}\left\{|c_x(\eta_x+1)-\nu_{\Lambda,r+1}[c_x(\eta_x)]|+|\eta_x+1-\tilde
r_x|\right\},
\end{eqnarray*}
for some constant $C'$ depending on $c_1$ and $c_2$.  Using the
uniform condtions (LG) and (M) we can show that for some $C>0$ we
have $$|\eta_x+1-\tilde r_x|\leq
C|c_x(\eta_x)-\nu_{\Lambda,r+1}[c_x(\eta_x)]|+C\sqrt{r},$$ which
implies that
\begin{eqnarray*}
|h_x(\eta_x)-\nu_{\Lambda,r}[h_x(\eta_x)]|&\leq& \frac{C}{\tilde
r_x} \left\{ C'
|c_x(\eta_x+1)-\nu_{\Lambda,r+1}[c_x(\eta_x)]|+C\sqrt{\tilde
r_x}\right\}\\
&\leq&\frac{C'C }{ r_x} \left\{
|c_x(\eta_x)-\nu_{\Lambda,r}[c_x(\eta_x)]|+\sqrt{ r}\right\}.
\end{eqnarray*}
We thus have that if
$|c_x(\eta_x)-\nu_{\Lambda,r}[c_x(\eta_x)]|\leq M$ then
$r|h_x(\eta_x)-\nu_{\Lambda,r}[h_x(\eta_x)]|\leq
\frac{Cr}{r_x}(M+\sqrt{r})\leq C(M+\sqrt{r})$.  Therefore, for $t$
in $(0,1]$
\begin{eqnarray}
\nu_{\Lambda,r}[e^{tr(h_x(\eta_x)-\nu_{\Lambda,r}[h_x(\eta_x)])}]&=&t\int
e^{tz}\nu_{\Lambda,r}[r_x(h_x(\eta_x)-\nu_{\Lambda,r}[h_x(\eta_x)])>z]dz\notag\\
&\leq&t\int
e^{tz}\nu_{\Lambda,r}[c_x(\eta_x)-\nu_{\Lambda,r}[c_x(\eta_x)]>\frac{z}{C}-C\sqrt{r}]dz\notag\\
&&+t\int
e^{tz}\nu_{\Lambda,r}[c_x(\eta_x)-\nu_{\Lambda,r}[c_x(\eta_x)]<-\frac{z}{C}+C\sqrt{r}]dz\notag\\\label{intlinabove}
\end{eqnarray}
By change of variable \reff{intlinabove} is equal to
\begin{eqnarray*} &&\hspace{-2cm}Ct\int
e^{t(CM+C\sqrt{r})}\nu_{\Lambda,r}[c_x(\eta_x)-\nu_{\Lambda,r}[c_x(\eta_x)]>M]dM\\
&&\hspace{-1.4cm}+\ Ct\int
e^{t(CM+C\sqrt{r})}\nu_{\Lambda,r}[c_x(\eta_x)-\nu_{\Lambda,r}[c_x(\eta_x)]<-M]dM\\
&=&Ce^{Ct\sqrt{r}}\left\{\nu_{\Lambda,r}[e^{Ct(c_x(\eta_x)-\nu_{\Lambda,r}[c_x(\eta_x)])}]\right.\\
&&\hspace{4cm}\left.+\nu_{\Lambda,r}[e^{-Ct(c_x(\eta_x)-\nu_{\Lambda,r}[c_x(\eta_x)])}]\right\}\\
&\leq & Ce^{Ct\sqrt{r}}e^{Crt^2}.
\end{eqnarray*}
Replacing $h_x(\eta_x)-\nu_{\Lambda,r}[h_x(\eta_x)]$ with its
negative gives us the same bound for $t$ in $[-1,0)$.
\end{proof}

Because the constant $C(N')$ is allowed to depend on $N'$ in any
way, we may extend these calculations of Proposition
\ref{roughexpbounds} using Cauchy-Schwarz to say
\begin{align}
\log\,\,\nu_{\Lambda',r'}\left(e^{t\cdot \sum_{x\in
\Lambda'}(c_x(\eta_x)-\nu_{\Lambda,r}(c_x(\eta_x)))}\right)\leq
C(N') \,\,r' \,\,t^2\label{abovebound c}\\
\nu_{\Lambda',r'}\left(e^{t\cdot r \cdot \sum_{x\in \Lambda'}
(h_x(\eta_x)-\nu_{\Lambda,r}(h_x(\eta_x)))}\right)\leq
C(N')e^{C(N')(r't^2+\sqrt{r'}|t|)}.\label{abovebound h}
\end{align}

\begin{prop}\label{rough1}
There exists a constant $C(N')$ (depending on $N'=|\Lambda'|$)
such that
\begin{align*}
\nu_{\Lambda',r'}\left[f;\sum_{x\in
\Lambda'}c_x(\eta_x)\right]^2\leq C(N')\cdot r'\cdot\nu_{\Lambda',r'}[f]H(\tilde f|\nu_{\Lambda',r'})\\
\nu_{\Lambda',r'}\left[f;\sum_{x\in
\Lambda'}h_x(\eta_x)\right]^2\leq
\frac{C(N')}{r'}\nu_{\Lambda',r'}[f]\left\{\nu_{\Lambda',r'}[f]+H(\tilde
f|\nu_{\Lambda',r'})\right\}
\end{align*}
\end{prop}

\begin{proof}[Sketch of proof]
The proposition follows from direct calculation if in each case we
insert the bounds \reff{abovebound c} and \reff{abovebound h} into
the entropy inequality \reff{entineq} and optimize over $s$.
\end{proof}

From the above bounds we have the following
\begin{eqnarray}
\nu_{\Lambda_1,r_1}\left[\nu_{\Lambda_2,r-r_1}[f];\sum_{x\in
\Lambda_1}c_x(\eta_x)\right]^2\leq C(N) r_1
\nu_{\Lambda,r}[f|R_1=r_1]H(\nu_{\Lambda_2,r-r_1}[f]|\nu_{{\Lambda_1,r_1}})\notag\\\label{line:B1}
\end{eqnarray}
as well as
\begin{align}
&\nu_{\Lambda_2,r-r_1}\left[\nu_{\Lambda_1,r_1}[f];\sum_{y\in
\Lambda_2}h_y(\eta_y)\right]^2\notag\\
&\leq \frac{C(N)} {r-r_1}
\nu_{\Lambda,r}[f|R_1=r_1]\left\{\nu_{\Lambda,r}[f|R_1=r_1]+H(\nu_{\Lambda_1,r_1}[f]|\nu_{{\Lambda_2,r-r_1}})\rule[-0.1cm]{0cm}{0.6cm}\right\}\label{line:B2}
\end{align}

We now insert \reff{line:B1}  and \reff{line:B2} into \reff{line4}
to obtain \reff{line:B3}.  Notice that we have also used that the
entropy is convex.

\vskip0.5cm

\noindent\textbf{Step 4: putting it all together. } We next use
Proposition \ref{gammabounds} which by symmetry also applies to
$\gamma_2$,
from which it follows that
\begin{align*}
\frac{\gamma_2^2(r-r_1)}{\gamma_2^2(r-r_1+1)}\frac{r_1^2}{(r_2+1)^3}\leq\frac{C}{r_2+1}.
\end{align*}
We insert this into \reff{line:B3} to obtain
\begin{align}
&Cr\,\,\frac{\gamma_1(r_1)\wedge\gamma_1(r_1-1)}{\nu_{\Lambda,r}(f|R_1=r_1)\vee\nu_{\Lambda,r}(f|R_1=r_1-1)}\,\,B^2(r_1)\notag\\
&= Cr\,\,\frac{\gamma_2(r_2)\wedge\gamma_2(r_2+1)}{\nu_{\Lambda,r}(f|R_1=r_1)\vee\nu_{\Lambda,r}(f|R_1=r_1-1)}\,\,B^2(r_1)\notag\\
&\leq
C(N)\gamma_1(r_1)\left(\nu_{\Lambda,r}(f|R_1=r_1)\right.\notag\\
&\left.\hspace{1cm}+\nu_{{\Lambda_2,r-r_1}}(H(f|\nu_{{\Lambda_1,r_1}}))+\nu_{{\Lambda_1,r_1}}(H(f|\nu_{{\Lambda_2,r-r_1}}))\right)\label{line:B4}
\end{align}
where we have also used that $r_1\leq\frac{r}{2}$.  We obtain a
similar answer in the case $r_1>\frac{r}{2}$.

This gives us the necessary bounds on the term $B$.  We combine
this with Proposition \ref{prop:A} which gives us bounds on the
term $A$, and insert into \reff{line3}:
\begin{align*}
&H(\nu_{\Lambda,R}(f|R_1)|\nu_{\Lambda,r})\notag\\
&\leq C r\sum_{r_1=1}^r
\frac{\gamma(r_1)\wedge\gamma(r_1-1)}{\nu_{\Lambda,r}(f|R_1=r_1)\vee\nu_{\Lambda,r}(f|R_1=r_1-1)}[A^2(r_1)+B^2(r_1)]\\
&\leq
C(N)\sum_{r_1=1}^r\gamma(r_1-1)\left\{ D_{\Lambda,r}(\sqrt{f})+\nu_{\Lambda,r}(f|r_1)\right.\notag\\
&\hspace{1in}\left.+\nu_{{\Lambda_2,r-r_1}}[H(f|\nu_{{\Lambda_1,r_1}})]+\nu_{{\Lambda_1,r_1}}[H(f|\nu_{{\Lambda_2,r-r_1}})]\rule[-0.1cm]{0cm}{0.5cm}\right\}
\end{align*}
for some possibly different constant $C(N)$ depending again only
on $N$.  We next apply the induction hypothesis to obtain the new
bound
\begin{align*}
C(N)[D_{\Lambda,r}(\sqrt{f})+\nu_{\Lambda,r}(f)+\kappa(N,r)D_{\Lambda,r}(\sqrt{f})]
\end{align*}
where $\kappa(N,r)$ was defined in \reff{def:kappa}.  This
completes the argument required to prove \reff{impline}.

\begin{prop}\label{iniinduct}
$$\sup_{r}\kappa(2,r)<\infty.$$
\end{prop}
\begin{proof}
We assume that $\Lambda=\{0,1\}$.  In this case, since there is a
total of $r$ particles, the function $f(\eta)=f(k,r-k)=\tilde
f(k)$.  We begin by calculating the Dirichlet form
\begin{eqnarray*}
D_{2,r}\left(\sqrt{\tilde f}\right)&=&\sum_{k=0}^r
\gamma_1(k)c_1(k)\lf[\sqrt{\tilde f(k-1)}-\sqrt{\tilde f(k)}\rg]\\
&&+\sum_{k=0}^r
\gamma_1(k)c_2(r-k)\lf[\sqrt{\tilde f(k+1)}-\sqrt{\tilde f(k)}\rg]\\
&=&\sum_{k=1}^r \gamma_1(k)c_1(k)\lf[\sqrt{\tilde
f(k-1)}-\sqrt{\tilde f(k)}\rg]
\end{eqnarray*}
using the relationship
$\gamma_1(k)c_2(r-k)=\gamma_1(k+1)c_1(k+1)$. We next prove that
there exists a finite constant $B$ so that $r\,\,
\gamma_1(k)\wedge\gamma_1(k-1)\leq B \gamma_1(k)c_1(k).$
\begin{eqnarray*}
\frac{\gamma_1(k)c_1(k)}{\gamma_1(k)\wedge\gamma_1(k-1)}=c_1(k)\vee
c_2(r-k-1)\geq \frac{1}{B} \frac{r}{2}.
\end{eqnarray*}
We may now put these results together to obtain
\begin{eqnarray*}
H(f|\nu_{2,r})&=&H(\tilde f|\gamma_1(\cdot))\\
&\leq&C r \sum_{k=1}^r
\gamma_1(k)\vee\gamma_1(k-1)\lf[\sqrt{\tilde
f(k-1)}-\sqrt{\tilde f(k)}\rg]\\
&\leq&C B \sum_{k=1}^r \gamma_1(k)c_1(k)\lf[\sqrt{\tilde
f(k-1)}-\sqrt{\tilde f(k)}\rg]\\
&=&CB D_{2,r}(\sqrt{f}).
\end{eqnarray*}
We also used Proposition \ref{bdlogsob} in the above.  This
completes the proof.
\end{proof}

\section{Tightening the Bounds\label{sec:tight}}

As discussed earlier, in this section we obtain improved bounds on
the covariances appearing in Proposition \ref{rough1}, which will
allow us to conclude that \reff{impline2} holds for large values
of $N$:
\begin{align*}
H(\nu_{\Lambda, R}(f|R_1=r_1)|\nu_{\Lambda,r})\leq C N^2
D_{\Lambda,r}(\sqrt{f})+C\nu_{\Lambda,r}[f] +
\kappa(N,r)D_{\Lambda,r}(\sqrt{f}).
\end{align*}
The two tighter bounds on the covariances in Proposition
\ref{rough1} are given below.
\begin{prop}  \label{tightc}
For every $\epsilon>0$ there exists a constant $C=C(\epsilon)>0$
and an $N_0=N_0(\epsilon)$ such that for all $|\Lambda'|=N\geq
N_0$, all $r'$ and all positive functions $f$
\begin{eqnarray*}
\nu_{\Lambda',r'}[f;\sum_{x\in\Lambda'}c_x(\eta_x)]^2\leq r'
\nu_{\Lambda',r'}[f]\left[C \nu_{\Lambda',r'}[f]+ C N^2
D_{\Lambda',r}(\sqrt{f})+\epsilon H(f|\nu_{\Lambda',r})\right].
\end{eqnarray*}
\end{prop}

\begin{prop}\label{tighth}
For every $\epsilon>0$ there exists a constant $C=C(\epsilon)>0$
such that for all $|\Lambda'|=N\geq N_0$, all $r'$ and all
positive functions $f$
\begin{eqnarray*}
\nu_{\Lambda',r'}[f;\sum_{x\in\Lambda'}h_x(\eta_x)]^2\leq
\frac{N^2}{r'} \nu_{\Lambda',r'}[f]\left[C \nu_{\Lambda',r'}[f]+ C
N^2 D_{\Lambda',r'}(\sqrt{f})+\epsilon
H(f|\nu_{\Lambda',r'})\right].
\end{eqnarray*}
\end{prop}

Let us first show how these bounds give the desired result. We may
assume that both $N_0$ are the same in the propositions. Using the
same argument as in Step 4 of the previous section, but replacing
the covariance bounds of Proposition \ref{rough1} with the bounds
of the above Propositions, we estimate $B^2(r_1)$. The new
estimates give
\begin{align*}
&r\,\,\frac{\gamma_1(r_1)\wedge\gamma_1(r_1-1)}{\nu_{\Lambda,r}(f|R_1=r_1)\vee\nu_{\Lambda,r}(f|R_1=r_1-1)}\,\,B^2(r_1)\\
&\leq
\gamma_1(r_1)\left(C(\epsilon)\nu_{\Lambda,r}(f|R_1=r_1)+C(\epsilon)N^2D_{\nu_{{\Lambda_2,r-r_1}}}(\sqrt{f})+\epsilon
H(f|\nu_{{\Lambda_2,r-r_1}})\right).
\end{align*}

We now combine the above result along with Proposition
\ref{prop:A} to continue with line \reff{line3}:
\begin{align*}
&H(\nu_{\Lambda,R}(f|R_1)|\nu_{\Lambda,r})\notag\\
&\leq C r\sum_{r_1=1}^r
\frac{\gamma_1(r_1)\wedge\gamma_1(r_1-1)}{\nu_{\Lambda,r}(f|R_1=r_1)\vee\nu_{\Lambda,r}(f|R_1=r_1-1)}[A^2(r_1)+B^2(r_1)]\\
&\leq
C\sum_{r_1=1}^r\gamma_1(r_1)\left\{N^2D_{\Lambda,r}(\sqrt{f})+C(\epsilon)\nu_{\Lambda,r}(f|r_1)\right.\notag\\
&\hspace{1in}\left.+C(\epsilon)N^2D_{\nu_{{\Lambda_2,r-r_1}}}(\sqrt{f})+\epsilon
H(f|\nu_{{\Lambda_2,r-r_1}})\right\}\\
&\leq C(\epsilon)N^2
D_{\Lambda,r}(\sqrt{f})+C(\epsilon)\nu_{\Lambda,r}(f)+C\cdot\epsilon
H(f|\nu_{{\Lambda_2,r-r_1}})
\end{align*}
We now fix an $\epsilon$ so that $C\cdot\epsilon<1$.  The above
together with \reff{def:kappa} gives the bound
\begin{align*}
C N^2
D_{\Lambda,r}(\sqrt{f})+C\nu_{\Lambda,r}(f)+\kappa(N,r)D_{\Lambda,r}(\sqrt{f}),
\end{align*}
for some new constant $C$, and $|\Lambda|\geq 2N_0$. This is
\reff{impline2} as required.

The rest of this section is divided as follows.  We describe in
detail the proof of Proposition \ref{tightc}, which is split into
two main cases: small and large density.  We then proceed with the
proof of Proposition \ref{tighth}, which follows by a similar
argument.  For ease of presentation, we will write $\Lambda$ for
$\Lambda'$ and $r$ for $r'$ in both proofs.

\subsection{Proof of Proposition \ref{tightc}} The proof of this
result is split into several lemmas.  We begin by partitioning
$\Lambda$ into $m$ disjoint blocks $\Lambda_1, \ldots, \Lambda_m$,
which we assume, without loss of generality, to be of equal size
$l=N/m$.  Denote by $\mathcal{G}$ the $\sigma$-field generated by
$R_1, \ldots, R_m$, where $R_i$ is the random number of particles
inside the subset $\Lambda_i$.  We thus obtain
\begin{align}
&\nu_{\Lambda, r}[f;\sum_{x\in\Lambda}c_x(\eta_x)]\notag\\
&= \nu_{\Lambda, r}[\nu_{\Lambda,
r}[f;\sum_{x\in\Lambda}c_x(\eta_x)|\mathcal{G}]]+\nu_{\Lambda,
r}[f;\sum_{k=1}^m\nu_{\Lambda_k,
R_k}[\sum_{x\in\Lambda_k}c_x(\eta_x)]]\label{tight1}
\end{align}
and we bound the left hand side and the right hand side of
\reff{tight1} separately.  The bound on the left hand side is
easier and its proof is essentially a restatement of the proof of
the first part of Proposition \ref{rough1} .
\begin{prop}\label{prop:tightc1}
There is a constant C, possibly depending on $l$, such that
$$\nu_{\Lambda, r}[\nu_{\Lambda,
r}[f;\sum_{x\in\Lambda}c_x(\eta_x)|\mathcal{G}]]^2\leq C r
\nu_{\Lambda, r}[f] \nu_{\Lambda,
r}[H(f|\nu_{\Lambda,r}(\cdot|\mathcal{G}))] .$$
\end{prop}

\begin{proof}
We begin with the entropy inequality; for any $t>0$ we have
\begin{eqnarray*}
\nu_{\Lambda,r}[f;\sum_{x \in \Lambda}c_x(\eta_x)|\mc G]&\leq&
\frac{\nu_{\Lambda,r}[f]}{t}\sum_{k=1}^m \log
\nu_{\Lambda_k,R_k}[\exp\{t \sum_{x \in \Lambda_k}
\left(c_x(\eta_x)-\nu_{\Lambda_k,R_k}[c_x( \eta_x)]\right)\}]\\
&&\hspace{1.5in}+\frac{1}{t}H(f|\nu_{\Lambda,r}[\cdot|\mc G]).
\end{eqnarray*}
Using the Cauchy-Schwarz inequality and Proposition
\ref{roughexpbounds} we have the following bound
\begin{eqnarray*}
\nu_{\Lambda_k,R_k}[\exp\{t \sum_{x \in \Lambda_k}
\left(c_x(\eta_x)-\nu_{\Lambda_k,R_k}[c_x(\eta_x)]\right)\}]\leq
\exp\{c(l) R_k t^2\},
\end{eqnarray*}
for some constant $c(l)$ depending on $l$. Combining the two
inequalities we then have for any $t>0$
\begin{eqnarray*}
\nu_{\Lambda, r}[\nu_{\Lambda,
r}[f;\sum_{x\in\Lambda}c_x(\eta_x)|\mathcal{G}]]^2\leq
\nu_{\Lambda,r}[f]^2 c(l)r^2 t^2 +
\frac{1}{t^2}\left(\nu_{\Lambda,r}[H(f|\nu_{\Lambda,r}(\cdot|\mc G
))]\right)^2
\end{eqnarray*}
The result follows if we optimize in $t$.
\end{proof}

The bounds on the right hand side of \reff{tight1} are
considerably more difficult.  These are given in the following
lemma.
\begin{prop}\label{tightc2}
For every $\epsilon>0$ there is an $l=l(\epsilon)$,
$N_0=N_0(\epsilon)$, and a constant $C=C(\epsilon)>0$ such that
for all $N\geq N_0$
\begin{align*}
\nu_{\Lambda, r}[f;\sum_{k=1}^m\nu_{\Lambda_k,
R_k}[\sum_{x\in\Lambda_k}c_x(\eta_x)]]^2\leq r
\nu_{\Lambda,r}[f]\left\{C\nu_{\Lambda,r}[f]+CN^2D_{\Lambda,r}(\sqrt{f})+\epsilon
H(f|\nu_{|\Lambda,r}) \right\}.
\end{align*}
\end{prop}

Notice that from Section \ref{sec:ind} we know that logarithmic
Sobolev constant $\kappa$ depends only the the size of the subset
(and not on the number of particles).  We apply this to obtain the
bound
\begin{eqnarray*}
\nu_{\Lambda,r}[H(f|\nu_{\Lambda,r}(\cdot|\mathcal{G}))]&=&\nu_{\Lambda,r}[\sum_{k=1}^m
H(f|\nu_{C_k,R_k}(\cdot)) ]\\
&\leq& \nu_{\Lambda,r}[\sum_{k=1}^m \kappa(l) D_{C_k,R_k}(\sqrt{f})]\\
&\leq& C D_{\Lambda,r}(\sqrt{f})
\end{eqnarray*}
for a constant $C$ depending on $l$. From this it follows that
Propositions \ref{prop:tightc1} and \ref{tightc2} together imply
Proposition \ref{tightc}. We next prove Proposition \ref{tightc2}.
We split it up into several cases, depending on the size of
$\rho$. Up to now our estimates have relied largely on either
one-site bounds or bounds using the local central limit theorem.
Because of this the proofs have been similar to the
non-homogeneous case. However, because of the two-blocks
estimates, the proofs now rely on the joint behaviour over the
boxes.  The methods developed in \cite{logsob2} still apply,
however, with slight modifications.  We begin with some initial
estimates.

\begin{lem}\label{expmom1&2}
\begin{enumerate}
\item For every $\varphi>0$ and $t\in\mathbb{R}$
$$\mu_{\Lambda,\varphi}[e^{t(c_x(\eta_x)-\varphi)}]\leq e^{\varphi a_1 t^2 e^{a_1 |t|}}$$
\item There exists a $C>0$ so that
$$\mu_{\Lambda,\varphi}[e^{t\eta_x}]\leq e^{Ct\rho e^{Ct}}.$$
\end{enumerate}
\end{lem}

\begin{proof}
In the first inequality we repeat the argument of
\reff{roughexpbounds}.  Let
$h(t)=\mu_{\Lambda,\varphi}[e^{tc_x(\eta_x)}]$. Also we remind the
reader of the inequality due to assumption (LG) \linebreak
$|c_x(k+1)-c_x(k)|\leq a_1$.  By a simple change of measure
$$\mu_{\Lambda,\varphi}[c_x(\eta_x)f(\eta_x)]=\varphi
\mu_{\Lambda,\varphi}[f(\eta_x+1)]$$ and Jensen's inequality we
obtain
\begin{eqnarray*}
t h'(t)-h(t)\log h(t)&=&t\mu_{\Lambda,\varphi}[c_x e^{t
c_x}]-\mu_{\Lambda,\varphi}[e^{t c_x}]\log \mu_{\Lambda,\varphi}[e^{t c_x}]\\
&\leq&\varphi t \mu_{\varphi}[e^{tc_x(\eta_x+1)}-e^{tc_x(\eta_x)}]\\
&\leq& a_1\varphi t^2 e^{a_1 t}\mu_{\varphi}[e^{tc_x}].
\end{eqnarray*}
We used the inequality $|e^{x}-e^{y}|\leq |x-y|e^{|x-y|}e^{|y|}$.
Because $$th'(t)-h(t)\log h(t)= t^2 h(t) \partial_t \frac{\log
h(t)}{t}$$ this translates to
$$\partial_t \frac{\log h(t)}{t} \leq \varphi a_1 e^{a_1 t},$$
where $\lim_{t\rightarrow 0} \frac{\log h(t)}{t} =\varphi$.
Integrating we thus have that
$$h(t)=\mu_{\varphi}[e^{tc_x}]\leq e^{\varphi t e^{a_1 t}},$$
which implies
\begin{eqnarray*}
\mu_{\varphi}[e^{t(c_x-\varphi)}]\leq e^{\varphi t ( e^{a_1
t}-1)}\leq e^{\varphi a_1 t^2 e^{a_1 t}}.
\end{eqnarray*}
The bounds on $h(t)$ along with the fact that $c_2\geq
\frac{c_x(k)}{k}\geq c_1$ imply the second inequality.
\end{proof}

We also need a result similar to (i) above for $\eta_x$.
\begin{lem}\label{expmom3}
There exists a $C>0$ so that
$$\mu_{\Lambda,\varphi}\left[e^{t|\eta_x-\rho_x|}\rule[-0.1cm]{0cm}{0.6cm}\right]\leq C e^{C (t\sqrt{\rho_x}+\rho_x t^2 e^{C |t|})}$$
for all $\varphi>0$.
\end{lem}
\begin{proof}
We first use conditions (LG) and (M) to obtain
\begin{eqnarray*}
|\eta_x-\rho_x|&\leq& C \{ |c_x(\eta_x)-c(\rho_x)|+1\}\\
&\leq&
C\{|c_x(\eta_x)-\varphi(\rho_x)|+|c_x(\rho_x)-\varphi(\rho_x)|+1\}\\
&\leq&C\{|c_x(\eta_x)-\varphi(\rho_x)|+\sqrt{\rho_x}+1\}\\
\end{eqnarray*}
where the last inequality
$$|c_x(\rho_x)-\varphi(\rho_x)|\leq \sqrt{\rho_x}$$
is proved as in \cite{lsv}:
\begin{eqnarray*}
|\varphi(\rho_x)-c_x(\rho_x)|\leq
\mu_{\varphi_x}[|c_x(\eta_x)-c_x(\rho_x)|]\leq a_1
\mu_{\varphi_x}[|\eta_x-\rho_x|] \leq a_1 \sigma_x(\rho_x).
\end{eqnarray*}
The remainder now follows from Proposition \ref{sigmabounds} and
Lemma \ref{expmom1&2}.
\end{proof}

We continue with the proof of Proposition \ref{tightc2}.  As
mentioned previously, we split this into two cases:  large density
$\rho>\rho_0$ and small density $\rho\leq\rho_0$.

\subsubsection{Case 1. large density: $\rho>\rho_0$}\label{case1}

For ease of calculation, and without loss of generality, we may
assume that $\nu_{\Lambda,r}[f]=1$. We begin with the entropy
inequality:
\begin{align}
&\nu_{\Lambda, r}[f;\sum_{k=1}^m\nu_{\Lambda_k,
R_k}[\sum_{x\in\Lambda_k}c_x(\eta_x)]]\notag\\
&\leq \frac{1}{t}\log
\nu_{\Lambda,r}[e^{t\sum_{k=1}^m\left\{\nu_{R_k,\Lambda_k}\left[\sum_{x\in\Lambda_k}c_x(\eta_x)\right]-\nu_{r,\Lambda}\left[\sum_{x\in\Lambda_k}c_x(\eta_x)\right]\right\}}]
-\frac{1}{t}H(f|\nu_{\Lambda,r})\label{line:ent}
\end{align}

The next steps will focus on bounding the expectation inside the
logarithm.  Here is where the first difference from the proof of
the homogeneous case appears. The tighter bounds are achieved by
applying a Taylor series type argument to the function $\varphi$
on the boxes $\Lambda_k$.

In what follows, unless otherwise specified, let
$\varphi=\varphi_\Lambda(\frac{r}{|\Lambda|})$ and
$\rho=\frac{r}{|\Lambda|}$.  We let $r_k=\nu_{\Lambda,r}[R_k]$,
where  $R_k$ is the number of particles in $\Lambda_k$ and we also
define \linebreak $\rho_k=\mu_{\Lambda,\varphi}[AV_{x \in
\Lambda_k}\eta_x]$.  We will also denote $\varphi_{\Lambda_k}$ as
$\varphi_k$.  Notice that $\varphi_k(\rho_k)=\varphi(\rho)$.

We define the function
\begin{eqnarray}
\tilde{c}_x(m)=c_x(m)-\varphi'_k(\rho_k)\cdot m,\label{def:newc}
\end{eqnarray}
for all $ x \in \Lambda_k $. For the time being it is enough to
know that $\varphi'_k(x) $ is a strictly positive quantity
uniformly bounded in $x$ for all $k$.  We next bound
\begin{align}
&\nu_{\Lambda,r}[e^{t\sum_{k=1}^m\left\{\nu_{R_k,\Lambda_k}\left[\sum_{x\in\Lambda_k}c_x(\eta_x)\right]-\nu_{r,\Lambda}\left[\sum_{x\in\Lambda_k}c_x(\eta_x)\right]\right\}}]\leq
H_1 \times H_2\notag
\end{align}
where
\begin{align}
H_1 &= e^{-t \sum_{k=1}^m
\left\{\left[\nu_{r,\Lambda}\left[\sum_{x\in\Lambda_k}\tilde{c}_x(\eta_x)\right]-\mu_{\varphi}\left[\sum_{x\in\Lambda_k}\tilde{c}_x(\eta_x)\right]\right]+\varphi'_k(\rho_k)\left[r_k-l
\rho_k \right]\right\}}\notag\\
H_2&=
\nu_{\Lambda,r}[e^{t\sum_{k=1}^m\left\{\left(\nu_{R_k,\Lambda_k}\left[\sum_{x\in\Lambda_k}\tilde{c}_x(\eta_x)\right]-\mu_{\varphi}\left[\sum_{x\in\Lambda_k}\tilde{c}_x(\eta_x)\right]\right)+\varphi'_k(\rho_k)(R_k-\rho_k)\right\}}].\label{line:AB1}
\end{align}
From Propositions \ref{sigmabounds} and \ref{submom}, and from the
Cauchy-Schwarz inequality we have that
\begin{eqnarray*}
|\nu_{\Lambda,r}[c_x(\eta_x)]-\varphi|&\leq&
C\frac{1}{|\Lambda|}\sqrt{\rho_x}, \ \ \mbox{and} \\
|r_k-l\rho_k|&\leq& C\frac{l}{|\Lambda|}\sqrt{l\rho_k},
\end{eqnarray*}
for $|\Lambda|$ sufficiently large, from which it follows that
\begin{align}\label{line:AB2}
H_1 \leq e^{Ct\sqrt{r}}.
\end{align}
Next, by the Cauchy-Schwarz inequality
\begin{eqnarray}
H_2
&\leq&
\nu_{\Lambda,r}[e^{2t\sum_{k=1}^m\left\{\nu_{R_k,\Lambda_k}\left[\sum_{x\in\Lambda_k}\tilde{c}_x(\eta_x)\right]-\mu_{\varphi}\left[\sum_{x\in\Lambda_k}\tilde{c}_x(\eta_x)\right]\right\}}]^{1/2}\notag\\
&&\hskip4cm\times\ \
\nu_{\Lambda,r}[e^{2t\sum_{k=1}^m\left\{\varphi'_k(\rho_k)(R_k-\rho_k)\right\}}]^{1/2}\label{line:AB3}.
\end{eqnarray}  The
second term of these satisfies the following inequality by
applying Cauchy-Schwarz  again
\begin{eqnarray*}
\nu_{\Lambda,r}[e^{t\sum_{k=1}^m\left\{\varphi'_k(\rho_k)(R_k-l\rho_k)\right\}}]
&\leq& \nu_{\Lambda,r}[e^{2t\sum_{k\leq m/2}\left\{\varphi'_k(\rho_k)(R_k-l\rho_k)\right\}}]^{1/2}\\
&&\hskip1cm\times\nu_{\Lambda,r}[e^{2t\sum_{k>m/2}\left\{\varphi'_k(\rho_k)(R_k-l\rho_k)\right\}}]^{1/2}.
\end{eqnarray*}
We next apply Proposition \ref{subprob} to obtain that this is
bounded above by
\begin{eqnarray*}
 C
\left\{\Pi_{k=1}^m\Pi_{x\in\Lambda_k}\mu_{\varphi}[e^{2t\varphi'_k(\rho_k)(\eta_x-\rho_x)}]\right\}^{\frac{1}{2}}.
\end{eqnarray*}
We bound this last quantity using the estimates of Lemma
\ref{expmom3} by
\begin{align}\label{line:tight1}
\nu_{\Lambda,r}[e^{t\sum_{k=1}^m\left\{\varphi'_k(\rho_k)(R_k-l\rho_k)\right\}}]
\leq e^{Ct\sqrt{r}+c\rho t^2e^{Ct}}.
\end{align}
Therefore it remains to bound the first part of line
\reff{line:AB3}. This is where the Taylor argument becomes
important.  Notice that because
$$\rho_k=\mu_{\varphi(\rho)}[AV_{x\in
\Lambda_k}\eta_x],$$ we have that
$\varphi_k(\rho_k)=\varphi(\rho)$.  We have set up a two-block
argument, and we would like to work in measures $\mu$ on
$\Lambda_k$ where the underlying density is $r_k=\frac{R_k}{l}$.
With this in mind, and in a slight abuse of notation, let
$\varphi_k$ denote $\varphi_k\left(\frac{R_k}{l}\right)=
\varphi_{\Lambda_k}\left(\frac{R_k}{l}\right) $.  By the same
argument as above, we may use the Cauchy-Schwarz inequality
together with Proposition \ref{subprob} to obtain
\begin{eqnarray*}
&&\hspace{-1.5cm}\nu_{\Lambda,r}\left[e^{t\sum_{k=1}^m\left\{\left(\nu_{R_k,\Lambda_k}\left[\sum_{x\in\Lambda_k}\tilde{c}_x(\eta_x)\right]-\mu_{\varphi}\left[\sum_{x\in\Lambda_k}\tilde{c}_x(\eta_x)\right]\right)\right\}}\right]\\
&\leq&
C\mu_{\varphi}\left[e^{2t\sum_{k=1}^m\left\{\nu_{R_k,\Lambda_k}\left[\sum_{x\in\Lambda_k}\tilde{c}_x(\eta_x)\right]-\mu_{\varphi}\left[\sum_{x\in\Lambda_k}\tilde{c}_x(\eta_x)\right]\right\}}\right]^{1/2}\\
&\leq& C \mu_{\varphi}\left[e^{4t\sum_{k=1}^m\bar
Y_k}\right]^{1/4}\mu_{\varphi}\left[e^{4t\sum_{k=1}^m\bar
W_k}\right]^{1/4},
\end{eqnarray*}
where we set
\begin{eqnarray*}
Y_k &=&
\nu_{R_k,\Lambda_k}\left[\sum_{x\in\Lambda_k}\tilde{c}_x(\eta_x)\right]-\mu_{\varphi_k}\left[\sum_{x\in\Lambda_k}\tilde{c}_x(\eta_x)\right]\hspace{1cm}\mbox{and}\\
W_k &=&
\mu_{\varphi_k}\left[\sum_{x\in\Lambda_k}\tilde{c}_x(\eta_x)\right]-\mu_{\varphi}\left[\sum_{x\in\Lambda_k}\tilde{c}_x(\eta_x)\right].
\end{eqnarray*}
We let $\bar{Y}_k$ and $\bar W_k$ denote the centered versions of
$Y_k$ and $W_k$ under the measure $\mu_{\varphi}$.  We first bound
$\mu_{\varphi}[e^{t\sum_{k=1}^m\bar Y_k}]$. Notice the following
inequality
\begin{eqnarray*}
e^{x}&=&\sum_{k\geq 0} \frac{x^k}{k!}\,\,\, \leq \,\,\,1+ x +
\frac{x^2}{2} \sum_{k\geq 0} \frac{|x|^k}{k!} \,\,\,\leq \,\,\,1+
x + \frac{x^2}{2}e^{|x|}.
\end{eqnarray*}
This implies that $\exp\{x\} \leq \exp\{x+\frac{x^2}{x}e^{|x|}\}$.
 Applying this inequality to $E[e^X]\leq e^{\{E[X]+E[X^2e^{|X|}]\}}$ we have
\begin{eqnarray*}
\mu_{\varphi}\left[e^{4t\bar Y_k}\right]&\leq& \exp\left\{Ct^2
\mu_{\varphi}\left[\bar Y_k^2e^{t|\bar Y_k|}\right]
\right\}\notag.
\end{eqnarray*}
We next bound $Y_k$ by $C \sqrt{1+\frac{R_k}{l}}$, for some
positive constant $C$ depending on $l$. This is a consequence of
Proposition \ref{submom2}, and holds for $l$ sufficiently large.
The presence of the extra $1$ comes from not being able to bound
$R_k$ from below.  We plug this into the above to obtain
\begin{eqnarray}
&& \hspace{-1cm}\exp\left\{Ct^2
\mu_{\varphi}\left[\left(1+\frac{R_k}{l}\right)e^{t\sqrt{1+\frac{R_k}{l}}}\right] \right\}\notag\\
&\leq&\exp\left\{Ct^2
\mu_{\varphi}\left[\left(1+\frac{R_k}{l})^2\right]^{1/2}\mu_{\varphi}[e^{2t\sqrt{1+\frac{R_k}{l}}}\right]^{1/2} \right\}\notag\\
&\leq&\exp\left\{C(\rho_0)t^2
\rho_k\mu_{\varphi}\left[e^{2t\sqrt{1+\frac{R_k}{l}}}\right]^{1/2}
\right\}.\label{Yline}
\end{eqnarray}
We next bound
$\mu_{\varphi}\left[e^{2t\sqrt{1+\frac{R_k}{l}}}\right]^{1/2}$. We
make use of the following lemma proved \linebreak in
\cite{logsob2}.

\begin{lem}\label{glemma}
Suppose that for a random variable $X\geq 0$ and a function $g:\bb
R_+ \mapsto \bb R_+$ we have that for all $t\geq 0$ we
\begin{eqnarray}\label{glemmaassume}
E[e^{tX}]\leq e^{tg(t)E[X]}.
\end{eqnarray}
Then, for all $t\geq 0$
$$E[e^{t\sqrt{X}}]\leq \exp\{t\sqrt{g(2t)+g(t)}\sqrt{E[X]}\}+e^t.$$
\end{lem}

\begin{proof}
Using the inequality $\sqrt{x}\leq x+1$ as well as Cauchy-Schwarz
and Chebychev's inequalities we have
\begin{eqnarray*}
E[e^{t\sqrt{X}}]&\leq&E[e^{t\sqrt{X}} \bb I [X< k E[X]
]]+E[e^{t\sqrt{X}} \bb I [X\geq k E[X] ]]\\
&\leq&e^{t\sqrt{kE[X]}}+e^t E[e^{tX} \bb I [X\geq k E[X] ]]\\
&\leq&e^{t\sqrt{kE[X]}}+e^t E[e^{2 tX}]^{1/2} P [X\geq k E[X] ]^{1/2}\\
&\leq&e^{t\sqrt{kE[X]}}+e^t E[e^{2 tX}]^{1/2}
\left[\frac{E[e^{tX}]}{e^{tkE[X]}}\right]^{1/2}.
\end{eqnarray*}
Applying twice the assumption \reff{glemmaassume} we have that
\begin{eqnarray*}
e^t E[e^{2 tX}]^{1/2}
[\frac{E[e^{tX}]}{e^{tkE[X]}}]^{1/2}&\leq&e^{\{t
+\frac{1}{2}g(2t)E[X]\}}\left [\frac{e^{t g(t)E[X]
}}{e^{tkE[X]}}\right]^{1/2} \\
&=& e^{t}e^{\frac{t}{2}E[X]\{g(2t)+g(t)-k\}}.
\end{eqnarray*}
The desired result follows if we choose $k=g(2t)+g(t)$.
\end{proof}

\bigskip

\noindent The second part of Lemma \ref{expmom1&2} gives that
\begin{eqnarray*}
\mu_{\varphi}[e^{2t\frac{R_k}{l}}]&\leq&e^{C\rho_k t e^{Ct}}.
\end{eqnarray*}
We may now use Lemma \ref{glemma} with $g(t)=e^{Ct}$.
\begin{eqnarray*}
\mu_{\varphi}[e^{2t\sqrt{1+\frac{R_k}{l}}}]&\leq&e^{2Ct}\mu_{\varphi}[e^{2Ct
\sqrt{\frac{R_k}{l}} }]\\
&\leq&e^{2Ct}(e^{Ce^{t}\sqrt{\rho_k}t}+e^{Ct}).
\end{eqnarray*}
We insert this result into \reff{Yline} and sum over $k$ to get
the following bound
\begin{eqnarray}\label{line:AB4}
\mu_{\varphi}[e^{t\sum_{k=1}^m\bar Y_k}]\leq e^{C\rho t^2
e^{Ct\sqrt{\rho}e^{Ct}}}.
\end{eqnarray}

We now consider the remaining term
$\mu_{\varphi}[e^{t\sum_{k=1}^m\bar W_k}].$ We first write $W_k=l
F_k$ where
\begin{eqnarray}
F_k &=&
\varphi_k\left(\frac{R_k}{l}\right)-\varphi(\rho)-\varphi'_k(\rho_k)
\left(\frac{R_k}{l}-\rho_k\right)\notag\\
&=&
\varphi_k\left(\frac{R_k}{l}\right)-\varphi_k(\rho_k)-\varphi_k'(\rho_k)
\left(\frac{R_k}{l}-\rho_k\right)\label{def:F}
\end{eqnarray}
and define
\begin{eqnarray}\label{def:Z}
Z_k=\frac{1}{l}\sum_{x\in\Lambda_k}\frac{\eta_x-\rho_x}{\sigma_k(\varphi)}.
\end{eqnarray}
Notice that $\sigma_k Z_k = \frac{R_k}{l}-\rho_k$.  Using the
inequalities $e^x \leq 1+x +\frac{1}{2} x^2 e^{|x|}$ and $1+x\leq
e^x$ as before, as well as Cauchy-Schwarz, we get
\begin{eqnarray}\label{line:Fbound}
\log \mu_{\varphi}[e^{tlF_k}]\leq \frac{1}{2}t^2l^2
(\mu_{\Lambda,\varphi}[F_k^4])^{1/2}(\mu_{\Lambda,\varphi}[e^{2tl|F_k|}])^{1/2}.
\end{eqnarray}
Fix a positive constant $B$.  If $|Z_k|> B$ we have that
$|F_k|\leq C \sigma_k |Z_k|$, for some $C$, in which case
\begin{eqnarray}
\mu_{\varphi}[F_k^4 \bb I\{|Z_k| > B\}]&\leq&C\sigma^4
\mu_{\Lambda,\varphi}[Z_k^4 \bb I\{|Z_k| > B\}]\notag\\
&\leq&C\frac{\sigma^4}{B^4}\mu_{\Lambda,\varphi}[Z^8_k]\label{useagain1}
\end{eqnarray}
where we write $\sigma^2$ in lieu of $\sigma^2_\Lambda$.  If
$|Z_k|\leq B$ we may bound $|F_k|\leq C(B)\sigma_k Z_k^2$ by using
a Taylor argument, which also gives us
\begin{eqnarray}
\mu_{\varphi}[F_k^4 \bb I\{|Z_k| \leq  B\}] &\leq&C(B)\sigma^4
\mu_{\Lambda,\varphi}[Z_k^8] .\label{useagain2}
\end{eqnarray}
From Proposition \ref{hatbounds} it follows that
$\mu_{\Lambda,\varphi}[Z_k^8]\leq C l^4$ for some constant $C$.
Using these last results now along with Proposition
\ref{sigmabounds} we obtain that the first part of
\reff{line:Fbound} is bounded by
\begin{eqnarray*}
\mu_{\Lambda,\varphi}[F^4]\leq  C(B)\frac{\rho^2}{l^4}.
\end{eqnarray*}
Using the estimate $|F_k|\leq C\sigma_k |Z_k|$ we get
\begin{eqnarray*}
\mu_{\Lambda,\varphi}[e^{2tl|F_k|}]&\leq&
\mu_{\Lambda,\varphi}[e^{2t\sum_{x\in
C_k}|\eta_x-\rho_x|}]\\
&\leq& e^{Ctl\sqrt{\rho}}e^{Ct^2 l \rho e^{Ct}}
\end{eqnarray*}
using \reff{bd:cond3}.  These last two statements imply that
\reff{line:Fbound} is bounded above by
\begin{eqnarray}\label{line:AB5}
\log \mu_{\varphi}[e^{tlF_k}]&\leq &C \rho
t^2e^{Ctl\sqrt{\rho}}e^{Ct^2 l \rho e^{Ct}}.
\end{eqnarray}

We now combine the bound on $H_1$ from \reff{line:AB2} and the
bounds on $H_2$ from \reff{line:AB3}, line \reff{line:AB4}, as
well as \reff{line:AB5}
 to obtain
\begin{align}
\nu_{\Lambda,r}[e^{t\sum_{k=1}^n\left\{\left(\nu_{R_k,\Lambda_k}\left[\sum_{x\in\Lambda_k}\tilde{c}_x(\eta_x)\right]-\mu_{\varphi}\left[\sum_{x\in\Lambda_k}\tilde{c}_x(\eta_x)\right]\right)\right\}}]
\leq C e^ {Ct\sqrt{r}} e^{Ct^2\frac{r}{l}
e^{Ctl\sqrt{\rho}}e^{Ct^2l^2\rho e^{Ctl}}}\label{line:tight2}
\end{align}
where $\rho=\frac{r}{N}$.  Combining \reff{line:tight2} together
with the entropy inequality \reff{line:ent} we obtain that for any
$t>0$
\begin{align}
\nu_{\Lambda,r}[f;\sum_{k=1}^m\nu_{\Lambda_k,
R_k}[\sum_{x\in\Lambda_k}c_x(\eta_x)]]^2\leq\frac{C}{t^2}+Cr+ct^2\frac{r^2}{l^2}e^{Ctl\sqrt{\rho}}e^{Ct^2l^2\rho
e^{Ctl}} + \frac{1}{t^2}H^2(f|\nu_{\Lambda,r})
\end{align}
We would now like to optimize the above inequality in $t$, as we
have done before.  However, due to the presence of the additional
exponential terms (in comparison with the initial bounds on the
covariances), this is considerably more difficult.  We hence use a
different approach.  We choose $t$ such that
$$t^2=\frac{1 \vee [M H (f| \nu_{\Lambda,r} )]}{r},$$ for any fixed
$M$, and obtain the necessary bounds for three regimes on $t$.

\bigskip
\pagebreak

 \noindent\textbf{Case 1(a):
$t\leq\frac{M}{l\sqrt{\rho}}\wedge M$.}

\medskip

\noindent In this case we obtain the bound
\begin{align}
\nu_{\Lambda,r}[f;\sum_{k=1}^m\nu_{\Lambda_k,
R_k}[\sum_{x\in\Lambda_k}c_x(\eta_x)]]^2\leq Cr +
C\left\{\frac{M}{l}e^{CM+CM^2e^{CM}}+\frac{1}{M}\right\}r
H(f|\nu_{\Lambda,r})\label{line:opt1}
\end{align}

\bigskip

\noindent\textbf{Case 1(b): $t> M$. }

\medskip

\noindent In this setting very rough approximations give
 \begin{align}
\nu_{\Lambda,r}[f;\sum_{x\in\Lambda}c_x(\eta_x)]^2\leq
Cr^2\leq\frac{Cr}{M}H(f|\nu_{\Lambda,r}),
\end{align}
since $ \sum_{x\in\Lambda}c_x(\eta_x) \leq C r$. We used that
$t>M$ implies that $rM \leq H$.  This is clearly a tighter bound
than the one in case 1 (a).

\bigskip

\noindent \textbf{Case 1(c): $\frac{M}{l\sqrt{\rho}}\wedge M<t<M$.
}

\medskip

\noindent Notice that this setting implies that
\begin{align}
1<l\sqrt{\frac{\rho}{Mr}H(f|\nu_{\Lambda,r})}.\label{caseimplies}
\end{align}
This last case is the most complicated of the three, and requires
its own approach.  As before we have
$\varphi_k=\varphi_k\left(\frac{R_k}{l}\right)$. In what follows
we continue to assume that $\nu_{\Lambda,r}[f]=1$. The function
$\tilde{c}_x(k)= c_x(k)-\varphi'_k(\rho_k) k$ is the same function
as before.

By carefully adding and subtracting terms, and noting that
$\nu_{\Lambda,r}[f]=1$, we obtain the following
\begin{align}
&\left|\nu_{\Lambda,r}\left[f;\sum_{k=1}^m\nu_{\Lambda_k,R_k}[\sum_{x\in\Lambda_k}c_x(\eta_x)]\right]\right|\notag\\
&\leq\lf|\nu_{\Lambda,r}\lf[f \cdot \sum_{k=1}^m\left\{\mu_{\varphi_k}[\sum_{x\in\Lambda_k}\tilde{c}_x(\eta_x)]-\mu_{\Lambda,\varphi}[\sum_{x\in\Lambda_k}\tilde{c}_x(\eta_x)]\right\}\rg]\rg|\label{case:3c1}\\
&\hspace{1cm}+\lf|\nu_{\Lambda,r}\lf[f \cdot \sum_{k=1}^m \varphi'_k(\rho_k)\lf\{R_k-\rho_k \cdot l\rg\}\rg]\rg|\label{case:3c2}\\
&\hspace{1cm}+\lf|\nu_{\Lambda,r}\lf[\sum_{k=1}^m\left\{\mu_{\varphi_k}[\sum_{x\in\Lambda_k}c_x(\eta_x)]-\mu_{\Lambda,\varphi}[\sum_{x\in\Lambda_k}c_x(\eta_x)]\right\}\rg]\rg|\label{case:3c3}\\
&\hspace{1cm}+\lf|\nu_{\Lambda,r}\lf[f\cdot\sum_{k=1}^m\left\{\nu_{\Lambda_k,
R_k
}[\sum_{x\in\Lambda_k}c_x(\eta_x)]-\mu_{\varphi_k}[\sum_{x\in\Lambda_k}c_x(\eta_x)]\right\}\rg]\rg|\label{case:3c4}
\end{align}

The terms \reff{case:3c2} and \reff{case:3c4} may each be bounded
by
$$\frac{C}{\sqrt{M}}\sqrt{r H(f|\nu_{\Lambda,r})}$$ in the
following way. We again use Proposition \ref{submom} to get that
$$\nu_{\Lambda_k, R_k
}[c_x(\eta_x)]-\mu_{\varphi_k}[c_x(\eta_x)]\leq C \sqrt
{1+\frac{R_k}{l}},$$ and hence we require $l=|\Lambda_k|$ to be
sufficiently large.  We begin with \reff{case:3c4}. Here we have
that
\begin{eqnarray*}
\nu_{\Lambda,r}\left[f\cdot\sum_{k=1}^m\left\{\nu_{\Lambda_k, R_k
}[\sum_{x\in\Lambda_k}c_x(\eta_x)]-\mu_{\varphi_k}[\sum_{x\in\Lambda_k}c_x(\eta_x)]\right\}\right]
\leq C \nu_{\Lambda,r}\left[f \cdot \sum_{k=1}^m
\sqrt{1+\frac{R_k}{l}} \right].
\end{eqnarray*} We proceed to bound this as
\begin{eqnarray*}
\nu_{\Lambda,r}\left[f \cdot \sum_{k=1}^m
\sqrt{1+\frac{R_k}{l}} \right]&\leq& C m \sqrt{1+\rho}\\
&\leq& C m l \sqrt{1+\rho}
\sqrt{\frac{\rho}{Mr}}\sqrt{H(f|\mu_{N,r})}\\
&\leq& \frac{C}{\sqrt{M}} \sqrt{r H(f|\mu_{N,r})},
\end{eqnarray*}
where we used the fact that under case 1(c) we have $1\leq
l\sqrt{\frac{\rho}{Mr}}\sqrt{H(f|\mu_{N,r})}.$  We handle
\reff{case:3c2} in exactly the same way.

It remains to place bounds on \reff{case:3c1} and \reff{case:3c3}.
In fact, the bounds obtained for \reff{case:3c1} will imply the
necessary bounds for \reff{case:3c3}.   We re-introduce the
notation $F_k$ defined in \reff{def:F} thus obtaining
\begin{align}
\nu_{\Lambda,r}\left[f\cdot\sum_{k=1}^m\left\{\mu_{\varphi_k}[\sum_{x\in\Lambda_k}\tilde{c}_x(\eta_x)]-\mu_{\Lambda,\varphi}[\sum_{x\in\Lambda_k}\tilde{c}_x(\eta_x)]\right\}\right]
=\nu_{\Lambda,r}[f\cdot\sum_{k=1}^m l F_k]\label{line:extra1}
\end{align}
Recall that
$$Z_k=\frac{1}{l}\sum_{x\in\Lambda_k}\frac{\eta_x-\rho_x}{\sigma_k(\varphi)},$$
and define a function $G$
\begin{eqnarray}
G(z)= \left\{\begin{array}{ll}
                C(B)z^2 & \mbox{for $|z|\leq B$}\\
                C(B)B^2 + 2BC(B)(|z|-B) &\mbox{for $|z|> B$.}
                \end{array}
         \right. \label{def:G}
\end{eqnarray}
In the definition of $G$ we choose the constant $B$ so that
$|F_k|\leq \sigma_k G(Z_k)$.  We use this bound in
\reff{line:extra1} along with the second part of Proposition
\ref{sigmabounds} to obtain
\begin{align}
\nu_{\Lambda,r}[f;\sum_{k=1}^m l F_k] &\leq l \sigma \sum_{k=1}^m
\nu_{\Lambda,r}[f\cdot
G(Z_k)]\notag\\
&\leq l\sigma \sum_{k=1}^m
\left\{\nu_{\Lambda,r}[f;G(Z_k)]+\nu_{\Lambda,r}[G(Z_k)]\right\}\label{gline1}
\end{align}

We next introduce the fields $\mathcal{F}_k=\sigma\{\eta_x; x\in
\Lambda_k\cup\Lambda_{k+1}\}$.  For $l$ sufficiently large we may
use Propositions \ref{submom} and \ref{sigmabounds} to get the
following:
\begin{eqnarray}
G_k&\equiv&\nu_{\Lambda,r}[G(Z_k)|\mathcal{F}_k]\leq C(B)
\frac{R_k+R_{k+1}}{l^2 \sigma^2(\varphi)}.\label{line:gbound1}
\end{eqnarray}
These calculations allow us to bound the right hand side of
\reff{gline1}:
\begin{eqnarray}
\ l \sigma \sum_{k=1}^m\nu_{\Lambda,r}[G(Z_k)]
&\leq& \frac{C r}{l \sigma}\notag\\
&\leq& C \sqrt{\frac{r}{M}H_{N,r}(f)}\label{gline2}
\end{eqnarray}
where we have used \reff{caseimplies} in the last line.  It
remains to estimate the left hand side of \reff{gline1}.
\begin{eqnarray}
l \sigma \sum_{k=1}^m \nu_{\Lambda,r}[f; G(Z_k)]
&=& l \sigma \sum_{k=1}^m \nu_{\Lambda,r}[f;G(Z_k)]\label{gline3}\\
&&+ l \sigma \sum_{k=1}^m
\nu_{\Lambda,r}[\nu_{\Lambda,r}[f;G(Z_k)|\mathcal{F}_k]]\label{gline4}
\end{eqnarray}
By an identical argument to that for \reff{gline2} and the fact
that $\nu_{\Lambda,r}[f]=1$ we bound \reff{gline4} by $C
\sqrt{\frac{r}{M}H_{N,r}(f)}$. It remains to study \reff{gline3}.
We use the notation $\nu_k[\cdot]$ for
$\nu_{\Lambda,r}[\cdot|\mathcal{F}_k]$.
\begin{eqnarray}
\nu_{k}[f;G(Z_k)]
&\leq&\left\{\nu_k\lf[\lf(\sqrt{f}-\nu_k[\sqrt{f}]\rg)^2\left|G(Z_k)-G_k\rule[-0.1cm]{0cm}{0.6cm}\right|\rg]\right\}^{1/2}\notag\\
&&\hspace{1cm}\times\left\{\nu_k\left[\left(\sqrt{f}+\nu_k[\sqrt{f}]\right)^2\left|G(Z_k)-G_k\rule[-0.1cm]{0cm}{0.6cm}\right|\right]\right\}^{1/2}\notag\\
&=&\sqrt{V_1 V_2}.\label{line:gbound2}
\end{eqnarray}
Arguing as in \reff{line:gbound1}, but using this time Lipshitz
bounds on $G$, we have that
$$\nu_k[G_k]\leq C'(B) \frac{\sqrt{R_k+R_{k+1}}}{l\sigma},$$
for some  constant $C'(B)$.  Hence
\begin{eqnarray*}
V_1 &\leq& C'(B)
\nu_k\left[(\sqrt{f}-\nu_k[\sqrt{f}])^2\left|\frac{R_k-l\rho_k}{l\sigma}+\frac{\sqrt{R_k+R_{k+1}}}{l\sigma}\right|\right]
\end{eqnarray*}
(where $C'(B)$ may be a new constant) which by the entropy
inequality may be bounded by (writing $C=C'$)
\begin{eqnarray*}
\frac{C}{l\sigma}\nu_k\left[(\sqrt{f}-\nu_k[\sqrt{f}])^2\right]\left(\frac{1}{s}\log\nu_k[e^{s|R_k-l\rho_k|}]+\sqrt{R_k+R_{k+1}}\right)\\
\hspace{2cm}+\ \ \frac{C}{s l \sigma}
H((\sqrt{f}-\nu_k[\sqrt{f}])^2|\nu_k).
\end{eqnarray*}
Notice that the above inequality holds because
$\nu_k\left[(\sqrt{f}-\nu_k[\sqrt{f}])^2\right]$ is not equal to
one. Using Lemma \ref{expmom3} in the above this is smaller than
\begin{eqnarray*}
V_1\leq\frac{C}{l\sigma}\nu_k[\sqrt{f};\sqrt{f}]\left(\frac{1}{s}+s (R_k+R_{k+1})+\sqrt{R_k+R_{k+1}}\right)\\
\hspace{1cm} +\ \ \ \frac{C}{s l \sigma}
H((\sqrt{f}-\nu_k[\sqrt{f}])^2|\nu_k),
\end{eqnarray*}
where the constant $C$ in front may depend on $M$. Optimizing over
$s$ we next get
\begin{eqnarray*}
V_1\ &\leq&\
\frac{C}{l\sigma}\nu_k[\sqrt{f};\sqrt{f}]\left(\sqrt{R_k+R_{k+1}}(1+\sqrt{\frac{H((\sqrt{f}-\nu_k[\sqrt{f}])^2|\nu_k)}{\nu_k[\sqrt{f};\sqrt{f}]}})\right)
\end{eqnarray*}
Let $D_{k,k+1}(\sqrt{f})$ denote the Dirichlet form of the process
defined over $\Lambda_k\cup\Lambda_{k+1}$.  We next apply the
logarithmic Sobolev inequality in the line above, with constant
$\kappa=\kappa(2l)$
\begin{eqnarray*}
\frac{C\sqrt{R_k+R_{k+1}}}{l\sigma}\nu_k[\sqrt{f};\sqrt{f}]+\frac{C\sqrt{R_k+R_{k+1}}}{l\sigma}\sqrt{\kappa(2l)\nu_k[\sqrt{f};\sqrt{f}]D_{k,k+1}(\sqrt{f})}.
\end{eqnarray*}
 We repeat the argument for $V_2$ to obtain
\begin{eqnarray*}
V_2&\leq&\frac{C\sqrt{R_k+R_{k+1}}}{l\sigma}\nu_k[f]+\frac{C\sqrt{R_k+R_{k+1}}}{l\sigma}\sqrt{\kappa(2l)\nu_k[f]D_{k,k+1}(\sqrt{f})}.
\end{eqnarray*}
We now use the spectral gap result $\nu_k[\sqrt f; \sqrt f]\leq
Cl^2D_k(\sqrt{f})$ and $\nu_k[\sqrt f; \sqrt f]\leq\nu_k[f]$ in
the above bounds to obtain for some constant $C$ depending on $l$
$$\nu_{\Lambda,r}[\nu_{k}[f;G(Z_k)]]\leq \frac{C}{\sigma}\sqrt{\nu_{\Lambda,r}[(R_k+R_{k+1})f]\nu_{\Lambda,r}[D_{k,k+1}({\sqrt{f}})]}. $$
We use this to compute the quantity of interest
\begin{align*}
l \sigma \sum_{k=1}^m \nu_{\Lambda,r}[f;G_k]
&= l \sigma \sum_{k=1}^m \nu_{\Lambda,r}[\nu_{k}[f;G_k]]\\
&\leq l \sigma \sum_{k=1}^m
\frac{C}{\sigma}\sqrt{\nu_{\Lambda,r}[(R_k+R_{k+1})f]\nu_{\Lambda,r}[D_k({\sqrt{f}})]}\\
& \leq C l \sum_{k=1}^m \sqrt{\nu_{\Lambda,r}[R_k
f]\nu_{\Lambda,r}[D_k({\sqrt{f}})]}\\
& \leq C l m \sqrt{r\nu_{\Lambda,r}[D_k({\sqrt{f}})]}\\
& \leq C N \sqrt{r D_{\Lambda,r}({\sqrt{f}})}
\end{align*}

We may now combine the above line together with the bounds
obtained for \reff{case:3c2},\reff{case:3c3} and \reff{case:3c4}
to obtain that
\begin{align*}
\nu_{\Lambda,r}[f;\sum_{k=1}^m\nu_{\Lambda_k,R_k}[\sum_{x\in\Lambda_k}c_x(\eta_x)]]^2
\leq  r \left(\frac{C}{M}H_{N,r}(f)+ C(M)N^2
D_{\Lambda,r}({\sqrt{f}})\right)
\end{align*}

We combine the results of Case 1, (a) through (c), and choose
$\epsilon = \frac{C}{M}$ to obtain Proposition \ref{tightc2} for
any large density. Notice that in the above work although we need
to choose $l$ sufficiently large for certain bounds to hold, once
we do so, it remains fixed.  We show next the necessary bounds for
small density with appropriate choice of cutoff $\rho_0$.

\subsubsection{Case 2. small density:
$\rho\leq\rho_0$}\label{case2}

\begin{lem}\label{tightsmall}
For every $\epsilon>0$ there exists a $\rho_0$ and a constant
$C=C(\epsilon)>0$ and an $N_0$ so that for $\frac{r}{N}\leq
\rho_0$ and $N\geq N_0$
$$\nu_{\Lambda,r}[f;\sum_{x\in\Lambda}c_x(\eta_x)]^2\leq r \nu_{\Lambda,r}[f](C\nu_{\Lambda,r}[f] +\epsilon H(f|\nu_{\Lambda,r})).$$
\end{lem}

\begin{proof}  We assume $\nu_{\Lambda,r}[f]=1$, again without loss of
generality.  By the entropy inequality
\begin{align}\label{yetanent}
\nu_{\Lambda,r}[f;\sum_{x\in\Lambda}c_x(\eta_x)]\leq
\frac{1}{t}\log\nu_{\Lambda,r}[e^{t\sum_{x\in
\Lambda}(c_x(\eta_x)-\nu_{\Lambda,r}[c_x(\eta_x)])}]+\frac{1}{t}H(f|\nu_{\Lambda,r})
\end{align}
Now,
\begin{align*}
\nu_{\Lambda,r}[e^{t\sum_{x\in
\Lambda}(c_x(\eta_x)-\nu_{\Lambda,r}[c_x(\eta_x)])}]\leq C
e^{Ct\sqrt{r}}\Pi_{x\in\Lambda}\mu_{\varphi(\frac{r}{N})}[e^{t(c_x(\eta_x)-\frac{\varphi}{\rho_x}\eta_x)}]
\end{align*}
by the Cauchy-Schwarz inequality, and Propositions \ref{subprob}
and \ref{submom2}. We next handle the term
$\mu_{\varphi(\frac{r}{N})}[e^{t(c_x(\eta_x)-\frac{\varphi}{\rho_x}\eta_x)}]$.
We have previously obtained bounds in Lemma \ref{expmom1&2}, but
these are not sufficient here.  Fix a constant $M$ and for $t$ in
$[0,M]$ and $\varphi$ in $[0,1]$ define
\begin{eqnarray}\label{newF}
F(t,\varphi)=\mu_{\varphi}[e^{t(c_x(\eta_x)-\frac{\varphi}{\rho_x}\eta_x)}].
\end{eqnarray}
We notice a few things about the function $F$:
\begin{enumerate}
\item[$\bullet$] $F(0,\varphi)=1$,
\item[$\bullet$] $\partial_tF(0,\varphi)=0$,
\item[$\bullet$] $\partial_t^2
F(0,\varphi)=\mu_{\varphi}[(c_x(\eta_x)-\frac{\varphi}{\rho_x}\eta_x)^2]$.
\end{enumerate}
We wish to bound these derivatives and integrate to obtain an
appropriate bound on the function $F$.  Because we are in the
setting of bounded densities and bounded $t$ the function $F$ as
well as its derivatives are well behaved.  We also have the
following:
\begin{enumerate}
\item[$\bullet$]   $\partial_t \partial_\varphi^k
F(0,0)=0$ for all $k$,
\item[$\bullet$] $\partial_\varphi \partial_t^2 F(0,0)=0$.
\end{enumerate}
By the above, there exists a constant $C(M)$, possibly depending
on $M$,  such that
\begin{eqnarray}\label{newF2}
F(t,\varphi)\leq 1+C(M)\varphi^2 t^2 \leq e^{C(M)\varphi^2t^2}
\end{eqnarray}
for $t$ in $[0,M]$ and $\varphi \leq 1$. We next replace $\varphi$
in the above by $C \rho$. Inserting these bounds into the entropy
inequality \reff{yetanent} we get
\begin{eqnarray*}
\nu_{\Lambda,r}[f;\sum_{x\in\Lambda}c_x(\eta_x)]^2\leq
\frac{C}{t^2}+Cr+C(M)r^2\rho^2
t^2+\frac{1}{t^2}H^2(f|\nu_{\Lambda,r}),
\end{eqnarray*}
for all $\rho$ such that $\varphi(\rho)\leq 1$.  We again wish to
optimize over $t$. As before, choose
$$t^2=\frac{1\vee[M H(f|\nu_{\Lambda,r})]}{r}.$$  As long as this $t \leq M$
we may plug it into the above bound to get
\begin{eqnarray*}
\nu_{\Lambda,r}[f;\sum_{x\in\Lambda}c_x(\eta_x)]^2\leq
2Cr+C(M)r\rho^2 (1\vee [M
H(f|\nu_{\Lambda,r})])+\frac{r}{M}H(f|\nu_{\Lambda,r}).
\end{eqnarray*}
Otherwise, because $\sum_{x\in\Lambda}c_x(\eta_x)\leq r$,  and
$t>M$ implies that $rM<H$, we have the easier bound
\begin{eqnarray*}
\nu_{\Lambda,r}[f;\sum_{x\in\Lambda}c_x(\eta_x)]^2\leq Cr^2 \leq
\frac{Cr}{M}H(f|\nu_{\Lambda,r}).
\end{eqnarray*}
We now choose $\epsilon=\frac{C}{M}$ and then $\rho_0$ small
enough so that the result follows.
\end{proof}

\subsection{Proof of Proposition \ref{tighth}}

\begin{proof}
In what follows, assume $\varphi=\varphi(\rho)$ with
$\rho=\frac{r}{N}$.   We will prove instead the inequality
$$\nu_{\Lambda,r}[f;\sum_{x\in\Lambda}\varphi h_x(\eta_x)]^2\leq r
\nu_{\Lambda,r}[f]\left[C \nu_{\Lambda,r}[f]+ C N^2
D_{\Lambda,r}(\sqrt{f})+\epsilon H(f|\nu_{\Lambda,r})\right],$$
from which the desired result follows.  Because $c_x$ and $\varphi
h_x$ are of the same order, we may follow the proof of Proposition
\ref{tight1} with only mild modifications.   We first consider the
case of $\rho$ bounded below.  All of the arguments go through as
before, until \reff{Yline}, where we need to know that
$$|\nu_{\Lambda,r}[\varphi h_x(\eta_x)]-\mu_{\varphi}[\varphi h_x(\eta_x)]|\leq \frac{C}{N}\sqrt{1+\rho}.$$
This is a consequence of Proposition \ref{submom2} as in the case
of $c_x$ as long as we can show that
$$\mu_{\Lambda,\varphi}[h_x;h_x]\leq \frac{C}{\rho}.$$
Because of Lemma \ref{Kawasakispecgap} we have
\begin{eqnarray*}
\mu_{\Lambda,\varphi}[h_x;h_x]&\leq&
B\mu_{\Lambda,\varphi}[c_x(\eta_x)(h_x(\eta_x-1)-h_x(\eta_x))^2],
\end{eqnarray*}
where $B$ holds for all $x$. This last quantity is less than
\begin{eqnarray*}
&&\hspace{-1cm}B
\mu_{\Lambda,\varphi}\left[\frac{1}{c_x(\eta_x)c^2_x(\eta_x+1)}(\eta_x(c_x(\eta_x+1)-c_x(
\eta_x)))^2\right]\\
&\leq & C B\mu_{\Lambda,\varphi}\left[\frac{1}{c_x(\eta_x+1)}\right]\\
&\leq& \frac{C B}{\rho},
\end{eqnarray*}
for some positive constant $C$. To continue with the argument we
need to also specify how to handle the Taylor approximation
arguments involving the functions $\tilde{c}_x$ and $F_k$ (cf.
\reff{def:newc} and \reff{def:F}).  Here we use
$$\varphi \tilde{h}_x = \varphi h_x - \varphi \gamma'_k(\rho_k) \eta_x,$$
for all $x$ in $\Lambda_k$ where
$\gamma_k(x)=\frac{x}{\varphi_k(x)}$. This implies that the new
version of $F_k$ becomes
\begin{eqnarray*}
F_k &=&\varphi\left[\frac{\frac{R_k}{l}}{\varphi_k(\frac{R_k}{l})}-\frac{\rho_k}{\varphi(\rho)}\right]-\varphi\gamma'_k(\rho_k)\left[\frac{R_k}{l}-\rho_k\right]\\
&=&\varphi\left[\frac{\frac{R_k}{l}}{\varphi_k(\frac{R_k}{l})}-\frac{\rho_k}{\varphi_k(\rho_k)}\right]-\varphi\gamma'_k(\rho_k)\left[\frac{R_k}{l}-\rho_k\right].
\end{eqnarray*}
So that the argument of the previous section goes through we need
to know two things: we need $\delta_k$ to be uniformly bounded and
we need to be able to choose a constant $B$ so that $|F_k|\leq
\sigma_k G(Z_k)$ with $G$ and $Z$ defined in \reff{def:G} and
\reff{def:Z}.

By direct calculation we have
\begin{eqnarray}\label{gammacalc}
\gamma_k'(x)=\frac{1}{\varphi_k(x)}\lf\{1-\frac{x}{\sigma^2_k(x)}\rg\},
\end{eqnarray}
which implies that $\varphi \gamma_k'(\rho_k)$ is uniformly
bounded by Proposition \ref{sigmabounds}.

We next need to show that $\varphi \gamma_k$ is Lipschitz and that
we can choose a $B$ so that $\gamma''_k(x)$ is bounded when
$|x-\rho_k|\leq B \sigma_k$.  We first show that
$$|\gamma_k(x)-\gamma_k(\rho_k)|\leq \frac{C}{\varphi}\ |x-\rho_k|.$$
Because $\gamma_k$ is uniformly bounded (in $k$) the inequality is
immediate for $|x-\rho_k|>\frac{\rho_k}{2}$. Also, by
\reff{gammacalc}, we know that $\gamma_k'(x)$ is bounded for $x$
away from zero, and this gives us the necessary bound for
$|x-\rho_k|$ small. Lastly, we need to show that the second
derivative of $\gamma_k$ is well behaved when $|x-\rho_k|\leq B
\sigma_k$ for some choice of $B$, so that we may obtain the
tighter bounds for smaller values of $Z_k$.  We calculate directly
$$\gamma''(x)=-2\frac{\varphi'_k(x)}{\varphi_k^2(x)}+2x\frac{\varphi_k'(x)^2}{\varphi_k^3(x)}-x\frac{\varphi_k''(x)}{\varphi_k^2(x)}.$$
Careful examination of the above reveals that this does indeed
remain bounded for $x$ away from zero, and hence we have the
necessary bounds as long as we consider $\rho$ away from zero and
$B$ sufficiently small. The rest of the proof of Proposition
\ref{tighth} follows as before for the case of $\rho$ bounded
below.

For the case of $\rho$ small (bounded above), we need to specify
the version of
$$c_x(\eta_x)-\frac{\varphi}{\rho_x}\eta_x$$ required in the definition of \reff{newF}.  We use in
this case
$$\varphi h_x(\eta_x)-\rho_x.$$
This function defines the new $F(t,\varphi)=\mu_{\varphi}[e^{t
(\varphi h_x(\eta_x)-\rho_x)}]$  in lieu of \reff{newF} from
before.  It is straightforward to check that this new function
satisfies all of the required properties so that we obtain the
appropriate bounds \reff{newF2}. The rest of the proof goes
through without further changes.

\end{proof}

\section{Proof of Spectral Gap for Inhomogeneous Zero Range.\label{sec:specgap}}

This section is dedicated to the proof of the Theorem
\ref{specgap}.  The method is the same as that used in \cite{lsv},
while carefully making sure all necessary bounds hold uniformly in
the sites.  We will frequently make use of the spectral gap
results obtained in section \ref{birthdeath} for certain birth and
death processes.  We present the proof for a general dimension
$d$, in order to highlight the changes needed to extend the proof
of the logarithmic Sobolev inequality to higher dimensions.

\subsection{Outline of Proof.}

Let $\omega(N,r)$ be the smallest constant such that
$$\nu_{\Lambda,r'}[f;f]\leq \omega(N,r) D_{\Lambda,r'}(f)$$
for all $|\Lambda|\leq N^d$ and all $r'\leq r$.  Also, let
$\omega(N)=\sup_k\omega(N,k)$.

The general approach here is quite similar to the one in the
logarithmic Sobolev inequality: using induction we establish two
recursive equations for $\omega$.  The first equation allows us to
establish that the constant is free of the number of particles,
while the second, valid only for sufficiently large $N$, gives the
$N^2$ order.  The difference with the approach of the logarithmic
Sobolev inequality is that the induction increment \emph{adds} one
to the side length of the cube $\Lambda$, and does not
\emph{double} it.

Assume then that $\Lambda$ is a set of size $|\Lambda|=N^d$, and
write $\Lambda = \Lambda_0 \cup \Lambda_1$, where
$|\Lambda_1|=(N-1)^d$ and $\Lambda_1$ is still a cube containing
one of the corner points.  That is, if $d=1$ then $\Lambda_1 =
\{1,\ldots,N-1\}$, if $d=2$ then $\Lambda_1= \{(z_1,z_2);
z_i=1,\ldots, N-1 \}$, and so forth. Denote by $R_0=R_0(\eta)$ the
random variable counting the number of particles in
$\Lambda_0=\Lambda \setminus\Lambda_1$. Also, enumerate the sites
of $\Lambda_0$ so that $\eta_{\Lambda_0} =
\{\eta_{z_k}\}_{k=1}^{N^*},$ where $N^*=N^d-(N-1)^d.$  For
simpilicity, we denote $z_k$ simply as $k$.  Lastly, let $\mc
F_k=\sigma\{\eta_1, \ldots, \eta_k\}$ denote the $\sigma$-algebra
generated by the first $k$ elements of $\Lambda_0$. Thus, $\mc
F_k$ forms an increasing filtration, and denoting
$\nu_{\Lambda,r}[f|\mc F_k]$ by $f_k$, we may write
\begin{eqnarray}\label{sg:line1}
\nu_{\Lambda,r}[f;f]=\nu_{\Lambda,r}[\nu_{\Lambda_1,r-R_0}[f;f]]+\sum_{k=0}^{N^*}\nu_{\Lambda,r}[(f_{k+1}-f_k)^2].
\end{eqnarray}
Here, $\mc F_0$ denotes the trivial $\sigma$-algebra.  By the
induction hypothesis, we may bound the first term above by
\begin{eqnarray}\label{sg:line1.1}
\omega(N-1)\nu_{\Lambda,r}[D_{\Lambda_1,r-R_0}(f)]\leq
\omega(N-1)D_{\Lambda,r}(f).
\end{eqnarray}
To bound the second term we write
\begin{eqnarray*}
\nu_{\Lambda,r}[(f_{k+1}-f_{k})^2]&=&\nu_{\Lambda,r}[\nu_{\Omega_k,R_k}[(f_{k+1}-f_{k})^2]]\notag\\
&=&\nu_{\Lambda,r}[\nu_{\Omega_k,R_k}[f_{k+1};f_{k+1}]]\notag
\end{eqnarray*}
where $\Omega_k=\Lambda\setminus\{z_m\in\Lambda_0\}_{m=1}^k$ and
$R_k$ is the number of particles there.  Restricting consideration
to the measure $\nu_{\Omega_k,R_k}$ we think of $f_{k+1}$, a
function of $\{\eta_m\}_{m=1}^{k+1}$ as a function only of
$\eta_{k+1}$, imagining the remaining sites to be fixed.  We thus
write $f_{k+1}$, with a slight abuse of notation, as
$\phi_k(\eta_{k+1})$. We obtain
\begin{eqnarray}
\nu_{\Omega_k,R_k}[f_{k+1},f_{k+1}]&=&\nu_{\Omega_k,R_k}[\phi_k,\phi_{k}]\notag\\
 &\leq&B_1\
\nu_{\Omega_k,R_k}[c_k(\eta_{k+1})\{\phi_k(\eta_{k+1}-1)-\phi_k(\eta_{k+1})\}^2],\label{sg:line2}
\end{eqnarray}
by Lemma \ref{bdspecgap}, where the constant $B_1$ does not depend
on the location of the site $k+1$.  We next write \reff{sg:line2}
as
\begin{eqnarray}\label{sg:line3}
B_1 \
\sum_{m=0}^{R_k-1}\nu_{\Omega_k,r_k}(\eta_{k+1}=m+1)c_{k+1}(m+1)\{\phi_k(m+1)-\phi_k(m)\}^2.
\end{eqnarray}
Using a calculation similar to the more general one of Proposition
\ref{above1} (indeed, this is just a special case of that result)
we obtain that
\begin{eqnarray*}
\phi_k(m+1)-\phi_k(m)&=&\ \nu_{\Omega_k,R_k}[f|\eta_{k+1}=m+1]-\nu_{\Omega_k,R_k}[f|\eta_{k+1}=m]\\
&=&\
\frac{1}{\nu_{\Omega_k,R_k}(\eta_{k+1}=m+1)c_k(m+1)}\left\{A_1+A_2\right\},
\end{eqnarray*}
where
\begin{eqnarray*}
A_1(k,\Lambda,r,f)&=&AV_{y\in\Omega_{k+1}}
\nu_{\Omega_k,R_k}[c_y(\eta_y)\nabla_{y,k+1}f \ \bb I
(\eta_{k+1}=m)],\\
A_2(k,\Lambda,r,f)&=&
\nu_{\Omega_k,R_k}[\eta_{k+1}=m]\nu_{\Omega_k,R_k}[f;AV_{y\in\Lambda_{k+1}}c_y(\eta_y)|\eta_{k+1}=m].
\end{eqnarray*}
This means that \reff{sg:line2} is bounded above by (a constant,
$2B_1$, times) the sum of $T_1$ and $T_2$, where
\begin{align*}
&T_1 = \sum_{m=0}^{R_k-1}\frac{\left\{AV_{y\in\Omega_{k+1}}
\nu_{\Omega_k,R_k}[c_y(\eta_y)\nabla_{y,k+1}f  \ \bb I
(\eta_{k+1}=m)]\right\}^2}{\nu_{\Omega_k,R_k}(\eta_{k+1}=m+1)c_k(m+1)},\  \mbox{  and}\\
&T_2 =
\sum_{m=0}^{R_k-1}\frac{\nu_{\Omega_k,R_k}(\eta_{k+1}=m)}{\nu_{\Omega_{k+1},R_k-m}[c_{k}(\eta_{k})]}\left\{\nu_{\Omega_k,R_k}[f;AV_{y\in\Lambda_{k+1}}c_y(\eta_y)|\eta_{k+1}=m]\right\}^2,
\end{align*}
and we have used the relation
$$\nu_{\Lambda\setminus\{z\},r-\eta_z}[c_y(\eta_y)]\nu_{\Lambda,r}(\eta_z=m)=\nu_{\Lambda,r}(\eta_z=m+1)c_z(m+1),$$
 for any $z,y \in \Lambda$, in the latter. The next steps establish bounds on these
terms.
\begin{prop}\label{T1} There exists a finite constant $C$ such
that
$$\sum_{k=1}^{N^*}\nu_{\Lambda,r}[T_1(k,\Lambda,r,f)] \leq C N D_{\Lambda,r}(f).$$
\end{prop}
This is a universal bound on $T_1$. We also need to establish both
a weak and a strong version of bounds on $T_2$, to be used in the
recursive equations.
\begin{prop}\label{T2}  There exists a finite constant $C$ such
that
$$T_2(k,\Lambda,r,f) \leq C \omega(N-1) D_{\Omega_k,R_k}(f).$$
\end{prop}
Using local limit theorems, this may be tightened for sufficiently
large values of $N$.
\begin{prop}\label{T2plus}
For all $\epsilon>0$, there exist finite constants $n_0(\epsilon)$
and $C(\epsilon)$ such that
$$T_2(k,\Lambda,r,f) \leq C(\epsilon)N^{-d}D_{\Omega_k,R_k}(f)+\epsilon N^{-d}\nu_{\Lambda_K,R_k}[f;f]$$
for all $n\geq n_0$.
\end{prop}

 We may now
combine these propositions to prove the result. First,
Propositions \ref{T1} and \ref{T2} applied in \reff{sg:line1}
together with \reff{sg:line1.1} give
\begin{eqnarray}
\nu_{\Lambda,r}[f;f]&\leq& \omega(N-1)D_{\Lambda,r}(f)+\sum_{k=0}^{N^*}\nu_{\Lambda,r}[\nu_{\Omega_k,R_k}[f_{k+1};f_{k+1}]]\notag\\
&\leq& \omega(N-1)D_{\Lambda,r}(f)+CND_{\Lambda,r}(f)+CN^{d-1}\sum_{k=0}^{N^*}\omega(N-1)\nu_{\Lambda,r}[D_{\Omega_k,R_k}(f)]\notag\\
&\leq&\left\{\left(1+CN^{d-1}\right)\omega(N-1)+CN\right\}D_{\Lambda,r}(f),\label{loose}
\end{eqnarray}
since $\nu_{\Lambda,r}[D_{\Omega_k,R_k}(f)]\leq D_{\Lambda,r}(f)$
and $N^*\leq CN^{d-1}$ for some constant $C$.  Tightening these
bounds using Proposition \ref{T2plus} we have for any $\epsilon$
\begin{align}\label{tight}
\nu_{\Lambda,r}[f;f]\leq
\left\{CN+C(\epsilon)+\omega(N-1)\right\}D_{\Lambda,r}(f)+B_1\epsilon
N^{-d}\nu_{\Lambda,r}[f;f] \end{align}
 for sufficiently large $N$.

Notice that the initial induction step, $\omega(2)<\infty$, is
established in Lemma \ref{bdspecgap}, because when $|\Lambda|=2$,
then $f(\eta)=f(\eta_1, r-\eta_1 )=\phi(\eta_1)$.  From
\reff{loose} we have that there is some constant, $C(N)$,
independent of the number of particles, such that
$$\omega(N)\leq C (N)\omega(N-1)+\frac{N}{2}.$$
This recursive equation implies that $\omega(N)$ is finite for
every $N$.

Similarly using the tighter bounds from \reff{tight} we obtain
also that for all $\epsilon$ and sufficiently large values of $N$
$$\omega(N)\leq (1-\epsilon/N^d)^{-1}\left\{B_1C(\epsilon)+\omega(N-1)+CN\right\},$$
which implies the required quadratic growth.  Thus the two
recursive formulae above, along with the initial induction step,
establish Theorem \ref{specgap}.

\subsection{Proof of Proposition \ref{T1}.}

From
$$\nu_{\Omega,m}[c_z(\eta_z)\mathbb
I(\eta_y=k)]=\nu_{\Omega,m}(\eta_y=k+1)c_{y}(k+1),$$ it follows by
the Schwarz inequality that
\begin{align}
&\sum_{k=1}^{N^*}\nu_{\Lambda,r}[T_1(k,\Lambda,r,f)]\notag\\
&=\sum_{k=1}^{N^*}\nu_{\Lambda,r}\left[\sum_{m=0}^{R_k-1}\frac{\left\{AV_{y\in\Omega_{k+1}}
\nu_{\Omega_k,R_k}[c_y(\eta_y)\nabla_{y,k+1}f \ \bb I
(\eta_{k+1}=m)]\right\}^2}{\nu_{\Omega_k,R_k}(\eta_{k+1}=m+1)c_k(m+1)}\right]\notag\\
&\ \leq\sum_{k=1}^{N^*}\nu_{\Lambda,r}\left[\nu_{\Omega_k,R_k}[AV_{y\in\Omega_{k+1}}c_y(\eta_y)\{\nabla_{y,k+1}f\}^2]\right]\notag\\
&\
=\frac{C}{N^d}\sum_{k=1}^{N^*}\sum_{y\in\Omega_{k+1}}\nu_{\Lambda,r}[c_y(\eta_y)\{\nabla_{y,k+1}
f \}^2].\label{sg:line55}
\end{align}
We next bound $\{\nabla_{y,k+1}f \}^2$ by
$$CN\sum_{\{e_z\}}\{f(\eta^{x,e_{z+1}})-f(\eta^{x,e_z})\}^2,$$
where the sum is over $\{e_z\}$: sites which form a path from $y$
to $k+1$.  We pick these paths in a particular way.  We number the
directions from $1$ to $d$.  The path from a site $y$ to a site
$k+1$ is a path such that we move maximally in the first
direction, then maximally in the second direction, etc.. For
example, in $d=2$ to join a site $x=\{x_1,x_2\}$ with a site
$y=\{y_1,y_2\}$ such that $x_1>y_1$ and $x_2<y_2$ we choose the
path such that there exists and $m$ so that $e_1 = \{x_1,x_2\},
e_2 = \{x_1-1,x_2\}, \ldots, \ e_m =\{y_1,x_2\}$ and
$e_{m+1}=\{y_1, x_2+1\}, e_{m+2}=\{y_1, x_2+2\}, \ldots, \
e_{m^*}=\{y_1,y_2\}$. Here $m^* = |x_1-y_1|+|x_2-y_2|+1.$ Using
this decomposition we have that
\begin{align*}
&\hspace{-.5cm}\nu_{\Lambda,r}\left[c_y(\eta_y)\{\nabla_{y,k+1}f\}^2\rule[-0.1cm]{0cm}{0.6cm}\right]\\
&\leq
CN\sum_{\{e_z\}}\nu_{\Lambda,r}[c_y(\eta_y)\{f(\eta^{x,e_{z+1}})-f(\eta^{x,e_z})\}^2]\\
& =
CN\sum_{\{e_z\}}\nu_{\Lambda,r}[c_{e_z}(\eta_{e_z})\{\nabla_{e_z,e_{z+1}}f\}^2].
\end{align*}
By changing the order of summation, we conclude that
\reff{sg:line55} is bounded above by
\begin{align*}
\frac{CN}{N^d}\sum_{k=1}^{N^*}\sum_{y\in\Omega_{k+1}}\sum_{\{e_z\}}\nu_{\Lambda,r}[c_{e_z}(\eta_{e_z})\{\nabla_{e_z,e_{z+1}}f\}^2]\\
\leq \frac{CN}{N^d} \sum_{w\sim v \in \Lambda}
\nu_{\Lambda,r}[c_{e_z}(\eta_{e_z})\{\nabla_{e_z,e_{z+1}}f\}^2] \
\sum _{k,x}1 \  ,
\end{align*}
where the last sum is taken over $k,x$ all sites $k$ and $x$ such
that both sites $w$ and $v$ are in the path $\{e_z\}$ from $x$ to
$k$.  Because of our construction, this last quantity is bounded
by a constant times $N^d$.  This concludes the proof.

\subsection{Proof of Proposition \ref{T2}.}

By the Schwarz inequality
\begin{eqnarray}
&&\hspace{-2cm}\left\{\nu_{\Omega_k,R_k}[f;AV_{y\in\Omega_{k+1}}c_y(\eta_y)|\eta_{k+1}=m]\right\}^2\notag\\
&\leq&
\nu_{\Omega_{k+1},R_k-m}[f;f]\nu_{\Omega_{k+1},R_k-m}[AV_{y\in\Omega_{k+1}}c_y(\eta_y);AV_{y\in\Omega_{k+1}}c_y(\eta_y)].\label{sg:line4}
\end{eqnarray}
For any set $\Omega$ and $R$, and fixed $z\in \Omega$, we apply
Lemma \ref{bdspecgap}
\begin{eqnarray*}
\nu_{\Omega,R}[AV_{y\in\Omega}c_y(\eta_y);AV_{y\in\Omega}c_y(\eta_y)]&\leq&\nu_{\Omega,R}[c_z(\eta_z);c_z(\eta_z)]\\
&\leq&B_1
\nu_{\Omega,R}[c_z(\eta_z)\{c_z(\eta_z-1)-c_z(\eta_z)\}^2].
\end{eqnarray*}
By assumption (LG) this last expression is smaller than $a_1^2
\nu_{\Omega,R}[c_z(\eta_z)]$.  Plugging this back into
\reff{sg:line4}, we obtain
$$T_2\leq \nu_{\Omega_k, R_k}[\nu_{\Omega_{k+1},R_k-\eta_{k+1}}[f;f]].$$
Applying the induction hypothesis to this last statement, we show
\begin{eqnarray*}
T_2&\leq& \nu_{\Omega_k, R_k}[\omega(N-1)D_{\Omega_{k+1}}(f)]\\
&\leq& \omega(N-1)D_{\Omega_k,R_k}(f),
\end{eqnarray*}
as desired.

\subsection{Proof of Proposition \ref{T2plus}.}

Set $\Omega^+=\Omega\cup\{z\}$ for a fixed site $z$, where
$\Omega$ is such that $|\Omega|\geq CN^d$.  Fix another site
$x\in\Omega$.  Using this notation we re-write our goal in simpler
form: for all $\epsilon>0$, there exist finite constants
$N_0(\epsilon)$ and $C(\epsilon)$ such that
\begin{align}
&\hspace{-1cm}\sum_{m=0}^{R-1}\frac{\nu_{\Omega^+,R}(\eta_z=m)}{\nu_{\Omega,R-m}[c_{x}(\eta_{x})]}\left\{\nu_{\Omega,R-m}[f;AV_{y\in\Omega}c_y(\eta_y)]\right\}^2\notag\\
&\leq \ \  C(\epsilon)N^{1-d}D_{\Omega^+,R}(f)+\epsilon
N^{-d}\nu_{\Omega^+,R}[f;f] \label{line:T2restate}
\end{align}
for all $N\geq N_0$.

The proof of this is split into two cases: that of small density
and that of large density.  To this end, let
$\rho=\frac{R}{|\Omega^+|}$, and fix $\rho_0>0$.
\medskip

\noindent\textbf{Case 1. $\rho\leq\rho_0$.} By the Schwarz
inequality we have that
$$\nu_{\Omega,R-m}[f;AV_{y\in\Omega}c_y(\eta_y)]^2\leq \nu_{\Omega,R-m}[f;f]\nu_{\Omega,R-m}[AV_{y\in\Omega}c_y(\eta_y);AV_{y\in\Omega}c_y(\eta_y)].$$
Using change of measure we may write
\begin{align}
&\frac{1}{\nu_{\Omega,R-m}[c_z(\eta_z)]}\nu_{\Omega,R-m}[AV_{y\in\Omega}c_y(\eta_y);AV_{y\in\Omega}c_y(\eta_y)]\notag\\
&=\nu_{\Omega,R-m-1}[c_z(\eta_z)]-\nu_{\Omega,R-m}[c_z(\eta_z)]+\frac{1}{|\Omega|}\nu_{\Omega,R-m-1}[c_z(\eta_z+1)-c_z(\eta_z)]\label{sg:line5}
\end{align}
We wish to bound the term in \reff{sg:line5} by
$\frac{\epsilon}{|\Omega|}$ for all sufficiently large $|\Omega|$.
For $\rho\geq \frac{A}{|\Omega|}$ this follows from Proposition
\ref{submom7}.  Otherwise the number of particles is bounded and
we have that the Poisson limit theorem holds.  Here again we
obtain the desired bounds from Lemma \ref{poissonLLT-2}.

\medskip

\noindent\textbf{Case 2. $\rho>\rho_0$.}  This is the more
involved case of the two, and it requires a ``two-block" argument.
That is, we write $\Omega$ as a union of smaller cubes $B_1 \cup
\cdots \cup B_K$, where for simplicity we assume that each cube is
exactly of size $l^d$.  As in the previous section we shall pick
$l$ to be a fixed quantity, however, sufficiently large so that
certain estimates hold.  We write as $R_j$ the number of particles
on cube $B_j$, $j=1, \ldots, K$. To simplify notation we also
write $\tilde R$ instead of $R-m$. By the Schwarz inequality we
write
\begin{eqnarray*}
\nu_{\Omega,\tilde R}[f;AV_{y\in\Omega}c_y(\eta_y)]^2&\leq&
2\nu_{\Omega,\tilde
R}\left[f;AV_{y\in\Omega}\{c_y(\eta_y)-\nu_{B_j,R_j}[c_y(\eta_y)]\mathbb
I(y\in B_j) \}\right]^2\\
&&\hspace{4cm}+\ \ \ 2\nu_{\Omega, \tilde R}\left[f;
AV_j\nu_{B_j,R_j}[c_y(\eta_y)] \right]^2
\end{eqnarray*}
We first handle the first term on the right hand side. We may
write this as
\begin{align}
&\hspace{-1cm}\left|\frac{1}{|\Omega|}\sum_{j=1}^K |B_k|\
\nu_{\Omega,\tilde
R} \left[\nu_{B_j,R_j}\left[f;AV_{y\in B_j}c_y(\eta_y)\right]\right]\right|\notag\\
&\leq \ \ \frac{\alpha}{2|\Omega|}\sum_k |B_k|\ \nu_{\Omega,\tilde
R}
\lf[\nu_{B_j,R_j}\lf[f;f\rg]\rg]\notag\\
&\ \ \ \ +\frac{1}{2\alpha|\Omega|}\sum_{j=1}^K |B_k|\
\nu_{\Omega,\tilde R} \lf[\nu_{B_j,R_j}\lf[AV_{y\in
B_j}c_y(\eta_y);AV_{y\in B_j}c_y(\eta_y)\rg]\rg],\label{sg:line6}
\end{align}
for any strictly positive $\alpha$. By the induction assumption we
have
\begin{align*}
\sum_k |B_k|\ \nu_{\Omega,\tilde R}
\lf[\nu_{B_j,R_j}\lf[f;f\rg]\rg]\leq l^d \omega(l+1)
D_{\Omega,\tilde R}(f).
\end{align*}
On the other hand, by the Schwarz inequality we have for the
second term
\begin{align*}
\frac{1}{|\Omega|}\sum_{j=1}^K |B_k|\nu_{\Omega,\tilde R}
\lf[\nu_{B_j,R_j}\lf[AV_{y\in B_j}c_y(\eta_y);AV_{y\in
B_j}c_y(\eta_y)\rg]\rg]\leq AV_{y\in\Omega}\nu_{\Omega,\tilde
R}\lf[c_y(\eta_y);c_y(\eta_y)\rg],
\end{align*}
which, by Proposition \ref{submom}, is bounded above by $C\rho$,
for some constant $C$, and sufficiently large $|\Omega|$.
Plugging these bounds into \reff{sg:line6} and optimising in
$\alpha$ we obtain that
\begin{align*}
&2\nu_{\Omega,\tilde
R}\left[f;AV_{y\in\Omega}\{c_y(\eta_y)-\nu_{B_j,R_j}[c_y(\eta_y)]\mathbb
I(y\in B_j) \}\right]^2\\
&\leq C l^d \omega(l+1)\frac{1}{|\Omega|} D_{\Omega,\tilde
R}(f)\nu_{\Omega,\tilde R}[c_x(\eta_x)],
\end{align*}
where we have also used the fact that there exists a positive
constant such that $\rho\leq C\nu_{\Omega,\tilde R}[c_x(\eta_x)].$

We next turn our attention to the second term.  Using a similar
argument to that of Section \ref{sec:tight} we write
\begin{align*}
&2\nu_{\Omega, \tilde R}\left[f; AV_j\nu_{B_j,R_j}[c_y(\eta_y)]
\right]^2\\
&=2\nu_{\Omega, \tilde R}\left[f;
\frac{1}{|\Omega|}\sum_{j=1}^K|B_j|\left\{\rule[-0.1cm]{0cm}{0.6cm}\nu_{B_j,R_j}[c_y(\eta_y)]-\varphi(\tilde\rho)-\varphi_j'(\rho_j)\left\{AV_{y\in
B_j }\eta_y-\rho_j\right\}\right\} \right]^2,
\end{align*}
where $\tilde\rho=\tilde R/|\Omega|$,
$\rho_j=\mu_{\Omega,\varphi(\tilde\rho)}[AV_{y\in B_j}\eta_y]$ and
$\varphi_j(\rho_j)=\mu_{B_j,\varphi(\tilde\rho)}[AV_{y\in
B_j}c_y(\eta_y)]$.  Notice that
$\varphi(\tilde\rho)=\varphi_j(\rho_j)$ (however, also note that
$\varphi'(\tilde\rho)$ is not equal to $\varphi_j'(\rho_j)$). We
let \linebreak $m_j=AV_{y\in B_j}\eta_y=R_j/l$ and set
$$F_j(m_j)=\nu_{B_j,R_j}[c_y(\eta_y)]-\varphi(\tilde\rho)-\varphi_j'(\rho_j)\left\{m_j-\rho_j\right\}.$$
 With this notation we bound the last line above
using the Schwarz inequality by
\begin{align*}
2\nu_{\Omega,\tilde R}[f;f]\nu_{\Omega,\tilde
R}\left[\left(\frac{1}{|\Omega|}\sum_{j=1}^K|B_j|F_j(m_j)\right)^2\right].
\end{align*}
We then write
\begin{align}
&\hspace{-2cm}\nu_{\Omega,\tilde
R}\left[\left(\frac{1}{|\Omega|}\sum_{j=1}^K|B_j|F_j(m_j)\right)^2\right]\notag\\
&=\frac{l^{2d}}{|\Omega|^2}\sum_{j=1}^K\nu_{\Omega,\tilde
R}[F_j^2(m_j)]+\frac{l^{2d}}{|\Omega|^2}\sum_{j\neq
i}\nu_{\Omega,\tilde R}[F_j(m_j)F_i(m_i)].\label{sg:line7}
\end{align}
We next use the second part of Proposition \ref{submom2} to switch
to the grand canonical measure.  This will allow us to take
advantage of the Taylor series expansion we have set up.  We do
not use Proposition \ref{submom} here as we will take advantage of
the freedom of making $l$ large.  Hence, we require the condition
$\tilde \rho>\rho_0$ to make the argument. Thus,
\begin{eqnarray*}
\nu_{\Omega,\tilde R}[F_j^2(m_j)]&\leq&
\mu_{\varphi(\tilde\rho)}[F_j^2(m_j)]+E_0(\rho_0)\frac{l^d}{|\Omega|}\{\mu_{\varphi(\tilde\rho)}[F_j^4(m_j)]\}^{1/2},
\end{eqnarray*}
and similarly
\begin{eqnarray*}
\hspace{-0.5cm}\nu_{\Omega,\tilde R}[F_j(m_j)F_i(m_i)]&\leq&
\mu_{\varphi(\tilde\rho)}[F_j(m_i)]\mu_{\varphi(\tilde\rho)}[F_i(m_i)]\\
&&\hspace{2cm}+\ \
E_0(\rho_0)\frac{l^d}{|\Omega|}\{\mu_{\varphi(\tilde\rho)}[F_j^2(m_j)]\mu_{\varphi(\tilde\rho)}[F_i^2(m_i)]\}^{1/2}\\
&=&E_0(\rho_0)\frac{l^d}{|\Omega|}\{\mu_{\varphi(\tilde\rho)}[F_j^2(m_j)]\mu_{\varphi(\tilde\rho)}[F_i^2(m_i)]\}^{1/2}.
\end{eqnarray*}
Applying the first part of Corollary \ref{submom2} to
$\nu_{B_j,R_j}[c_y(\eta_y)]$ , and using Propositions
\ref{sigmabounds}, \ref{hatbounds}, and Corollary \ref{phideriv}
to bound the resulting moments, we obtain that
\begin{eqnarray*}
\mu_{\varphi(\tilde\rho)}[F_j^4]&\leq&
\mu_{\varphi(\tilde\rho)}\left[\left\{\frac{C(\rho_0)}{l^d}\sqrt{1+\tilde\rho}+\tilde
c
(m_j-\rho_j) \right\}^4\right]\\
&\leq& C'(\rho_0) \frac{1}{l^{\,2d}}\tilde\rho^2,
\end{eqnarray*}
for $l$ and $|\Omega|$ sufficiently large.  For the quadratic term
we have
\begin{eqnarray*}
2^{-1}\mu_{\varphi(\tilde\rho)}[F_j^2]&\leq&\mu_{\varphi(\tilde\rho)}\left[\nu_{B_j,R_j}\left[c_y(\eta_y)-\varphi(m_j)\rule[-0.1cm]{0cm}{0.6cm}\right]^2\right]\\
&&+\mu_{\varphi(\tilde\rho)}\left[\left\{\varphi(m_j)-\varphi_k(\rho_k)-\varphi_k'(\rho_k)(m_k-\rho_k)\rule[-0.1cm]{0cm}{0.6cm}\right\}^2\right]
\end{eqnarray*}
By Corollaries \ref{submom2} and \ref{phideriv} we bound the first
term above by  $C \frac{1}{l^{\,2d}}(1+~\tilde\rho)$, for some
constant $C$. The second term may be bounded by
$C\frac{1}{l^{\,2d}}\tilde\rho$ using \reff{useagain1} and
\reff{useagain2} (noting that those particular arguments do not
depend on the dimension). We now put all of the above work
together to obtain the bound in \reff{sg:line7}
\begin{align}
\nu_{\Omega,\tilde
R}\left[\left(\frac{1}{|\Omega|}\sum_{j=1}^K|B_j|F_j(m_j)\right)^2\right]
\leq C(\rho_0)\frac{1}{|\Omega|}\frac{1}{l^d}\tilde\rho.
\end{align}
We may now select $l=\epsilon^{-d}$.  This completes the proof of
\reff{line:T2restate}.

Notice that in the above arguments it is only Proposition \ref{T1}
which is sensitive to the geometry of the problem induced by
change in dimension.


\begin{thebibliography}{9}

\bibitem[A]{a} Andjel E.; Invariant measures for the zero range
processes.  Annals of Probability {\bf 10}, 525--547, (1982).

\bibitem[B]{B} Billinglsey, P.; \textit{Convergence of Probability
Measures.} John Wiley \& Sons, Inc., (1968)

\bibitem[CMR]{cmr} Canrini N., Martinelli F., Roberto C.; The logarithmic Sobolev constant of Kawasaki dynamics under a
mixing condition revisited.  Ann. Inst. H. Poincar\'{e} Probab.
Statist. {\bf 38}, No. 4, 385--234, (2002)


\bibitem[CP]{cp} Caputo P., Posta G.; Entropy Dissipation Estimates in
a Zero-Range Dynamics.  Preprint.

\bibitem[CY]{yau92} Chang C.C., Yau H.T.; Fluctuations of one dimensional {G}insberg-{L}andau models in
nonequilibirium. Communication in Mathematical Physics {\bf 69},
No. 145, 209--234, (1992).

\bibitem[DGS]{dgs} Davies E.B., Gross L., Simon B.; Hypercontractivity: A bibliographical review. In: Ideas
and Methods of Mathematics and Physics. In Memoriam of Raphael
Hoegh-Krohn. S. Albeverio, J.E. Fenstand, H. Holden, T. Lindstrom
(eds.) Cambridge: Cambridge University Press, (1992)

\bibitem[DPP1]{logsob1} Dai Pra P., Posta G.; Logarithmic Sobolev
Inequality for Zero-Range Dynamics: Independence of the number of
particles.  Electronic Journal of Probability {\bf 10}, 525--576,
(2005).


\bibitem[DPP2]{logsob2} Dai Pra P., Posta G.; Logarithmic Sobolev
Inequality for Zero-Range Dynamics.  To Appear in:  Annals of
Probability .


\bibitem[DS]{persi} Diaconis P., Saloff-Coste L.; Logarithmic
Sobolev inequalities for finite Markov chains.  Annals of Applied
Probability {\bf 6}, No. 2, 696--759, (1996).

\bibitem[G]{gross} Gross L., Logarithmic Sobolev Inequalities, Am.
J. Math. {\bf 97}, 1061--1083 (1976)

\bibitem[GK]{gk} Gnedenko, B.V., Kolmogorov A.N.; \textit{Limit Distributions
for Sums of Independent Random Variables.}  Addison-Wesley
Publishing Company, (1962).

\bibitem[KL]{kl} Kipnis C., Landim C.; \textit{Scaling
Limits of Interacting Particle Systems.} Springer Verlag, Berlin,
(1999).

\bibitem[L]{lig} Ligget, T., \textit{Interacting Particle
Systems.} Springer-Verlag, (1985).

\bibitem[LSV]{lsv} Landim C., Sethuraman S., Varadhan S. R. S.;
Spectral Gap for Zero-Range Dynamics. Annals of Probability {\bf
24}, No.4, 1871--1902, (1996).

\bibitem[LY]{ly} Lu S.L., Yau H.T.; Spectral gap and logarithmic
{S}obolev inequality in {K}awasaki and {G}lauber dynamics.
Communications in Mathematical Physics, {\bf 156}, 399--433,
(1993).

\bibitem[J]{me} Jankowski, H.; Nonequilibrium density fluctuations
for the zero range process with colour. Preprint (2006).

\bibitem[M1]{m1} Miclo L.; An Example of Application of Discrete
Hardy's Inequalities.  Markov Processes Relat. Fields {\bf 5},
319--330, (1999).

\bibitem[M2]{m2} Miclo L.; Relations entre isp\'{e}rim\'{e}trie et trou
spectral pour les cha\^{i}nes de Markov finies.  Probability
Theory and Related Fields {\bf 114}, No. 4, 431--485 ,(1998).

\bibitem[P]{pet} Petrov V.V.; \textit{Sums of Independent Random Variables.}
Springer-Verlag, New York, 1975.

\bibitem[Q]{q} Quastel J.; Diffusion of color in the simple exclusion
process.  Communications in Pure and Applied Mathematics, {\bf45},
623--679, (1992).

\bibitem[R]{Rot} Rothaus O.S.; Analytic inequalities,
isoperimetric inequalities, and logarithmic Sobolev inequalities.
Journal of Functional Analysis, {\bf 64}, 296--313, (1985).


\end{thebibliography}
\end{document}